\documentclass[11pt,reqno,letterpaper]{amsart}

\usepackage{mathtools, amsthm, verbatim, amscd, amssymb, setspace, enumitem, exscale, color, slashed, arydshln, accents, scalerel}
\usepackage{esint}
\usepackage[all,pdf]{xy}
\usepackage{stackengine}
\usepackage{hyperref}

\newcommand{\CC}{\mathbb{C}}

\newcommand{\RR}{\mathbb{R}}

\newcommand{\im}{\operatorname{Im}}
\newcommand{\psh}{\operatorname{Psh}}
\newcommand{\C}{\mathcal{C}}

\newcommand{\F}{\mathcal{F}}
\newcommand{\HH}{\mathcal{H}}
\newcommand{\W}{\mathcal{W}}

\newcommand{\LL}{\mathcal{L}}
\newcommand{\I}{\mathcal{I}}

\newcommand{\re}{\operatorname{Re}}
\newcommand{\bff}{\mathfrak{b}}

\newcommand{\ZZ}{\mathbb{Z}}

\newcommand{\vol}{\operatorname{vol}}

\newcommand{\diam}{\operatorname{diam}}

\newcommand{\B}{\mathcal{B}}
\newcommand{\D}{\mathcal{D}}
\newcommand{\BMO}{\mathrm{BMO}}

\newcommand{\Hilb}{\operatorname{Hil}}

\newtheorem{thm}{Theorem}[section]
\newtheorem{pro}[thm]{Proposition}
\newtheorem{cor}[thm]{Corollary}
\newtheorem{lem}[thm]{Lemma}

\newtheorem{mtheorem}{Theorem}

\theoremstyle{definition}
\newtheorem{define}[thm]{Definition}
\newtheorem{rem}[thm]{Remark}

\newtheorem*{ac}{Acknowledgements}
\newtheorem*{aideclaration}{Declaration of AI use}

\numberwithin{equation}{section}

\title[On the Local Regularity of HCMA Equations]%
      {On the Local Regularity of Homogeneous Complex Monge-Ampère Equations}

\author{Qi Yao}
\address{IASM, Zhejiang University, Hangzhou, Zhejiang, China}
\email{qiyao11@zju.edu.cn}
\date{}
\begin{document}
\raggedbottom

\begin{abstract}
We establish a local higher regularity theorem for the homogeneous complex Monge--Amp\`ere equation on $X\times D$, where $X$ is a complete, possibly noncompact K\"ahler manifold and $D\subset\mathbb C$ is the unit disc. For a bounded weak solution given by the Perron envelope, we prove higher regularity on a smaller spatial domain when the smooth boundary data are locally small. The geometric hypotheses require uniformly convex local K\"ahler potentials and a local barrier threshold sufficiently large relative to the global $L^\infty$ bound of the solution. We also obtain quantitative local H\"older estimates. 

The proof combines a perturbative construction of local holomorphic disc foliations with a construction of global subsolutions using convexity and barriers. A key analytic difficulty is the loss of H\"older regularity in the construction of holomorphic disc families, which we address using BMO estimates and Nash--Moser--Zehnder theorems. Comparison along the discs identifies the local smooth solution by holomorphic discs foliation with the global weak solution.

We apply the theorem and its proof to compact K\"ahler manifolds and to ALE, ALF, and asymptotically complex hyperbolic (ACH) ends. Under the stated geometric and decay assumptions, we obtain higher regularity of each prescribed order sufficiently far out, together with weighted H\"older estimates preserving the boundary decay.

\end{abstract}

\maketitle

\section{Introduction}

Initiated by the works of Mabuchi \cite{mabuchi1987some}, Semmes \cite{semmes1992complex}, and Donaldson \cite{donaldson1999symmetric}, the HCMA equation has emerged as a central object in the study of canonical K\"ahler metrics, particularly in relation to the constant scalar curvature K\"ahler (cscK) problem on compact K\"ahler manifolds.
Fixing a reference K\"ahler form $\omega$ on a compact K\"ahler manifold, the space of all K\"ahler metrics in $[\omega]$ can be identified with the space of K\"ahler potentials,
\begin{align*}
\HH(X, \omega) = \{\varphi \in \C^\infty (X); \ \omega + i \partial \overline{\partial} \varphi > 0\},
\end{align*}
modulo constants. Given two potential functions in $\HH(X, \omega)$, $\psi_0$, $\psi_1$,
we connect $\psi_0$ and $\psi_1$ with a path $\varphi(t): [0,1] \rightarrow \HH(X,\omega)$ with $\varphi(0) = \psi_0$, $\varphi(1) = \psi_1$. We say $\varphi(t)$ is a geodesic in $\HH(X,\omega)$, if 
\begin{align} \label{geoequintro}
    \ddot{\varphi} - \frac{1}{2} |\nabla_{\omega_{\varphi(t)}} \dot{\varphi}|_{\omega_{\varphi(t)}}^2 = 0.
\end{align}
The above equation \eqref{geoequintro} is called the geodesic equation in the space of K\"ahler potentials.
As observed by Donaldson \cite{donaldson1999symmetric} and Semmes \cite{semmes1992complex}, the geodesic equation is equivalent to the HCMA equation in the product space $X \times \Sigma$, where $\Sigma \cong [0,1]\times S^1$ can be embedded as an annulus in $\CC$. 
Notice that any path in $\HH(X,\omega)$ can be viewed as a function defined in $X \times \Sigma$, $\Phi (\cdot, t, e^{is}) = \varphi(t)$. Let $\Omega =\pi^* \omega $, where $\pi$ is the natural projection from $X\times \Sigma$ to $X$. Then the equation (\ref{geoequintro}) can be rewritten on $X \times \Sigma$ as follows,
\begin{align*}
    (&\Omega +i \partial \overline{\partial} \Phi)^{n+1} = 0, \\
    &\Omega +i \partial \overline{\partial} \Phi \geq 0, \\
    &\Phi|_{X \times \partial\Sigma } =  \psi_{0,1}. 
\end{align*}

In the cases of compact K\"ahler manifolds, Chen \cite{chen2000space} showed that for any $\psi_0$, $\psi_1 \in \HH$, the geodesic equation has a unique solution up to $dd^c$-regularity. Blocki \cite{blocki2009gradient} and He \cite{he2015space} built up direct calculations to prove the gradient estimate. The full $\C^{1,1}$ estimate was proved by Chu-Tosatti-Weinkove in \cite{chu2017regularity}. In the other direction, Lempert-Vivas \cite{lempert2013geodesics} and Darvas-Lempert \cite{darvas2012weak} constructed counterexamples to assert that $\partial \overline{\partial} \varphi$ is not continuous in general, hence the $\C^{1,1}$ regularity is optimal in general.

While the counterexamples cited above preclude higher regularity in the general case, global smoothness can still be established within perturbative regimes (cf. \cite{Lempert1983, donaldson2002holo, CHEN2020108603, hu2024metric}). However, these results invariably rely on global smallness assumptions on the boundary data or slight perturbations of the trivial solution on the boundaries.

A more geometrically significant question arises in the context of complete (possibly non-compact) Kähler manifolds: if the boundary data satisfies a smallness condition only on $U \times \partial D$, where $U$ is an open domain in $X$, does the solution exhibit higher regularity within $V \times D$ for some $V \Subset U$? This question is particularly relevant for the study of open Kähler manifolds (such as ALE, AC, ALF, $\ldots$), where the metric is asymptotic to a fixed model. In these settings, the smallness assumption emerges naturally near spatial infinity. Consequently, the local regularity theory established in this work serves as the essential bridge for understanding the asymptotic behavior of Kähler metrics along $\C^{1,1}$ geodesics.

\subsection{Setup and Main Theorem}

Let $(X, g)$ be a complete (possibly noncompact) K\"ahler manifold of complex dimension $n$ and let $D \subseteq \CC$ be the unit disc. In this paper, we consider the HCMA equation on the product space $X \times D$. Let $ \omega(\cdot, \cdot) = g(J\cdot, \cdot)$ be the reference K\"ahler form on $X$, and let $\Omega = \pi^* \omega$ be the pullback of $\omega$ to $X\times D$, where $\pi: X \times D \rightarrow X$ is the natural projection. We denote by $\partial$ and $\overline{\partial}$ the complex differential operators on the product space $X \times D$. We consider the following Dirichlet-type problem for the HCMA equation:
\begin{align}\label{introHCMA}
    \begin{cases}
        (\Omega + i\partial \overline{\partial} \Phi)^{n+1} = 0, \hspace{1cm} \text{ in } X \times D;\\
        \Omega + i\partial \overline{\partial} \Phi \geq 0, \hspace{1.9cm} \text{ in } X \times D;\\
        \Phi(x, \tau) = \Psi(x, \tau), \hspace{1.35cm} \text{ on } X \times \partial D,
    \end{cases}
\end{align}
where $\Psi$ is called the \textit{boundary data} of the above HCMA equation and satisfies
\begin{flalign}\label{cond:boundary_data}
    \begin{aligned}
      \quad  &\bullet \ \Psi \in \mathcal{C}^\infty (X \times \partial D); \\
       \quad &\bullet \ \psi_\tau (\cdot) = \Psi (\cdot, \tau) \in \mathcal{H}(X, \omega), \quad \text{for all } \tau \in \partial D.
    \end{aligned} &&
\end{flalign}

To establish local higher regularity for \eqref{introHCMA}, we impose the following assumptions on the background geometry $(X, \omega)$, the boundary data $\Psi$, and the global weak solution \(\Phi\).

We first formulate the global solution in the framework of pluripotential theory. Let
\(\psh_{\Omega}(X\times D)\) denote the space of \(\Omega\)-plurisubharmonic ($\Omega$-psh) functions on \(X\times D\). Given boundary data \(\Psi\), define the associated class of bounded subsolutions by:
\begin{align*}
    \B_{\Omega, \Psi} = \{u \in &\psh_{\Omega}(X \times D);\  u \text{ is bounded, and } 
    \\
    &\limsup_{(x',\tau')\rightarrow (x, \tau)} u(x', \tau') \leq \Psi (x, \tau), \text{ for } (x', \tau') \in X \times  D, \ (x, \tau) \in X \times \partial D   \}.
\end{align*}
We say that $\Phi$ is a \textit{weak solution} to \eqref{introHCMA} if $\Phi$ is continuous on $X \times D$ and satisfies $(\Omega + i\partial \overline{\partial} \Phi)^{n+1} = 0$ in the sense of currents; specifically, in the sense of Bedford-Taylor.
Throughout this paper, we adopt the \emph{Perron Envelope Assumption}: we assume that \(\Phi\) is given by the Perron envelope associated with \((\Omega,\Psi)\), namely
\begin{align}\label{cond:perron}
\Phi (x, \tau) = \Phi_{\Omega, \Psi} (x, \tau) := \sup_{u \in \mathcal{B}_{\Omega, \Psi}} u(x, \tau), \qquad \text{ for } (x, \tau) \in X \times D.
\end{align}

Furthermore, we assume that the global solution $\Phi$ given as above satisfies the \emph{${L^\infty}$ boundedness}:
\begin{align}\label{cond:Linftybdd}
    \| \Phi\| _{L^\infty (X \times D)} \le M.
\end{align}

Let
\(
V\Subset U\subset X
\)
be open domains. We assume that \(\Psi\) is small on \(U\times \partial D\) in the {\(C^{k+3,\alpha}\)}-norm:
\begin{align}\label{cond:localsmall}
\|\Psi\|_{{C^{k+3,\alpha}(U\times \partial D)}}\leq \varepsilon.
\end{align}
No smallness or other constraint is imposed on \(\Psi\) outside \(U\times \partial D\).

Finally, we impose a local convexity condition on the reference Kähler metric \(\omega\) on the pair \(V\Subset U\). We assume that, for every \(x\in V\), there exists a pair of concentric holomorphic balls
\(
V_x\Subset U_x\subset U
\)
centered at \(x\), together with a holomorphic coordinate map
\(
b_x:U_x\longrightarrow B_{R_x}\subset \mathbb C^n,
~ b_x(x)=0,
\)
such that \(b_x(V_x)=B_{r_x}\), and
\[
R_x-r_x\geq R_0
\]
for some uniform constant \(R_0>0\). We further assume that, in these coordinates, \(\omega\) admits a uniformly convex local potential: there exists a smooth function \(\rho\) on \(B_{R_x}\) such that
\begin{align}\label{cond:localcx}
(b_x^{-1})^*\omega
=
i\partial\overline\partial \rho,
\qquad
\nabla_0^2\rho\geq \lambda I
\quad \text{on } B_{R_x},
\end{align}
where \(\nabla^2_0\rho\) denotes the standard Euclidean Hessian on \(B_{R_x}\). The constants \(R_0\) and \(\lambda\) are independent of \(x\in V\).

The relevant scale-invariant quantity is \(\lambda R_0^2\), which we call the \emph{local barrier threshold} of the pair \(V\Subset U\) in \((X,\omega)\). In order to obtain local higher regularity on \(V\), we assume
\begin{align}\label{thresholdcon}
\lambda R_0^2\geq 16M.
\end{align}

We can now state the main local regularity theorem.

\begin{mtheorem}\label{mthmlocalreg}
Let \((X,\omega)\) be a complete, possibly non-compact, Kähler manifold, and let \(\Psi\) be boundary data on \(X\times \partial D\) satisfying \eqref{cond:boundary_data}. Let \(\Phi\) be the corresponding weak solution of the HCMA equation \eqref{introHCMA}. Assume that \(\Phi\) is given by the Perron envelope \eqref{cond:perron} and satisfies the global bound \eqref{cond:Linftybdd}.

Let \(V\Subset U\subset X\) be open domains. Assume that the reference metric \(\omega\) satisfies the local convexity condition \eqref{cond:localcx} on the pair \(V\Subset U\), with uniform constants \(R_0\) and \(\lambda\), and that the corresponding local barrier threshold satisfies \eqref{thresholdcon}.

{Fix \(k\ge1\) and \(0<\beta<\alpha<1\).} Then there exists \(\varepsilon>0\), depending only on {\(n,k,\alpha,\beta,\lambda,R_0,M\)}, {the fixed chart bounds,} and {the uniform potential bounds in \eqref{cond:potentialbds}}, such that if
\[
\|\Psi\|_{{C^{k+3,\alpha}(U\times \partial D)}}\leq \varepsilon,
\]
then \(\Phi\) is locally smooth on \(V\times D\). Moreover, {the} estimate
\begin{align}\label{estilocal}
\|\Phi\|_{C^{k,\beta}(V\times D)}
\leq
C\,\|\Psi\|_{C^{k+1,\alpha}(U\times \partial D)}
\end{align}
holds, where \(C\) is independent of \(\Phi\) and \(\Psi\).
\end{mtheorem}

Theorem~\ref{mthmlocalreg} is a flexible criterion for local higher regularity. In Section~\ref{subsec:geometric_applications}, we discuss applications to compact K\"ahler manifolds and to ALE, ALF, and asymptotically complex hyperbolic (ACH) ends. On a compact manifold, we obtain higher regularity from a finite covering when the boundary data are sufficiently small. On the noncompact ends, under the stated geometric and decay hypotheses, the local barrier thresholds increase towards infinity while the boundary data become locally small. These two properties imply higher regularity sufficiently far out without a smallness restriction on the global bound $M$.  For an ALE or ALF manifold, let $X_\infty$ denote the end, let $r\ge1$ be its radius function extended smoothly over the compact part, and write $X_l=\{x\in X_\infty:r(x)>l\}$. Under the stated geometric assumptions, we can choose local domains with arbitrarily large barrier thresholds as $l\to\infty$, while decay of the boundary data ensures local smallness. These properties yield higher regularity sufficiently far out, without requiring the global bound $M$ to be small. We express the decay using weighted H\"older norms $\|\cdot\|_{-\gamma;k,\alpha}$, obtained by taking the supremum of $r(p)^\gamma$ times the local H\"older norms. In the ALE case, these local norms use rescaled end coordinates, as in \eqref{eq:ale_application_weighted_norm}. The similar techniques apply to the ALF case, with the rescaled coordinates adapted to the ALF asymptotics. In the ACH case, we use the hyperbolic distance barrier in place of the quadratic barrier in Theorem~\ref{mthmlocalreg}; see Lemma \ref{lem:ach_distance_barrier}. 

\subsection{Strategy of the Proof}

The proof of Theorem~\ref{mthmlocalreg} has two main steps. We first construct a local family of holomorphic discs from the small boundary perturbation. We then use the convexity and barrier assumptions to construct a global $\Omega$-psh subsolution which agrees with the local disc solution on the smaller domain. The Perron envelope construction and comparison along the discs identify this local solution with $\Phi$. 

\smallskip\noindent\textit{Step 1. Construction of local holomorphic discs.}
We recall the notation for the holomorphic discs families by Semmes \cite{semmes1992complex} and Donaldson \cite{donaldson2002holo}. Let $N\subset N'\subset X$ be an inner parameter domain and an outer domain, with a positive spatial margin in the chosen charts. On each holomorphic chart, choose a potential $\rho$ for $\omega$ and write $E$ for the local holomorphic cotangent bundle. These local bundles form an affine holomorphic fiber bundle $\pi:\mathcal W\to X$: on an overlap, the fibre coordinates satisfy $\xi_j=\xi_i+\partial(\rho_j-\rho_i)$, with the usual cotangent coordinate change understood. In this subsection, $\pi$ denotes this bundle projection.

For each $\tau\in\partial D$, write $\psi_\tau=\Psi(\cdot,\tau)$. The local graphs of $\partial(\rho+\psi_\tau)$ satisfy the transition rules of $\mathcal W$ and therefore define a submanifold $\Lambda_\tau\subset\mathcal W$. These submanifolds prescribe the boundary conditions for the holomorphic discs. When $\Psi=0$, the boundary condition reduces to $\Lambda_\rho=\operatorname{graph}\partial\rho$, and the constant discs form the trivial family
\[
 G_0(w,\tau)=(w,\partial\rho(w)),\qquad w\in N.
\]
For small $\Psi$, we seek a family $G:N\times\overline D\to\mathcal W$ close to $G_0$, holomorphic in $\tau\in D$, satisfying
\[
 G(w,\tau)\in\Lambda_\tau\quad(\tau\in\partial D),
 \qquad \pi\circ G(w,-i)=w.
\]
The projected family defines a foliation by holomorphic discs over a smaller spatial domain, as described in Section~\ref{subsec:hcma_disc_families}. In the estimates below, we write $\widetilde\Psi=\rho+\Psi$ and $\widetilde\Psi_0=\rho$ for the perturbed and reference local potentials, respectively.

The linearized boundary problem is a family of Riemann--Hilbert problems with parameter $w$. The main analytic difficulty is that the fibrewise Hilbert transform need not preserve the same joint H\"older exponent in the disc and parameter variables. We obtain estimates with an arbitrarily small loss of H\"older exponent using BMO estimates, and solve the nonlinear problem by the Nash--Moser--Zehnder iteration. Local uniqueness near $G_0$ allows the families constructed in overlapping charts to be patched together.

\begin{mtheorem}[{see Theorem \ref{thmexistbypertb}}] \label{mthmexistholodiskfoli}
    Let $N\Subset N'\subset X$ satisfy the uniform local geometric bounds described above, and let $G_0$ be the trivial HCMA disc family. Fix $k\ge1$ and $0<\beta<\alpha<1$. There exists $\varepsilon>0$ such that, if $\Psi\in C^\infty(N'\times\partial D,\mathbb R)$ satisfies
\[
 \|\partial\Psi\|_{k+2,\alpha;N'\times\partial D}\le\varepsilon,
\]
then there exists an HCMA disc family $G:N\times D\to\mathcal W$ adapted to $\Psi$, satisfying
\[
 \|G-G_0\|_{k,\beta;N\times D}
 \le C\|\partial\Psi\|_{k,\alpha;N'\times\partial D}.
\]
The constants $\varepsilon$ and $C$ depend only on $k,\alpha,\beta$, the chart margins, and the uniform local geometric bounds.
\end{mtheorem}

\smallskip\noindent\textit{Step 2. Comparison with the global weak solution.}
The disc family determines a local solution of HCMA. To prove Theorem~\ref{mthmlocalreg}, we must show that this solution agrees with the given global solution $\Phi$. This is where the convexity and barrier threshold enter the argument.

For each leaf labelled by $x$, we construct a real pluriharmonic function $L_x$ on the local product chart. Along the leaf, it is the harmonic extension of the boundary values $\rho+\psi_\tau$; in the spatial coordinates, it is affine with differential given by the fibre component of the disc. On $\partial D$, it is therefore an affine support of $\rho+\psi_\tau$. Convexity implies $L_x-\rho\le\Psi$ on the boundary, while the barrier estimate places $L_x-\rho$ below $-M$ near the outer spatial boundary. The latter bound extends into $D$ by the maximum principle.

Taking the local supremum over the leaf labels and the maximum with $-M$, we obtain an admissible subsolution which extends across the outer spatial boundary. The compatibility of the local disc families allows these constructions to be patched into a global function $F\in\mathcal B_{\Omega,\Psi}$. Since $\Phi$ is the Perron envelope, $F\le\Phi$. On each leaf through the smaller domain, {$L_x$} is harmonic and $\Phi+\rho$ is subharmonic, with the same boundary values. The maximum principle implies {$\Phi\le L_x-\rho\le F$} there. Hence $F=\Phi$ on the smaller domain. Section~\ref{secconstructF} carries out this comparison, and the local H\"older estimates follow by solving the Dirichlet problem along the discs.

The following statement establishes the existence of a global subsolution that is smooth on the smaller domain and agrees there with the weak solution given by the Perron envelope.

\begin{mtheorem}[{see Theorem \ref{thmconstpshsolu}}] \label{mthmconstrucsubsol}
Let $(X,\omega)$ be a complete K\"ahler manifold, and let $\Psi$ satisfy \eqref{cond:boundary_data}. Let $\Phi$ be the weak solution of \eqref{introHCMA}, given by the Perron envelope \eqref{cond:perron} and satisfying the global bound \eqref{cond:Linftybdd}.

{Fix $k\ge1$ and $0<\beta<\alpha<1$.} Suppose that $N\Subset N'\subset X$ satisfies the local convexity condition \eqref{cond:localcx}, the uniform potential bounds \eqref{cond:potentialbds}, and the barrier condition \eqref{thresholdcon}. Assume also that $\Psi$ satisfies the local smallness condition \eqref{cond:localsmall} on $N'\times\partial D$, with $\varepsilon$ sufficiently small as in Theorem~\ref{mthmlocalreg}.

Then there exists a bounded $\Omega$-plurisubharmonic subsolution
$F\in\mathcal B_{\Omega,\Psi}$ that is continuous on $X\times D$,
smooth on $N\times D$, and satisfies
\[
 F=\Phi\qquad\text{on }N\times D.
\]
\end{mtheorem}

The paper is organized as follows. Section~\ref{sec:notations} introduces the coordinate and norm conventions. Section~3 describes HCMA disc families and their relation to local solutions, and Section~4 establishes the linear estimates and the loss of H\"older regularity. Section~\ref{secholodiscfoliexist} proves existence, local uniqueness, and patching of the disc families. Section~\ref{secmaxrkthm} constructs the global subsolution and establishes the Perron comparison. Section~7 completes the proof of Theorem~\ref{mthmlocalreg} and discusses the geometric applications. Appendix~\ref{app:campanato} contains the Campanato estimates, and Appendix~\ref{app:ach_barriers} proves the barrier and distance estimates used for ACH ends.

\begin{ac}
This paper was largely completed during the author's stay at Stony Brook University and finalized at IASM. The author is grateful to Hans-Joachim Hein, Xiuxiong Chen, Claude LeBrun, and Song Sun for their help and support.
\end{ac}

\begin{aideclaration}
The mathematical ideas, results, and proofs were developed independently by the author, who also wrote the original manuscript without AI assistance. AI tools were used only during the final revision to assist with editing and checking for possible errors. The author takes full responsibility for the final manuscript.
\end{aideclaration}

\section{Notations and Norms} \label{sec:notations}

Let $(X, \omega)$ be a complete (possibly noncompact) K\"ahler manifold of complex dimension $n$. We consider a large open domain $U \subset X$ within which we shall establish our estimates. To formalize the interior regularity argument, we introduce a pair of open domains $V \subset U \subset X$. 
We fix a locally finite open covering $\{V_i\}_{i \in \mathcal{I}}$ of the smaller domain $V$ by holomorphic coordinate charts. Associated with each index $i \in \mathcal{I}$, we assume there exists a larger holomorphic chart $U_i$ with $V_i \subset U_i \Subset U$.

For each chart, we associate the local data $(V_i,U_i, b_i, \rho_i)$, where $b_i: U_i \to \mathbb{C}^n$ is a biholomorphism onto its image and $\rho_i$ is a local potential for the K\"ahler metric $\omega$ (i.e., $\omega = i\partial\bar{\partial}\rho_i$). We assume the following uniform geometric bounds for this atlas:
\begin{itemize}
\item \textit{Uniform Size:} There exists a uniform constant $R_0 > 0$ such that for each $i \in \mathcal{I}$, the chart map $b_i: U_i \to \mathbb{C}^n$ is a biholomorphism onto a domain containing the ball $B_{R}(0)$, and maps $V_i$ onto $B_r (0)$ with $R - r \ge R_0$. 
\item \textit{Uniform Convexity:} The local potentials are strictly convex with a uniform lower bound on the Hessian. Let $\nabla_0$ denote the Levi-Civita connection of the Euclidean metric. There exists a constant $\lambda > 0$ such that for all $i \in \mathcal{I}$ and $z \in b_i(U_i)$,
\begin{align} \label{est:convexity}\nabla_0^2 \rho_i(z) \ge \lambda I_{2n},\end{align}
where $I_{2n}$ denotes the identity matrix.
\item \textit{Potential Bounds:} The second derivatives of the potentials are uniformly bounded in the $C^{k+2, \alpha}$ topology. There exists a constant $C > 0$ such that for all $i \in \mathcal{I}$,
\begin{align} \label{cond:potentialbds}
\|\nabla_0^2\rho_i\|_{C^{k+2, \alpha}(U_i)} \le C.
\end{align}
\end{itemize}

In the norm definitions below, $N$ denotes a domain covered by the chosen charts, and $N_i:=N\cap U_i$. We define global H\"older norms for scalar functions and sections of vector bundles by pulling back to these local charts. First, recall the standard H\"older norm on a domain $ U \subseteq \CC^n $. For a function $u$ on $ U $ and $k \in \mathbb{N}, \alpha \in (0,1)$, we define
\begin{align*}
\| u\| _{{k, \alpha}; U} &\coloneqq \sum_{|\beta| \le k} \sup_{x \in U} |D^\beta u(x)| + \sum_{|\beta| = k} \sup_{x, y \in U, x \neq y} \frac{|D^\beta u(x) - D^\beta u(y)|}{|x - y|^\alpha},
\end{align*}
where $\beta$ denotes a multi-index. For a function $f: N \to \mathbb{C}$, we define its local norm on the chart $N_i$ by regarding it as a function on the Euclidean domain $b_i(N_i)$:
\begin{equation*}
\| f\| _{k, \alpha; N_i} \coloneqq \|  f \circ b_i^{-1} \| _{{k, \alpha};b_i(N_i)}.
\end{equation*}
The global $C^{k, \alpha}$ norm on $N$ is defined as the supremum over all charts:
\begin{equation*} 
\| f\| _{k, \alpha; N} \coloneqq \sup_{i \in \mathcal{I}} \| f\| _{k, \alpha; N_i}.
\end{equation*}
A tensor field $\mathbf{T}$ defined in $N$ can be represented in the standard coordinate frame of $\CC^n$, 
\begin{align*}
    \mathbf{T} = \sum_{\substack{i_1,\ldots, i_q\\j_1,\ldots, j_p}}T^{j_1\ldots j_p}_{i_1\ldots i_q} \frac{\partial}{\partial x^{j_1}} \otimes \ldots \otimes \frac{\partial}{\partial x^{j_p}} \otimes d x^{i_1} \otimes \ldots \otimes d x^{i_q}.
\end{align*}
Then, the $\C^{k,\alpha} $ norm of $\mathbf{T}$ on $N$ is
$$\| \mathbf{T}\| _{k,\alpha; N} = \sup_{i\in\I} \sum_{\substack{i_1,\ldots, i_q\\j_1,\ldots, j_p}}\| T^{j_1\ldots j_p}_{i_1\ldots i_q}\| _{k,\alpha; N_i} $$

Throughout the paper, we require Hölder norms for maps into specific fiber bundles constructed from the local geometry of $(N, \omega)$. We introduce a fiber bundle $\mathcal{W} \xrightarrow{\pi} N$ constructed by gluing local cotangent bundles $E_j \cong N_j \times \mathbb{C}^n$ (where $\mathbb{C}^n$ represents the fiber of $(1,0)$-forms) over the intersections $N_i \cap N_j$. The transition rule is defined by the difference of the reference potentials $\rho$. Specifically, a point $(z, \xi) \in E_j$ is identified with $(z, \xi') \in E_i$ via:
\begin{align} \label{def:bundle_W_transition}
(z, \xi) \sim \left(z, \xi + \partial(\rho_i - \rho_j)(z)\right).
\end{align}
{The} {formula} {is written after expressing both covectors in} the {same} base {coordinates. In distinct holomorphic coordinates, the usual cotangent} coordinate {change} {is} {also} {applied.}

We consider smooth maps $G: N \times {D} \to \mathcal{W}$. In a local chart $N_i \times {D}$, such a map is represented by a pair of functions $(z_i(w, \tau), \xi_i(w, \tau))$, where $z_i$ maps to the base $N_i$ and $\xi_i$ maps to $\mathbb{C}^n$. Write $|u|_{k,\alpha}(p)$ for the pointwise quantity in the standard Hölder norm above: fix the first point at $p$ and take the Hölder supremum over the second point in the same coordinate domain. Derivatives and distances use the real chart coordinates (and the disc variables for $G$), with $|(a,b)|:=|a|+|b|$. The global Hölder norm of $G$ is defined as:
\begin{align*}
\| G\| _{k, \alpha; N \times D}
\coloneqq \sup_{x\in N}\sup_{i:\,x\in N_i}\sup_{\tau\in D}
|(z_i,\xi_i)|_{k,\alpha}(x,\tau).
\end{align*}

{At the trivial family \(G_0\),} the linearized {bundle} is {\(\mathcal W'=G_0^*T\mathcal W\)} over {\(N\).} {Write} {\(E'_j\)} {for} {its} {local trivialization over \(N_j\). Its local fibres have} components {\((\hat z,\hat\xi)\).} {When} {both} {sides} {are} {expressed} in the {same} base {coordinates,} {the transition is
\begin{align*}
(w;\hat z,\hat\xi)\in E'_j
\sim
\left(w;\hat z,\hat\xi+D_\partial^2(\rho_i-\rho_j)(w)\cdot\hat z\right)\in E'_i,
\end{align*}
where
\[
\big(D_\partial^2(\rho_i-\rho_j)\big)_{kl}
=\partial_k\partial_l(\rho_i-\rho_j).
\]
In distinct coordinates one differentiates the full transition \eqref{def:bundle_W_transition}, including its cotangent Jacobian. Thus \(\hat z\)} and {\(\hat\xi\)} {are} the {components} {of} {a} {tangent vector to \(\mathcal W\).} For a section {\(\widehat G\in\Gamma(N,\mathcal W')\)}, represented {by} {\((\hat z_i,\hat\xi_i)\)}, we define
\begin{align*} 
\| \widehat{G}\| _{k, \alpha; N}
\coloneqq \sup_{x\in N}\sup_{i:\,x\in N_i}
|(\hat z_i,\hat\xi_i)|_{k,\alpha}(x).
\end{align*}
{These} norms are {taken in the fixed atlas; comparisons between charts use the full bundle transitions} and {the} {uniform} {chart} {and potential bounds.}

\section{Formalism via Families of Holomorphic Discs}

\subsection{Family of Holomorphic Discs and HCMA equation}\label{subsec:hcma_disc_families} The equivalence between smooth, non-degenerate solutions of the homogeneous complex Monge-Ampère (HCMA) equation and families of holomorphic discs was established by Semmes \cite{semmes1992complex} and Donaldson \cite{donaldson2002holo}. This equivalence extends naturally to the local setting (see also e.g., \cite{Chen2008geofoli, CHEN2020108603}).

Recalling the open domain $U \subset X$ defined in Section 2, we introduce a slightly larger open domain $U^+$ such that $U \Subset U^+ \subseteq X$. In what follows, $U$ will serve as the parameter space for the family of holomorphic discs, while $U^+$ will denote the total ambient region in $X$ containing the fibre of all holomorphic discs. We assume that $U^+$ is covered by a locally finite collection of holomorphic coordinate balls $\{U^+_i\}_{i \in \mathcal{I}}$. Specifically, for each $i \in \mathcal{I}$, the set $U^+_i$ is an enlarged neighborhood containing our previously defined patch $U_i$ (such that $U_i \Subset U^+_i$), and there exists a biholomorphism $b_i: U^+_i \to B^+_{R_i}$ onto a ball of radius $R_i+ \delta_i$ in $\mathbb{C}^n$.

Working within a fixed local chart (omitting the index $i$ for brevity), let $E = T^*U^+$ be the holomorphic cotangent bundle of $U^+$ with canonical projection $\pi: E \rightarrow U^+$. Since the domain is local, the bundle is trivial, $E \cong U^+ \times \mathbb{C}^n$. We denote the local holomorphic coordinates on $U^+$ by $z = (z^1, \ldots, z^n)$ and the fiber coordinates on $E$ by $\xi = (\xi^1, \ldots, \xi^n)$. The total space $E$ carries a canonical complex symplectic form defined by:
 \begin{align*}
     \Xi = dz^i \wedge d\xi^i.
 \end{align*}
 Locally, let $\rho$ be a potential for the reference metric, such that $\omega|_{U^+} = i\partial \overline{\partial} \rho$. We can associate $\omega|_{U^+}$ with a submanifold $\Lambda_\rho \subset E$ defined as the graph of the gradient section $\partial \rho$. Explicitly, $\Lambda_\rho = \{ (z, \xi) \in E \mid \xi^k = \frac{\partial \rho}{\partial z^k}(z) \}$.
 By restricting $\Xi$ on $\Lambda_\rho$, a direct calculation shows that,
 \begin{align*}
     \Xi|_{\Lambda_\rho} =  \partial \overline{\partial} \rho = -\sqrt{-1} \omega.
 \end{align*}
 Consequently, we have $\operatorname{Re} \Xi|_{\Lambda_\rho} = 0$ and $\operatorname{Im} \Xi|_{\Lambda_\rho} = -\omega$. This implies that $\Lambda_\rho$ is an exact Lagrangian submanifold with respect to the real symplectic structure $\operatorname{Re} \Xi$ on $E$. Consider the boundary data for the HCMA equation given by $\Psi$.  For each $\tau \in \partial D$, we define the potential $\tilde{\psi}_\tau = \rho + \Psi(\cdot, \tau)$. The associated metric is $\omega_\tau = \omega + i\partial \overline{\partial} \Psi(\cdot, \tau) = i\partial \overline{\partial} \tilde{\psi}_\tau$.
Let $\Lambda_{\tilde{\psi}_{\tau}} \subset E$ be the exact Lagrangian submanifold defined by the graph of $\partial \tilde{\psi}_\tau$.

We now extend this construction to the possibly non-local open domain $U^+ = \bigcup_{j \in \mathcal{I}} U^+_j$. We construct a total holomorphic fiber bundle $\mathcal{W}$ over $U^+$ by gluing the local cotangent bundles $E_j = T^* U^+_j$.
Let $\omega|_{U^+_j} = i\partial\overline{\partial}\rho_j$ be the local potential on each chart $U^+_j$. We glue the bundles $E_j$ and $E_k$ over their intersection $U^+_j \cap U^+_k$ via the transition rule:
\begin{align*} 
(z, \xi) \in E_j \sim (z, \xi + \partial(\rho_k - \rho_j)) \in E_k.
\end{align*}
This defines $\mathcal{W}$ as an affine bundle over $U^+$. Since the transition functions involve the gradient of the potential difference $\rho_k - \rho_j$ (which is holomorphic since $\partial\bar{\partial}(\rho_k - \rho_j) = 0$), the complex structure is well-defined. Within this global bundle $\mathcal{W}$, the boundary data $\Psi$ defines a global exact Lagrangian submanifold. Although the potentials $\rho_j$ are local, the sections defined by $\xi = \partial(\rho_j + \Psi)$ transform consistently under $\partial(\rho_j + \Psi) + \partial(\rho_k - \rho_j) = \partial(\rho_k + \Psi).$ Thus, these local sections patch together to form a single global submanifold $\Lambda_{\Psi(\cdot, \tau)} \subset \mathcal{W}$ for each $\tau \in \partial D$.

With this target space established, we formally define the family of discs corresponding to a solution.

\begin{define}[HCMA disc family] \label{def:disc_family}
Let $U \subseteq X$ and $\Psi$ be as above. A smooth map $G: U \times D \to \mathcal{W}$ is called an HCMA Disc Family adapted to $\Psi$ over $U$ if it satisfies the following conditions:
 \begin{itemize}
     \item let $g_x(\cdot) = G(x, \cdot): D \to \mathcal{W}$. $g_x (\tau)$ is holomorphic.
     \item for each $\tau \in \partial D $ and each $x\in U$, 
     \begin{align*}
         g_x(\tau) \in \Lambda_{\Psi(\cdot, \tau)};
     \end{align*}
     \item for each $x\in U$, we have $\displaystyle \pi\circ g_x (-i) = x$;
     \item let $Z (x,\tau) = \pi \circ G(x, \tau) $ and $z_\tau(x) = Z(x, \tau)$. For every $\tau \in D$, the map $z_\tau: U \to X$ is a diffeomorphism onto its image, and there exists a subdomain $U^- \Subset U$ such that $U^- \subseteq z_\tau(U)  \subseteq U^+$ for all $\tau$.
 \end{itemize}
\end{define}

To facilitate the analysis, we express the HCMA disc family in terms of holomorphic coordinates. Let $w = (w^1, \ldots, w^n)$ be complex coordinates on the parameter domain $U$. We denote the components of the map $g_w(\tau)$ locally in the base space coordinates and the fiber bundle coordinates by
$$ G(w, \tau)|_{E} = (z(w, \tau), \xi(w, \tau)). $$
When \(U\) is not contained in a single coordinate chart, the above description is understood locally on each chart of a covering of \(U\), and the local representations are patched together via the transition maps of the bundle \(\mathcal W\).

Since the projection of the HCMA disc family, $Z(w,\tau) = z(w, \tau)$, defines a foliation of the image domain by holomorphic discs, the map defines a coordinate transformation
$$(w, \tau) \longrightarrow (z, \tau).$$ 
We refer to the system $(w, \tau)$ as the \textit{foliated coordinates}, while $(z, \tau)$ denotes the standard product coordinates on $X \times D$.

Locally, the geometric conditions translate into the following system of partial differential equations. For $k = 1, \ldots, n$:
\begin{align}\label{holodiscpde}
\begin{split}
    \frac{\partial}{\partial \bar{\tau}} z^k (w, \tau) = \frac{\partial}{\partial \bar{\tau}} \xi^k (w, \tau) = 0, \qquad & (w, \tau) \in U \times D; \\
    \xi^{k} (w, \tau) = \Big(\frac{\partial}{\partial{z^k}} (\rho+ {\psi}_\tau)\Big) (z(w, \tau), \tau), \qquad & (w, \tau) \in U \times \partial D;\\
    {z} (w, -i) =  w, \qquad \qquad \qquad & \ w \in U,
\end{split}
\end{align}
and $z(w, \tau)$ is a family of diffeomorphisms of $U$ satisfying, 
\begin{align*}
U^- \subseteq \{z(w, \tau)| \ w\in U\} \subseteq U^+,\qquad \text{ for }  \tau \in D.
\end{align*}

\begin{define}
    \label{def:smooth_sol}
    We say that $\Phi$ is a smooth nondegenerate solution to the HCMA equation \eqref{introHCMA} on $U^- \times D$ if $\Phi \in C^\infty(U^- \times D)$ and the transversality condition holds:$$\omega + i\partial \overline{\partial} \varphi_\tau > 0 \quad \text{on } U^-$$
    for each $\tau \in \overline{D}$, where $\varphi_\tau = \Phi(\cdot, \tau)$.
\end{define}

The correspondence between smooth nondegenerate solutions and HCMA disc families is based on the kernel foliation associated with the degenerate form
\(
\Omega_\Phi:=\Omega+i\partial\bar\partial\Phi.
\) 
Indeed, the HCMA equation {\((\Omega+i\partial\bar\partial\Phi)^{n+1}=0\)} {and} {nondegeneracy} of the spatial slices {imply} {that} {\(\Omega_\Phi=\Omega+i\partial\bar\partial\Phi\)} has {a} complex {one-dimensional} {kernel} transverse to the {slices.} {Its} integral curves {form} {a} {local} {foliation} {by} holomorphic {curves} \cite{bedford1977Foli}{,} along which {\(\rho+\Phi\)} {is} {harmonic.} {When} the {leaves} {issuing} {from} {\(U\times\{-i\}\)} {extend} {over} the {whole} {disc} {inside} {\(U^+\times D\),} {with} the {projection property in Definition~\ref{def:disc_family}, they lift to an} HCMA disc family {by}
\begin{align}\label{construcHCMADiscFamilybundle}
\xi(w, \tau) = \partial (\rho + \Phi) (z(w, \tau), \tau), \qquad \tau \in D.
\end{align}
Finally, by patching the local data $(z, \xi)$ defined on each holomorphic chart, one obtains the HCMA disc family, $G: U \times D \rightarrow \W$ with $G(w, \tau)=(z(w, \tau), \xi(w, \tau)) $ locally, \textit{associated with} the smooth nondegenerate solution $\Phi$.

Conversely, starting from an HCMA disc family in the sense of Definition \ref{def:disc_family}, one can reconstruct a smooth nondegenerate solution of the HCMA equation locally. Since the reference potential \(\rho\) is only locally defined, the construction is carried out on each holomorphic coordinate chart.
On a chart with reference potential $\rho$, we define the local total potential $\widetilde{\Phi}(z, \tau) = (\rho + \Phi)(z, \tau)$ as the solution to the following Dirichlet problem along each disc:
\begin{align} \label{explicitconstructionsolu}
\begin{split}
\frac{\partial}{\partial {\tau}} \frac{\partial}{\partial \bar{\tau}} \big(\widetilde{\Phi}(z(w, \tau), \tau)\big) = 0,\quad\qquad & \tau \in D,\\
    \widetilde{\Phi}(z(w, \tau), \tau) = (\rho +\Psi) (z(w, \tau), \tau), \quad &\tau \in \partial D.
    \end{split}
\end{align}
Geometrically, this defines $\widetilde{\Phi}$ by requiring it to be harmonic along each leaf $z(w, \cdot)$ and to agree with the prescribed boundary potentials. 

The desired local solution is then obtained by setting $\Phi = \widetilde{\Phi} - \rho$. Note that on the intersection of two charts, the difference between reference potentials, $\rho_j - \rho_k$, is pluriharmonic. Since the discs are holomorphic curves, this difference remains harmonic when restricted to the leaves, ensuring that $\Phi$ is consistently defined across charts.

In the following Proposition, we will complete the proof of the equivalence:

\begin{pro}\label{propexplicitconstructionsolu}
Let \(G:U\times D\longrightarrow \mathcal W\) be an HCMA disc family adapted to the boundary data \(\Psi\). In a local holomorphic coordinate chart, write
\(G(w,\tau)=(z(w,\tau),\xi(w,\tau)),\) and let \(\widetilde\Phi\) be the local total potential constructed by the Dirichlet problem \eqref{explicitconstructionsolu}. Then, for each \(i=1,\ldots,n\),
\begin{align}\label{holoconditionfirstderivatives}
\frac{\partial}{\partial\bar\tau}
\left[
\frac{\partial \widetilde\Phi}{\partial z^i}
\bigl(z(w,\tau),\tau\bigr)
\right]
=0,
\qquad \tau\in D .
\end{align}

More precisely, for each \(\tau\in D\), let
\(\Lambda_{\Phi(\cdot, \tau)} \subset \W\) locally be the graph of $\partial_z \widetilde{\Phi}(\cdot, \tau)  $ in $\W$. 
Then, \(\Lambda_{\Phi(\cdot, \tau)}\) is an exact Lagrangian submanifold of \(\mathcal W\). Moreover, after patching the local potentials, the function
\[
\Phi=\widetilde\Phi-\rho
\]
is globally well-defined on \(U^-\times D\) and is a smooth nondegenerate solution of the HCMA equation \eqref{introHCMA} there.

Conversely, suppose that \(\Phi\) is a smooth nondegenerate solution of {\eqref{introHCMA}} on \(U^+\times D\){, and that the kernel leaves with fixed starting points in \(U\times\{-i\}\) extend over all of \(D\) inside \(U^+\times D\), with the projection property in Definition~\ref{def:disc_family}}. Then \(\Phi\) determines an HCMA disc family
\[
G:U\times D\longrightarrow \mathcal W
\]
adapted to the boundary data \(\Psi\).
\end{pro}

\begin{proof}
We present a concise geometric proof based on a change of coordinates, avoiding the more involved direct computations. Starting from the HCMA disc family
\(
G:U\times D\longrightarrow \mathcal W,
\)
we first prove \eqref{holoconditionfirstderivatives}, and hence shows that the locally constructed potential \(\Phi=\widetilde\Phi-\rho\) is a smooth nondegenerate solution of the HCMA equation.

Let $h(w, \tau) = \widetilde{\Phi}(z(w, \tau), \tau)$. By construction, for each fixed $w$, the function $\tau \mapsto h(w, \tau)$ is harmonic in $D$. Differentiating the equation $\Delta_\tau h = 0$ with respect to the parameter $w^j$, we see that the variation
$$A_j(w, \tau) := \frac{\partial}{\partial w^j} \big( h(w, \tau) \big)$$
is also harmonic with respect to $\tau$.

Consider the comparison function defined by the contraction of the fiber coordinates and the tangent of the discs:
$$B_j(w, \tau) := \xi^k(w, \tau) \frac{\partial z^k}{\partial w^j}(w, \tau) + \overline{\xi^k(w, \tau)} \frac{\partial \bar{z}^k}{\partial w^j}(w, \tau).$$
Since the HCMA disc family requires $z(w, \cdot)$ and $\xi(w, \cdot)$ to be holomorphic in $\tau$, the term $\xi^k \frac{\partial z^k}{\partial w^j}$ is holomorphic (and thus harmonic). Similarly, the conjugate term is anti-holomorphic (and thus harmonic). Therefore, $\Delta_\tau B_j = 0$.

We now compare the boundary values on $\partial D$.
Recall that $h|_{\partial D} = (\rho + \Psi)(z(w, \tau), \tau)$. Differentiating with respect to $w^j$:
$$A_j|_{\partial D} = \frac{\partial (\rho + \Psi)}{\partial z^k} \frac{\partial z^k}{\partial w^j} + \frac{\partial (\rho + \Psi)}{\partial \bar{z}^k} \frac{\partial \bar{z}^k}{\partial w^j}.$$
From the boundary condition of the disc family, we have $\xi^k = \frac{\partial (\rho + \Psi)}{\partial z^k}$ on $\partial D$. Substituting this into the expression for $B_j$, we see that
$$B_j|_{\partial D} = \frac{\partial (\rho + \Psi)}{\partial z^k} \frac{\partial z^k}{\partial w^j} + \overline{\left(\frac{\partial (\rho + \Psi)}{\partial z^k}\right)} \frac{\partial \bar{z}^k}{\partial w^j} = A_j|_{\partial D}.$$
Since $A_j$ and $B_j$ are both harmonic in $D$ and agree on $\partial D$, by the uniqueness of the Dirichlet problem, they must be identical in $D$:
\begin{align*} 
    \frac{\partial}{\partial w^j} \big(\widetilde{\Phi}(z(w, \tau), \tau)\big) = \xi^k \frac{\partial z^k}{\partial w^j} + \bar{\xi}^k \frac{\partial \bar{z}^k}{\partial w^j}.
\end{align*}
Equivalently, as an identity of one-forms in the \(w\)-variables,
\[
d_w\big(\widetilde{\Phi}(z(w,\tau),\tau)\big)
=
z_\tau^*
\left(
\xi^k\,dz^k+\bar\xi^k\,d\bar z^k
\right),
\]
where \(z_\tau(w)=z(w,\tau)\). On the other hand, note that
\(
d_w\big(\widetilde{\Phi}(z(w,\tau),\tau)\big)
=
z_\tau^*
\left(
\frac{\partial\widetilde\Phi}{\partial z^k}\,dz^k
+
\frac{\partial\widetilde\Phi}{\partial \bar z^k}\,d\bar z^k
\right).
\)
Since \(z_\tau\) is a diffeomorphism onto its image, its pullback is injective on one-forms. Hence,
\[
\frac{\partial\widetilde\Phi}{\partial z^k}
\big(z(w,\tau),\tau\big)
=
\xi^k(w,\tau).
\]
Therefore, using the holomorphicity of \(\xi^k(w,\tau)\) in \(\tau\), we obtain
\[
\frac{\partial}{\partial\bar\tau}
\left[
\frac{\partial\widetilde\Phi}{\partial z^k}
\big(z(w,\tau),\tau\big)
\right]
=
\frac{\partial\xi^k}{\partial\bar\tau}(w,\tau)
=0,
\]
which is \eqref{holoconditionfirstderivatives}.

The condition \eqref{holoconditionfirstderivatives} follows immediately from the holomorphicity of $\xi^k$:
$$\frac{\partial}{\partial \bar{\tau}} \left( \frac{\partial \widetilde{\Phi}}{\partial z^k} \right) (z(w, \tau), \tau) = \frac{\partial \xi^k}{\partial \bar{\tau}} (w, \tau) = 0.$$

The smoothness of $\widetilde{\Phi}$ directly follows from \cite[Lemma~2.8]{yao2025asymptotic} by solving a family of Dirichlet problems with smooth boundary data.

To see the nondegeneracy of $\partial\overline{\partial}\widetilde{\Phi}$ in the spatial direction, we calculate the complex symplectic form restricted to each slice, $\Lambda_{\Phi(\cdot, \tau)}$. Note that
$$-\sqrt{-1} \omega_{\tau} = \Xi|_{\Lambda_{\Phi(\cdot, \tau)}} = dz \wedge d\xi |_{\Lambda_{\Phi(\cdot, \tau)}} = \partial \overline{\partial} \widetilde{\Phi}|_{z_\tau(U) \times \{\tau\}}, \quad \tau \in D. $$
Let $z_\tau : w \rightarrow z(w, \tau)$. The pullback of the form to the foliated coordinates is given by $z_\tau^* (\omega_\tau) = \sqrt{-1} d_w z(w, \tau) \wedge d_w \xi(w, \tau)$.
The spatial nondegeneracy is a consequence of the holomorphic dependence of the foliation on $\tau$. Differentiating the pullback form with respect to $\bar{\tau}$, and using that both $z(w, \tau)$ and $\xi(w, \tau)$ are holomorphic in $\tau$, we obtain:
\begin{align*}
    \frac{\partial}{\partial \bar{\tau}} (z_\tau^* (\omega_\tau)) = \sqrt{-1} d_w \Big( \frac{\partial z}{\partial \bar{\tau}}(w, \tau) \Big) \wedge d_w \xi(w, \tau) +\sqrt{-1} d_w z(w, \tau) \wedge d_w \Big( \frac{\partial \xi}{\partial \bar{\tau}} (w, \tau) \Big) = 0.
\end{align*}
This implies that the coefficients of $z_\tau^* (\omega_\tau)$, when expanded in the $\tau$-independent basis $\{dw^j, d\bar{w}^j\}$, are holomorphic functions of $\tau$.
Since $\omega_\tau$ is real and nondegenerate for $\tau \in \partial D$, we derive that these coefficients are constants. Hence, $\omega_\tau$ is nondegenerate for each $\tau \in D$ (see also \cite{donaldson1999symmetric}).

It suffices to check that $\Phi$ is the solution to the HCMA equation. Let $X_\tau$ be the $(1,0)$ vector field tangential to each leaf by $w \rightarrow z(w, \cdot)$. We write $\partial \overline{\partial} \widetilde{\Phi}$ in terms of the dual basis of $\{X_\tau, \frac{\partial}{\partial z^1}, \ldots, \frac{\partial}{\partial z^n}\}$, denoted by $\{\theta^\tau, \theta^1, \ldots, \theta^n\}$. Note that 
\begin{align*}
    d\tau = \theta^\tau, \qquad dz^k = \theta^{k} + \frac{\partial z^k}{\partial \tau} \theta^\tau, \quad \text{for } k =1,\ldots, n.
\end{align*}
Then,
\begin{align}\label{Omegaformonleafbasis} 
\begin{split}
\partial \overline{\partial} \widetilde{\Phi} = &   \widetilde{\Phi}_{\tau\bar{\tau}} d\tau \wedge d\bar{\tau}  +  \widetilde{\Phi}_{j\bar{\tau}} dz^j \wedge d\bar{\tau} +  \widetilde{\Phi}_{\tau\bar{k}} d\tau \wedge d\bar{z}^k + \widetilde{\Phi}_{j\bar{k}} dz^j \wedge d\bar{z}^k \\
 = & \Big(\widetilde{\Phi}_{\tau\bar{\tau}} + \widetilde{\Phi}_{\tau \bar{k}}\frac{\partial \bar{z}^k}{\partial \bar{\tau}} + \widetilde{\Phi}_{j\bar{\tau}} \frac{\partial {z}^j}{\partial {\tau}} + \widetilde{\Phi}_{j\bar{k}} \frac{\partial {z}^j}{\partial {\tau}} \frac{\partial \bar{z}^k}{\partial \bar{\tau}} \Big) \theta^\tau \wedge \bar{\theta}^\tau \\
&+ \Big( \widetilde{\Phi}_{j\bar{\tau}} + \widetilde{\Phi}_{j\bar{k}} \frac{\partial \bar{z}^k}{\partial {\bar{\tau}}} \Big) \theta^j \wedge \bar{\theta}^\tau 
+  \Big( \widetilde{\Phi}_{\tau\bar{k}} + \widetilde{\Phi}_{j\bar{k}} \frac{\partial {z}^j}{\partial {{\tau}}} \Big) \theta^\tau \wedge \bar{\theta}^k \\
&+ \widetilde{\Phi}_{j\bar{k}}  \theta^j \wedge \bar{\theta}^k
\end{split}
\end{align}

Based on the equation (\ref{holoconditionfirstderivatives}), we have,
\begin{align*}
    \widetilde{\Phi}_{j\bar{\tau}} + \widetilde{\Phi}_{j\bar{k}} \frac{\partial \bar{z}^k}{\partial {\bar{\tau}}} = 0.
\end{align*}
In other words, the coefficients of tangential vector $X_\tau$ is given by $ \frac{\partial z^j}{\partial \tau} = -\widetilde{\Phi}^{j\bar{k}} \widetilde{\Phi}_{\tau\bar{k}}$. 
According the construction of $\widetilde{\Phi}$, (\ref{explicitconstructionsolu}), we have,
\begin{align}
\begin{split}\label{con:harmoniconleaf}
    0 &= \frac{\partial}{\partial {\tau}} \frac{\partial}{\partial \bar{\tau}} \big(\widetilde{\Phi}(z(w, \tau), \tau)\big)\\
    & = \widetilde{\Phi}_{\tau\bar{\tau}} + \widetilde{\Phi}_{\tau \bar{k}}\frac{\partial \bar{z}^k}{\partial \bar{\tau}} + \widetilde{\Phi}_{j\bar{\tau}} \frac{\partial {z}^j}{\partial {\tau}} + \widetilde{\Phi}_{j\bar{k}} \frac{\partial {z}^j}{\partial {\tau}} \frac{\partial \bar{z}^k}{\partial \bar{\tau}}.
\end{split}
\end{align}
Hence, given the HCMA Disc Family $G$ and local potential $\widetilde{\Phi}$ constructed in \eqref{explicitconstructionsolu}, we have $\partial \overline{\partial} \widetilde{\Phi} = \widetilde{\Phi}_{j\bar{k}}  \theta^j \wedge \bar{\theta}^k $; hence $\big(\sqrt{-1}\partial \overline{\partial} \widetilde{\Phi}\big)^{n+1} =0 $. 

Conversely, {under the} {continuation} {and} {projection} {hypotheses,} {we have} {the} {family of base leaves \(z(w,\tau)\) over \(U\times D\). It remains} to {verify that} the {lift} {\(\xi(w,\tau)\)} {in} {\eqref{construcHCMADiscFamilybundle}} is holomorphic in {\(\tau\)}. {The} {form} {\(\partial\bar\partial\widetilde\Phi\)} {is} written as {in} \eqref{Omegaformonleafbasis}.
Since by restricting to each leaf $w \rightarrow z(w, \cdot)$, $\sqrt{-1} \partial \overline{\partial}{\widetilde{\Phi}} =0$, the condition \eqref{con:harmoniconleaf} holds here. Therefore, 
\begin{align*}
    0=\big(\sqrt{-1}\partial \overline{\partial} \widetilde{\Phi}\big)^{n+1} = - \widetilde{\Phi}^{j\bar{k}} \Big( \widetilde{\Phi}_{\tau\bar{k}} + \widetilde{\Phi}_{l\bar{k}} \frac{\partial {z}^l}{\partial {{\tau}}} \Big) \Big( \widetilde{\Phi}_{j\bar{\tau}} + \widetilde{\Phi}_{j\bar{p}} \frac{\partial \bar{z}^p}{\partial {\bar{\tau}}} \Big) \det \big(\widetilde{\Phi}_{j\bar{k}}\big) \vol.
\end{align*}
Then, by \eqref{construcHCMADiscFamilybundle}, we have
$$\frac{\partial\xi^j }{\partial \bar{\tau}}  (w, \tau) = \widetilde{\Phi}_{j\bar{\tau}} + \widetilde{\Phi}_{j\bar{p}} \frac{\partial \bar{z}^p}{\partial {\bar{\tau}}} =0.$$
The $z(w, \tau)$ is holomorphic in $\tau$, due to the fact that $z(w, \cdot) $ is a holomorphic integral curve by \cite{bedford1977Foli}. 
\end{proof}

\subsection{Perturbation Theory and Linearized Operator}
We now develop the perturbation theory for HCMA disc families. Starting from the trivial HCMA disc family, we show that any sufficiently small perturbation of the boundary data produces a unique nearby family of holomorphic discs. This provides the basic stability result needed for the local regularity argument.

We begin by explicitly writing down the trivial foliation in local coordinates. In each local holomorphic chart $U_i$, $i \in \I$, we fix a K\"ahler potential $\rho$. We consider the trivial case where the perturbation (boundary data) is identically zero, $\Psi \equiv 0$. In this setting, the trivial HCMA disc family is locally given by
$$(z(w, \tau), \xi(w, \tau)) = (w, \partial \rho(w)).$$
Note that in this base case, the coordinate $z$ is independent of $\tau$, the projection down of HCMA disc family defines a trivial foliation of $X\times D$ by horizontal discs.

To analyze {perturbations} of {the} {boundary} {data,} {consider} the {following} space of {normalized holomorphic} disc {maps,} {restricted} {to} {those satisfying the projection condition in Definition~\ref{def:disc_family}}:
\begin{align} \label{defspaceholodiscs}
    \F =\{ G \in \C^\infty(U \times D, \W); \ \overline{\partial}_\tau G  =0,\  \pi\circ G (x, -i) =x\}.
\end{align}
The boundary conditions are the smooth families of exact Lagrangian subspaces of $\W$ parametrized by $\partial D$. Recall that the notation $\Lambda_{\Psi(\cdot, \tau)}$, $\tau \in \partial D$, denotes the graph of $\partial (\rho+{\psi}_\tau)$ locally, and patches the local pieces to get an exact Lagrangian subspace of $\W$. We denote the space of boundary conditions as follows
    \begin{align*}
    \B = \{ \Lambda_{{\Psi}} \subseteq \W \times \partial D; \ {\psi}_\tau (z)= {\Psi}(z, \tau)  \in \C^{\infty} (U^+ \times \partial D), \ \omega + i\partial \overline{\partial} {\psi}_\tau >0, ~ \tau \in \partial D \}.
    \end{align*}
{Write \(Z=\pi\circ G\). In a fixed local trivialization, let \(\mathcal D=\C^\infty(U\times\partial D,\CC^n)\) denote the space of boundary mismatch data. The mismatching operator \(T:\mathcal F\times\mathcal B\to\mathcal D\) is defined by}
\begin{align} \label{defopdiffholodiscbdry}
    T(G, \Lambda_{{\Psi}})(w, \tau) = \xi(w,\tau) - \partial \big( \rho + {\Psi}\big) (z(w, \tau), \tau).
\end{align}
{The equation \(T(G,\Lambda_\Psi)=0\) means that the discs satisfy the prescribed boundary condition, independently of the choice of local coordinates.} {Derivatives} of {its} {components} {below are computed in a fixed local trivialization. We} have the following lemma:

\begin{lem} Let $G \in \F$ and $\Lambda_{{\Psi}} \in \mathcal{B}$. Then, $(G, \Lambda_{{\Psi}})$ lies in the kernel of $ T$
    if and only if 
$G(w,\tau) = (z(w, \tau), \xi(w, \tau))$ is a solution to the equation (\ref{holodiscpde}) with  $ G(w,\tau ) \in \Lambda_{{\Psi}(\cdot, \tau)}$, $\tau \in \partial D$; or equivalently, $G$ is an HCMA disc family adapted to $\Psi$.
\end{lem}

{Fix} {\((G,\Lambda_\Psi)\in\mathcal F\times\mathcal B\)}. The {infinitesimal} {variations} of {\(G\) are sections of \(G^*T\mathcal W\). To retain the notation \(\mathcal W'\), regard the tangent bundle along \(G\) as a bundle over the projected family: its fibre} at {\((Z(w,\tau),\tau)\)} {is \(T_{G(w,\tau)}\mathcal W\). This is well-defined because every \(z_\tau\) is a diffeomorphism onto its image. We write \(Z^*\mathcal W'\) for its pullback} by {\((w,\tau)\mapsto(Z(w,\tau),\tau)\);} {thus} {\(Z^*\mathcal W'=G^*T\mathcal W\).} {The} {transition rules are the differentials} of the {full} {affine} {cotangent} {transitions} of {\(\mathcal W\)}. {In} {particular,}
{\[
\mathcal F'_G
=\left\{\widehat G\in\Gamma(U\times D,Z^*\mathcal W'):
\bar\partial_\tau\widehat G=0,\quad
p_1\circ\widehat G(w,-i)=0\right\},
\]}
where {\(p_1=d\pi\)} is the projection to {the base variation.}
We now compute the derivative of the operator $T$ with respect to the disc configuration $G$ in the direction of $\widehat{G} \in \mathcal{F}'_G$. Let the variation be written in local coordinates as
\(\widehat{G}(w, \tau) = (\hat{z}(w, \tau), \hat{\xi}(w, \tau)).\)
Then the derivative of $T$ in the direction of $\F$ can be written in the complex coordinates by, 
\begin{align*}
    \big(D_{\F} T |_{(G, \Lambda_{{\Psi}})} (\widehat{G})\big)^{i} = \hat{\xi}^i (w, \tau) - \sum_{j=1}^n \Big\{ (\partial_i\overline{\partial}_{j} (\rho+ {\Psi}) ) (z(w, \tau), \tau) \overline{\hat{z}}^j + (\partial_i \partial_j (\rho +{\Psi}))(z(w, \tau), \tau) \hat{z}^j\Big\}, 
\end{align*}
for each $\tau \in \partial D$, where $ (D_{\F} T |_{(G, \Lambda_\Psi)} (\widehat{G}))^{i}$ represents the $i$-th coordinate of the standard bundle coordinates of $E$.
Consequently, {in the chosen local trivialization} the linearized problem {has} the following {form. At a zero} of {\(T\),} the {component} {derivatives} {define} {the intrinsic differential. Away from a zero, differentiating a coordinate change also differentiates the covector transformation acting on the residual; the displayed formula is therefore used} in {a} {fixed} {chart.}
    \begin{align}\label{holodisclinearpde}
    \begin{split}
        &\overline{\partial}_\tau \hat{z}^i (w, \tau) = \overline{\partial}_\tau \hat{\xi}^i (w, \tau) = 0, \hspace{5.7cm} \text{ in }  U \times D,\\
        \begin{split}
         &\hat{\xi}^i(w,\tau) - \big(\partial_i \overline{\partial}_{{j}} (\rho+ {\Psi}) \big)(z(w, \tau), \tau) \overline{\hat{z}}^j(w, \tau)\\ 
         &\hspace{2cm}- \big(\partial_{i}\partial_{j} (\rho+ {\Psi})\big) (z(w, \tau), \tau) \hat{z}^j(w, \tau)=\bff^i(w,\tau),
         \end{split}
         \hspace{0.5cm}  \text{ in } U \times \partial D,\\
         &\hat{z}^i(w, -i) = 0,  \hspace{8cm}  \ w  \in U,
    \end{split}
\end{align}
where $\bff= (\bff^1, \ldots, \bff^n) \in \C^{\infty}(U \times \partial D, \CC^n)$.  

We now turn to the solvability of the linearized problem in the case of trivial data. We consider the linearization around the trivial solution $G_0(w, \tau) = (w, \partial \rho(w))$ with $\Psi \equiv 0$. In this setting, the coefficients of the linearized operator simplify significantly. Since the background discs are constant maps $z(w, \tau) = w$, the derivatives of the potential $\rho$ are independent of $\tau$. The general linearized system (\ref{holodisclinearpde}) reduces to:
\begin{align}\label{holodisclinearpde-trivial}
    \begin{split}
        &\overline{\partial}_\tau \hat{z}^i (w, \tau) = \overline{\partial}_\tau \hat{\xi}^i (w, \tau) = 0, \hspace{5.2cm} \text{ in }U \times D,\\
         &\hat{\xi}^i(w,\tau) - A_{i\bar{j}}(w) \overline{\hat{z}}^j(w, \tau)- B_{ij}(w) \hat{z}^j(w, \tau)=\bff^i(w,\tau),
         \hspace{0.8cm}  \text{ in } U \times \partial D,\\
         &\hat{z}^i(w, -i) = 0,  \hspace{7.5cm}  \ w  \in U,
    \end{split}
\end{align}
where $A(w)_{i\bar j} = \partial_i \overline{\partial}_j \rho(w)$ is the complex Hessian of $\rho$, and $B(w)_{ij} = \partial_i \partial_j \rho(w)$.

To solve the system, we introduce a transformation that decouples the equations. Following the approach in Chen-Feldman-Hu \cite{CHEN2020108603} (Lemma A.3), we define two auxiliary vector-valued functions $\hat{g}_1$ and $\hat{g}_2$ on $U \times D$ via the following linear combination:
\begin{align}\label{defg1g2}
\begin{split}
\hat{g}_{1}(w, \tau) &= \hat{\xi}(w, \tau) - B(w) \hat{z}(w, \tau) - \bar{A}(w) \hat{z}(w, \tau),  \\
\hat{g}_2(w, \tau) &= \hat{\xi}(w, \tau) - B(w) \hat{z}(w, \tau) + \bar{A}(w) \hat{z}(w, \tau).
\end{split}
\end{align}
Using this change of variables, the linearized system (\ref{holodisclinearpde-trivial}) decouples into the following problem:
\begin{align} \label{equstandardfbpros}
\begin{split}
&\overline{\partial}_\tau \hat{g}_1  = \overline{\partial}_ 
\tau \hat{g}_2  = 0, \hspace{3.15cm} \text{ in } U \times D,\\
&\text{Re}( \hat{g}_1) = \re(\bff), \hspace{3.4cm} \text{ in } U \times \partial D,\\
&\text{Im}(\hat{g}_2) = \im(\bff),  \hspace{3.4cm}\text{ in } U \times \partial D,\\
&\hat{g}_1 (w, -i) = \hat{g}_2 (w, -i) = \bff(w,-i), \hspace{0.5cm} \text{ for } w \in U.
\end{split}
\end{align}
This formulation represents a family of double Riemann-Hilbert problems parametrized by $w \in U$. The strategy for solving such disc problems has evolved through several key works. Bishop \cite{Bishop1965} originally established the existence theory for individual holomorphic discs. Later, Donaldson \cite{donaldson2002holo} and Chen-Feldman-Hu \cite{CHEN2020108603} extended the analysis to families of discs, addressing the regularity of the solution in terms of the families of linearized problems. However, a significant gap remains in the existing literature: these works did not fully realize and treat the loss of regularity in the parameter direction. This phenomenon is critical to the analysis of the linearized operator and will be discussed in detail in the next section.

\section{Analysis of the Linearized Problem}

\subsection{Loss of Regularity and Solution to the Linearized Problem at Trivial Data} \label{seclocalexholodiscs} 
In this subsection, our primary objective is to analyze the linearized equation \eqref{holodisclinearpde-trivial} at the trivial data $(G_0, \Psi\equiv 0)$ over $U$. Recall that the transformation \eqref{defg1g2} decouples this system into the pair of independent equations given in \eqref{equstandardfbpros}. Consequently, the fundamental analytical building block is the scalar problem for a holomorphic function with a prescribed real part on the boundary. We refer to this as a family of Riemann-Hilbert problems parametrized by $w \in U$. Precisely, we consider the following model problem:
 \begin{align}\label{dbardisk}
 \begin{split}
     \overline{\partial}_\tau u =0, \qquad &\text{ in } U \times  D,\\
     \re u = f, \qquad &\text{ on } U \times \partial D.
 \end{split}
 \end{align}
 Here we assume $f$ is a smooth real function defined on $U \times \partial D$. Based on complex Fourier expansion on the unit circle, we can explicitly write down the solutions to (\ref{dbardisk})
 \begin{align*}
 u(w,\tau) = \frac{1}{2\pi} \int^{2\pi}_0 f(w, e^{i\theta}) \frac{1+e^{-i\theta}\tau}{1- e^{-i\theta} \tau} d\theta + i A(w),
 \end{align*}
 where $A(w)$ is an arbitrary real function on $U$. The function $A(w)$ is the imaginary part of $u$ at $\tau =0$. In what follows, we refer to $A(w)$ as a fixed-point data (or a fixed point condition). Moreover, once $A(w)$ is fixed, the solution $u$ is uniquely determined.
 It is easily observed that the real part of the solution, $\re u$, is harmonic by integrating with the Poisson kernel. The imaginary part $\im u$ is a harmonic conjugate of $\re u$. By restricting $\im u$ on the boundary $\partial D$, the function $\im u |_{\partial D}$ is the \textit{Hilbert transform} of $f$ plus a fixed-point data $A(w)$. The \textit{Hilbert transform} is defined by the principal value of a singular integral:
 \begin{align*}
     \Hilb f (\xi) = \frac{1}{2\pi} \lim_{\varepsilon \rightarrow 0} \int_{\varepsilon <|\xi-\theta|< \pi} \cot \left(\frac{\xi-\theta}{2}\right) f(\theta) d\theta.
 \end{align*}

 It is well-known that the Hilbert transform preserves the H\"older classes on the unit circle;  precisely, $\Hilb: \C^{k,\alpha} (\partial D) \rightarrow \C^{k,\alpha} (\partial D) $ is a bounded linear operator (see \cite[Chapter III]{garnett2006bounded} for details). However, the situation is different for a family of functions on $U \times \partial D$. Although the Hilbert transform acts fiber-wise, preserving regularity in the $\tau$-direction, the Hilbert transform on $\partial D$ might cause the loss of regularity in the parameter direction $U$. To describe the loss of regularity, we introduce the notion of a \textit{modulus of continuity} for a function $g$ defined on a metric space $(M, d)$:
 $$\upsilon_{g} (t) = \sup_{d(p,p')\leq t} |g(p)- g(p')|.$$
Then, we have the following lemma:
 \begin{lem}\label{lemhilbtransesti}
     Given a function $f(x, \tau) \in \C^{k,\alpha}(U \times \partial D) $, the Hilbert transform of $f(x, \tau)$ in $\tau$ satisfies:
     \begin{align*}
         \sup_{|x-x'| \leq t}\big|\nabla_0^k \big(\Hilb_\tau f\big)(x,\tau) -\nabla_0^k \big(\Hilb_\tau f\big)(x', \tau)\big| \leq C t^\alpha(1+|\log t|) \| f\| _{k,\alpha; U \times \partial D}.
     \end{align*}
     Furthermore, for each $0<\beta < \alpha$,
     \begin{align*}
         \| \Hilb_\tau f\| _{k,\beta; U \times \partial D} \leq C (\alpha-\beta)^{-1} \| f\| _{k,\alpha; U \times \partial D}
     \end{align*}
     where $C$ is a uniform constant depending on $n$, $k$ and $\alpha$.
 \end{lem}
 \begin{proof} Without loss of generality, we can assume $U$ is a holomorphic chart hence a domain contained in $\CC^n$.
     Given a function $f\in \C^{k,\alpha}(U \times \partial D)$, we discuss the continuity of modulus in both parameter and circle directions as follows:
     \begin{align*}
         \upsilon^U_{f} (t) = \sup_{|w- w'|\leq t,\ \tau \in \partial D} |f(w, \tau) - f(w', \tau)|, \qquad \upsilon^{\partial D}_{f} (t) = \sup_{w\in U,\ |\tau-\tau'| \leq t} |f(w, \tau) - f(w, \tau')|,
     \end{align*}
     where the distances are induced by the standard metrics of $U\subseteq \CC^n$ and $\partial D$.
     
     If \(f\in C^{k,\alpha}(U\times \partial D)\), then the Hilbert transform
     commutes with differential operators up to order \(k\) in both the parameter and circle directions. Indeed, for derivatives in the circle variable this follows from the Fourier multiplier representation,
     while for derivatives in the parameter variable it follows from difference quotients and the boundedness of
     \(\Hilb_\tau:C^\alpha(\partial D)\to C^\alpha(\partial D)\).
     Therefore it suffices to prove the result for \(k=0\). For the sake of simplicity, we denote $\tilde{f} = \Hilb_\tau f$. 
     In the circle direction, according to \cite[Chapter III, Theorem 1.3]{garnett2006bounded}, we have the following inequality: for a sufficiently small constant, $0<\delta \ll 1 $,
     \begin{align} \label{esticircle}
         \upsilon^{\partial D}_{\tilde{f}(w, \cdot)} (\delta) \leq C (\delta^\alpha +  \delta) \| f(w, \cdot)\| _{0,\alpha; \partial D}, \qquad \text{for  } w \in U.
     \end{align}
     where $ C$ is a uniform constant depending only on $\alpha$. In the parameter direction, we have
     \begin{align} \label{estipara}
     \begin{split}
         |\tilde{f} (w,\tau) - \tilde{f}(w',\tau)| \leq & \ C \int_{\tau-\delta}^{\tau + \delta} \frac{|f(w, \theta) - f(w, \tau) |}{|\theta -\tau| }  d\theta  +C \int_{\tau-\delta}^{\tau + \delta} \frac{|f(w', \theta) - f(w', \tau) |}{|\theta -\tau| }  d\theta 
         \\
         &+ \frac{1}{2\pi} \int_{\delta<|\theta-\tau|< \pi} {|f(w, \theta) - f(w', \theta)|} \left|\cot \big(\frac{\theta-\tau}{2} \big)\right|d\theta.
    \end{split}
     \end{align}
The first two integrals are bounded by $O(\delta^\alpha)$ using the $C^\alpha$-regularity of $f$ in the $\tau$-direction. The third integral is bounded by $\displaystyle\upsilon_f^U(\delta) \int_{\delta<|\theta-\tau|<\pi} \frac{d\theta}{|\theta-\tau|} \approx \delta^\alpha \log(1/\delta)$. Hence, we have
     \begin{align*} 
         \upsilon^U_{\tilde{f}} (\delta) \leq C (1- \log \delta) \delta^{\alpha} \| f\| _{0,\alpha; U \times \partial D}.
     \end{align*}
The presence of the logarithmic factor $\log(1/\delta)$ prevents $\tilde{f}$ from retaining the full $C^\alpha$ regularity. However, since $\delta^\alpha \log(1/\delta) \le C (\alpha- \beta)^{-1}\delta^\beta$ for any $\beta < \alpha$, we retain $C^\beta$ regularity. If we take any constant $\beta$ with $0<\beta < \alpha$, \eqref{estipara} implies
     \begin{align} \label{estipara1}
         \upsilon^U_{\tilde{f}} (\delta) \leq C (\alpha-\beta)^{-1} \delta^{\beta} \| f\| _{0,\alpha; U \times \partial D}.
     \end{align}
     From \eqref{esticircle} and \eqref{estipara1}, we have 
     \begin{align*}
         \| \tilde{f}(w, \cdot)\| _{0,\beta; \partial D}  \leq C \| \tilde{f} (w, \cdot)\| _{0,\alpha; \partial D} \leq C \|  f (w, \cdot)\| _{0,\alpha; \partial D} \leq C \| f\| _{0,\alpha; U \times \partial D}, \qquad \text{ for all } w\in U,
     \end{align*}
     and 
     \begin{align*}
         \| \tilde{f} (\cdot, \tau)\| _{0,\beta; U} \leq C (\alpha-\beta)^{-1} \| f\| _{0,\alpha; U \times \partial D}, \qquad \text{for all } \tau \in \partial D,
     \end{align*}
     By \cite[Lemma~2.7]{yao2025asymptotic} (Bernstein's Theorem, concerning joint regularity via partial moduli of continuity), we complete the proof of the lemma.
\end{proof}

\begin{rem} \label{remcounterex}
    In this remark, we construct a counterexample to see that the Hilbert transform cannot be a bounded operator in $\C^{k,\alpha} (U \times \partial D)$. Let $U = (-\pi,\pi) \subseteq \RR$. Consider the following function, $f(x, \theta)$ defined on $U \times \partial D$:
    \begin{align*}
    f(x, \theta) = \begin{cases}
        -s(|\theta|)|\theta|^\alpha, \qquad &|x| \geq |\theta|, \quad 0\leq \theta\leq \pi;\\
        -s(|\theta|)|x|^\alpha, \qquad &|x| \leq |\theta|, \quad 0\leq \theta\leq \pi;\\
        s(|\theta|)|\theta|^\alpha, \qquad &|x| \geq |\theta|, \ -\pi\leq \theta\leq 0;\\
        s(|\theta|)|x|^\alpha, \qquad &|x| \leq |\theta|, \ -\pi\leq \theta\leq 0,
        \end{cases}
    \end{align*}
    where $0 \leq s(t) \leq 1$ is a smooth scaling function defined on $[0,\pi]$ such that $s(t) \equiv 1$ for $t \in [0,1]$ and $s(t) \equiv 0$ near $t = \pi$. 
    One can check that $\| f\| _{0,\alpha; U \times \partial D}\leq C $ for a uniform constant $C$. If we denote $\tilde{f} = \Hilb_\theta f$, for $0<x\ll1$, we have
    \begin{align*}
        |\tilde{f} (x, 0) - \tilde{f}(0, 0)| &= \frac{1}{\pi}\lim_{\varepsilon \rightarrow 0} \int_{\varepsilon<\theta< \pi}  \big|f(x, \theta)-f(0,\theta)\big| \cot\Big(\frac{\theta}{2}\Big) d\theta
        \\
        &=\frac{1}{\pi}\int_{0}^x |\theta|^\alpha \cot\Big(\frac{\theta}{2}\Big) d\theta + \frac{1}{\pi} \int_{x}^\pi s(|\theta|) |x|^\alpha \cot\Big(\frac{\theta}{2}\Big) d\theta
        \\
        & \geq c |x|^\alpha \big(  -\log |x| \big)
     \end{align*}
    It is clear that \( [\tilde{f}]_{0,\alpha; U \times \mathbb{T}^1} \) is unbounded. More precisely, the \(\alpha\)-Hölder seminorm of \(\tilde{f}\) in the \(U\)-direction is unbounded at \(\theta=0\).

If we instead consider the \(\alpha\)-Hölder norm in the \(U\)-direction at $(0, \theta_0)$ for some fixed \(\theta_{0}\neq 0\), then for sufficiently small \(x\) with \(0<x<|\theta_{0}|/4\), subtract $f(x,\theta_0)$ in the symmetric principal value integral near $\theta_0$. On that arc the $\theta$-derivative of $f$ is bounded by $Cx^\alpha$, so we have
\begin{align*}
\big|\tilde{f}(x,\theta_{0})-\tilde{f}(0,\theta_{0})\big|
&=\frac{1}{2\pi}\bigg|\operatorname{p.v.}\!\int_{-\pi}^{\pi}
\big(f(x,\theta)-f(0,\theta)\big)
\cot\!\Big(\frac{\theta_{0}-\theta}{2}\Big)\,d\theta\bigg|\\[0.7em]
&\le C|x|^\alpha\left(\int_0^{|\theta_0|/4}dt
+\int_{|\theta_0|/4}^{\pi}\frac{dt}{t}\right)\\[0.7em]
&\le C\Bigl(1-\log\frac{|\theta_{0}|}{2}\Bigr)|x|^{\alpha}.
\end{align*}
From this calculation we see that the \(\alpha\)-Hölder norm of \(\tilde{f}\) in the \(U\)-direction at each fixed \(\theta_{0}\neq 0\) is bounded, but the bound diverges as \(\theta_{0}\to 0\). 

Consider the Riemann-Hilbert problem \eqref{dbardisk} with the boundary data $f$ given as above. As discussed, the Hilbert transform of \(f\) in \(\theta\) exhibits a single singularity in its \(\alpha\)-H\"older seminorm at \(\theta=0\) on the circle when \(x=0\). The imaginary part of the solution then has boundary values \(\tilde{f}(x,\theta)+A(x)\), where \(A(x)\) is an arbitrary fixed-point function. Even if one chooses \(A(x)\) to cancel the singularity at \(\theta=0\) (for instance \(A(x)=-\tilde{f}(x,0)\), which is itself singular at \(x=0\)), the function \(\tilde{f}(x,\theta)+A(x)\) remains singular in the \(U\)-direction for every \(\theta\neq 0\). Thus no choice of fixed-point function restores uniform \(\alpha\)-H\"older regularity in \(x\) along \(\partial D\); the loss of regularity persists.
\end{rem}

We now establish the estimates for the solution $u$ to the model problem (\ref{dbardisk}). Recall that since $u$ is holomorphic in $\tau$, its real and imaginary parts are harmonic functions in the disc $D$. Their boundary values are given by:
$$\text{Re}(u)|_{U \times \partial D} = f, \qquad \text{Im}(u)|_{U \times \partial D} = \Hilb f + A(w).$$
Consequently, the regularity of $u$ in the interior domain $U \times D$ is determined by the regularity of these boundary traces via solving a family of Dirichlet problems. Combining the Hilbert transform estimates (Lemma \ref{lemhilbtransesti}) with standard harmonic extension estimates (\cite[Lemma~2.8]{yao2025asymptotic}), we derive the bound for the solution. Assuming the fixed-point data $A(w)$ belongs to $C^{k,\beta}(U)$, we have
 \begin{align} \label{estisolutionbabyversionfreebdry}
     \| u\| _{k,\beta; U \times D} \leq C \big\{ (\alpha - \beta)^{-1} \| f \| _{k,\alpha; U \times \partial D} + \| A\| _{k,\beta; U} \big\},
 \end{align}
 where $C$ is a uniform constant depending on $n$, $k$ and $\alpha$.

With the scalar estimates established, we now turn to the linearized problem at the trivial data. Recall that the trivial data corresponds to the holomorphic disc solution $G_0(w, \tau) = (w, \partial \rho(w))$ with vanishing pertubation $\Psi \equiv 0$. In this setting, the coefficients of the linearized operator depend only on the parameter $w \in U$ and are defined by the geometry of the defining function $\rho$: specifically, $A(w)$ represents the complex Hessian, $A_{i\bar{j}} = \partial_i \overline{\partial}_j \rho$, and $B(w)$ represents the remaining second-order derivatives, $B_{ij} = \partial_i \partial_j \rho$. By applying the decoupling transformation \eqref{defg1g2} and the scalar estimates derived in \eqref{estisolutionbabyversionfreebdry}, we obtain the following existence and regularity result for the coupled vector-valued system. We emphasize that the loss of regularity observed in the scalar case propagates to this full system, as captured by the change in indices from $\alpha$ to $\beta$ in the final estimate.

\begin{lem} \label{IFTlinearpart} Let $A= (A_{i\bar{j}}) \in \C^{\infty}(U, \CC^{n \times n})$ be a nowhere degenerate hermitian matrix satisfying
\begin{align*}
    \det A \geq \sigma>0
\end{align*}
and $ B= (B_{ij}) \in \C^{\infty} (U, \CC^{n\times n} )$, a symmetric matrix. For each $\bff = (\bff^1,\ldots,\bff^n) \in \C^{\infty}(U \times \partial D, \CC^n)$, there exists a unique solution, $(\hat{z}^1,\ldots, \hat{z}^n; \hat{\xi}^1, \ldots, \hat{\xi}^n)$, with each $\hat{z}^i$, $\hat{\xi}^i\in \C^{\infty} (U \times D, \CC)$, to the follow equation
 \begin{align}\label{holodisclinearpde-lemma}
    \begin{split}
        &\overline{\partial}_\tau \hat{z}^i (w, \tau) = \overline{\partial}_\tau \hat{\xi}^i (w, \tau) = 0, \hspace{5.2cm} \text{ in }U \times D,\\
         &\hat{\xi}^i(w,\tau) - A_{i\bar{j}}(w) \overline{\hat{z}}^j(w, \tau)- B_{ij}(w) \hat{z}^j(w, \tau)=\bff^i(w,\tau),
         \hspace{0.6cm}  \text{ in } U \times \partial D,\\
         &\hat{z}^i(w, -i) = 0,  \hspace{7.5cm}  \ w  \in U,
    \end{split}
\end{align}
Moreover, for each $k \geq 2, \ k\in \ZZ$ and $0<\beta< \alpha <1$, there is a uniform constant $C$ only depending on $n$, $k$, $\alpha$, $\sigma^{-1},\| A\| _{k,\beta; U}$, $\| B\| _{k, \beta; U}$ such that
\begin{align} \label{estiperturbholodiscs}
   \| (\hat{z},\hat{\xi})\| _{k,\beta; U \times D} \leq C (\alpha-\beta)^{-1} \| \bff\| _{k, \alpha; U \times \partial D}.
\end{align}
\end{lem}

\begin{proof}
The proof relies on the decoupling transformation introduced in \eqref{defg1g2}. As established in the previous subsection, this change of variables reduces the coupled linearized system \eqref{holodisclinearpde-lemma} to a pair of independent Riemann-Hilbert problems \eqref{equstandardfbpros}. Applying the scalar estimate \eqref{estisolutionbabyversionfreebdry} to each component of the decoupled variables $\hat{g}_1$ and $\hat{g}_2$, we obtain the following estimates:  
\begin{align} \label{estig1g2D}
    \| (\hat{g}_1,\hat{g}_2)\| _{k,\beta; U \times D} \leq C(\alpha-\beta)^{-1}\| \bff\| _{k,\alpha; U \times \partial D}
\end{align}
where $C$ is a uniform constant independent of $w$ and depending only on $n, \ k, \ \text{and } \alpha$.
Then, the solution, $(\hat{\xi}(w, \tau), \hat{z}(w, \tau))$, of (\ref{holodisclinearpde-lemma}) is given by
\begin{align} \label{inverserepzxibyg}
\begin{split}
    \hat{z} &= \frac{1}{2} (\bar{A})^{-1} (\hat{g}_2- \hat{g}_1),\\
    \hat{\xi} &= \frac{1}{2} (\hat{g}_1 + \hat{g}_2) + \frac{1}{2} B (\bar{A})^{-1} (\hat{g}_2- \hat{g}_1).
\end{split}
\end{align}
By inverse transformation \eqref{inverserepzxibyg}, along with \eqref{estig1g2D}, we complete the proof of (\ref{estiperturbholodiscs}). 

Finally, the uniqueness of $(\hat{z}, \hat{\xi})$ follows directly from the uniqueness of the scalar Riemann-Hilbert problem under the fixed-point normalization condition. \end{proof}

\begin{rem}\label{remregg1g2}
    The estimates in \eqref{estig1g2D} can be stated more precisely by distinguishing the regularity in the $U$-direction from that in the $D$-direction. By applying the estimate for the Hilbert transform in the $\partial D$ direction (see \ref{esticircle}), and together with solving the Dirichlet problems separately for the real and imaginary parts, we obtain the following estimates for $\hat{g}_{1,2}$ in the $ D$ direction:
\begin{align}\label{estiRHbdaryporblemD}
\| \hat{g}_1(w, \cdot)\| _{k,\alpha; D}, \quad \| \hat{g}_2(w, \cdot)\| _{k,\alpha; D} \leq C \| \bff(w, \cdot)\| _{k,\alpha; \partial D}, \qquad \text{for each } w \in U,
\end{align}
where the constant $C$ is a uniform constant only relying on $n$, $k$ and $\alpha$. By Lemma \ref{lemhilbtransesti}, for any $w, w' \in U$ with $|w - w'| \leq t$, the regularity of $\hat{g}_{1,2}$ in the $U$-direction can be described as follows:
\begin{align}\label{estiRHbdaryporblemN}
\sup_{|w - w'| \leq t,\ \tau \in D} \left| \nabla_0^{\mathbf{k}} \hat{g}_i(w', \tau) - \nabla_0^{\mathbf{k}} \hat{g}_i(w, \tau) \right| \leq C t^\alpha (1+|\log t|) \cdot\| \bff\| _{k,\alpha; U \times \partial D}, \qquad i=1,2.
\end{align}
where $\nabla_0^{\mathbf{k}}$ denotes a differential operator with multi-indices $\mathbf{k}$ of order $|\mathbf{k}| = k$ with respect to the coordinates on $U$, and $C$ is a uniform constant. One can easily observe that the estimates \eqref{estiRHbdaryporblemD} and \eqref{estiRHbdaryporblemN} together imply \eqref{estig1g2D}. Moreover, the blow-up rate of the coefficient in \eqref{estig1g2D} as $\beta \to \alpha$ is sharp as
\[
\sup_{0\leq t\leq 1} t^{\alpha-\beta} |\log t| = C (\alpha -\beta)^{-1}.
\]
In the following subsections, we study the regularity of the linearized problem near the trivial data. We will show that if a family of elliptic systems is solved in the $D$ direction—exhibiting a loss of regularity in the $U$ direction with the same behavior as \eqref{estiRHbdaryporblemN} at $(G_0, \partial \Psi_0)$—then, the general linearized problem in a neighborhood of $(G_0, 0)$ satisfies a comparable blow-up rate. Specifically, we will show that the estimate for the linearized operator involves a constant growing as $(\alpha - \beta)^{-1}$ as $\beta \to \alpha$. 
\end{rem}

\subsection{BMO Norms and the Loss of Regularity in H\"older Spaces} In the previous subsection, we observed the loss of regularity in H\"older norms in the parameter direction when solving the linear elliptic boundary problems. By Lemma \ref{lemhilbtransesti}, The critical issue is the following: given a function $f \in \C^{k,\alpha} ( U \times \partial D)$, 
let \(\tilde f:=\Hilb_\tau f\).
For a multi-index \(\mathbf k\) with \(|\mathbf k|=k\), we define the Hölder quotient in the \(U\)-direction by
\begin{align}\label{defckalphainN}
h_{\tilde f}^{\mathbf k,\alpha}(\tau;w,w')
:=
\frac{
\nabla_0^{\mathbf k}\tilde f(w',\tau)
-
\nabla_0^{\mathbf k}\tilde f(w,\tau)
}
{|w'-w|^\alpha}.
\end{align}
Then \(\left| h_{\tilde f}^{\mathbf k,\alpha} \right|\) may grow like
\(
\big|\log |w'-w|\big|
\)
as \(w'\to w\), and hence is not uniformly bounded in \(L^\infty\).

In this subsection, we reinterpret this phenomenon from the viewpoint of \emph{bounded mean oscillation} (BMO) spaces. There is an extensive literature on this subject; here we only mention a few classical references that may help the reader gain further insight \cite{David1984, grafakos2009modern, JNthm1961, SteinbookHA}. In this subsection, we avoid developing subtle technical details on the theorems concerning BMO spaces. Instead, we present only the essential definitions and lemmas that will be used in establishing the perturbation theorem for the linear elliptic boundary problem discussed in the next subsection.

\begin{define}
    Let $u$ be a locally integrable function defined in $\RR^n$. The bounded mean oscillation (BMO) seminorm of $u$ is defined to be
    \begin{align*}
        \| u\| _{\BMO} = \sup_{Q\subseteq \RR^n} \fint_Q |u-u_Q|,
    \end{align*}
    where the supremum is taken over all cubes $Q \subseteq \RR^n$, and $u_Q$ is the average of $u$ in $Q$. We say $u$ is a bounded mean oscillation (BMO) function if $\| u\| _{\BMO} < +\infty$.
\end{define}

\begin{define} \label{defCZO}
    Let $T$ be a linear operator $L^2(\RR^n) \rightarrow L^2(\RR^n)$ given by an integral of singular kernel as follows:
    \begin{align*}
        Tf(x) = \operatorname{p.v.}\!\int K(x, y) f(y) dy.
    \end{align*}
    We say $T$ is Calder\'on-Zygmund operator (CZO) if $T$ is bounded in $L^2$, and the kernel $K(x, y)$ is defined away from the diagonal $\{x=y\}$ satisfying:
    \begin{itemize}
        \item Decay control at singularities: 
        \[|K(x,y)| \leq  \frac{C}{|x-y|^n}.  \]
        \item H\"older regularity: for $0<\alpha\leq 1$,
        \[ |K(x,y) - K(x', y)| \leq C \frac{|x-x'|^\alpha}{|x-y|^{n+\alpha}}, \qquad \text{if }  |x-x'| < \frac{1}{2} |x-y|, \]
        and 
        \[|K(x,y) - K(x, y')| \leq C \frac{|y-y'|^\alpha}{|x-y|^{n+\alpha}}, \qquad \text{if }  |y-y'| < \frac{1}{2} |x-y|. \]
    \end{itemize}
    Furthermore, we say a CZO $T$ is cancellative if $T(1) =0$.
\end{define}
A more general definition of Calder\'on–Zygmund operators can be formulated in the framework of the $T(1)$ theorem (see, e.g., \cite{David1984}). Since our discussion centers on the Hilbert transform, it is enough for us to adopt the simpler $L^2$–bounded setting given above.
Recall that the Hilbert transform defined on the unit circle can be viewed roughly as a convolution with a singular kernel as 
\( \displaystyle
\Hilb f (\tau) = \frac{1}{2\pi}\operatorname{p.v.}\!\int  \cot \big(\frac{\tau-\theta}{2}\big)  f(\theta) d\theta
\). 
It is easy to check that the Hilbert transform is cancellative CZO with $\|\Hilb f\|_{L^2(\partial D)}=\|f-f_{\partial D}\|_{L^2(\partial D)}\leq\|f\|_{L^2(\partial D)}$.

It is known that CZOs are bounded maps from $L^\infty$ to BMO spaces; precisely, if $T$ is a CZO, we have
\begin{align} \label{bddCZOBMOLinfty}
\| Tf\| _{\BMO} \leq C\| f\| _{L^\infty},
\end{align}
where $C$ is a constant depending on $L^2$ norm of $T$, the decay control and H\"older regularity in Definition \ref{defCZO}. Note that $h_{\tilde{f}}^{\mathbf{k}, \alpha}(\tau; w, w')$ in \ref{defckalphainN} is the Hilbert transform of 
\[
h_{{f}}^{\mathbf{k},\alpha} (\tau; w,w') = \frac{\nabla_0^{\mathbf{k}}f(w', \tau) -  \nabla_0^{\mathbf{k}}f (w, \tau) } {|w'-w|^\alpha},
\]
in $\tau$. Consequently, $h_{\tilde{f}}^{\mathbf{k},\alpha}(\tau; w, w')$ is uniformly bounded in the BMO norm. 

The key result that we will apply later is the John–Nirenberg theorem, which provides a precise description of how BMO bounds control the size of a function.

\begin{lem}[John--Nirenberg Theorem \cite{JNthm1961}] \label{lemJNT}
Let $f \in \BMO(\RR^n)$, and let $Q \subset \RR^n$ be any cube. Then there exist constants $c_1, c_2 > 0$, depending only on the dimension $n$, such that
\[
 \big| \{ x \in Q : |f(x) - f_Q| > \lambda \} \big|
\;\leq\; c_1 \exp\!\left(- \frac{c_2 \lambda}{\|f\|_{\BMO}}\right) |Q|,
\qquad \forall \lambda > 0,
\]
where $f_Q $ is the average of $f$ over $Q$. 
\end{lem}

The following is a quick corollary from John-Nirenberg, which implies the alternate characterisation of BMO norm:

\begin{cor}
    For any $1\leq p< \infty$, and any locally integrable function $f$ we have 
    \begin{align*}
        \| f\| _{\BMO} \sim  \sup_{Q \subseteq \RR^n}\bigg( \fint_{Q} \big|f-  f_Q\big|^p\bigg)^{1/p} \sim \sup_{Q \subseteq \RR^n} \inf_{c} \bigg( \fint_{Q} \big|f-  c\big|^p\bigg)^{1/p}
    \end{align*}
\end{cor}

The next lemma is a generalization of \eqref{bddCZOBMOLinfty}, which is useful in the iterative process.

\begin{lem} \label{lembddBMOitera}
    Let $T$ be a cancellative CZO in $\RR^n$. Then, $T$ is a bounded operator in $\BMO(\RR^n)$ with
    \begin{align*}
        \| Tf\| _{\BMO} \leq C \| f\| _{\BMO},
    \end{align*}
    where $C$ is a uniform constant depending on $n$, the $L^2$ operator norm of $T$, the H\"older index $\alpha$ and the constant in H\"older regularity estimate of Definition \ref{defCZO}.
\end{lem}

\begin{proof}
    Given a cube $Q \subseteq \RR^n$, we estimate the mean oscillation of $Tf$ on $Q$: 
    \begin{align*}
        \frac{1}{|Q|}\int_Q |Tf-c_Q|,
    \end{align*}
    for a fixed constant $c_Q \in \CC$. Let $mQ$ denote the dilation of $Q$ about its center, $x_Q$, by a factor of $m$. To prove the estimate, we split the function $Tf$ to be the near and the far part:
    \[
    Tf = T \big( (f-f_{2Q})\mathbf{1}_{2Q} \big) + T \big( (f-f_{2Q}) \mathbf{1}_{(2Q)^c} \big),
    \]
    where $f_{2Q}$ is the average of $f$ in $2Q$. 
    For simplicity, we write $\tilde{f}_{\text{near}} := T \big( (f-f_{2Q})\mathbf{1}_{2Q} \big)$ and $\tilde{f}_{\text{far}} := T \big( (f-f_{2Q})\mathbf{1}_{(2Q)^c} \big)$. If we take the constant $c_Q$ to be, \(
    c_Q = \tilde{f}_{\text{far}} (x_Q),
    \)
    we have
    \begin{align*}
        \frac{1}{|Q|}\int_Q |Tf-c_Q| \leq  \fint_Q |\tilde{f}_{\text{near}}| + \fint_{Q} |\tilde{f}_{\text{far}} - \tilde{f}_{\text{far}} (x_Q)|.
    \end{align*}
    For the near part, we have
    \begin{align*}
        \fint_Q |\tilde{f}_{\text{near}}| 
        &\leq \frac{1}{|Q|^{1/2}} \| \tilde{f}_{\text{near}}\| _{L^2; ~ Q} \\
        &\leq \frac{1}{|Q|^{1/2}} \| T\| _{L^2 \rightarrow L^2} \| f-f_{2Q}\| _{L^2; ~ 2Q} \\
        &\leq C_n\| T\| _{L^2 \rightarrow L^2} \| f\| _{\BMO}.
    \end{align*}
    For the far part, we have 
    \begin{align*}
        \fint_Q |\tilde{f}_{\text{far}} - \tilde{f}_{\text{far}}(x_Q)| &\leq \frac{1}{|Q|} \int_{Q} dx \int_{(2Q)^c} |K(x,y) - K(x_Q, y)|\cdot |f(y) - f_{2Q}| dy \\
        & \leq \frac{C_h}{|Q|} \int_Q |x-x_Q|^\alpha dx \int_{(2Q)^c} \frac{|f(y) - f_{2Q}|}{|x_Q-y|^{n+\alpha}} dy
        \\
        & \leq C_h |\diam Q|^\alpha \int_{(2Q)^c} \frac{|f(y) - f_{2Q}|}{|x_Q-y|^{n+\alpha}} dy
    \end{align*}
    where $C_h$ represents the constant in H\"older regularity of Definition \ref{defCZO}. The integral region $(2Q)^c$ can decompose to $(2Q)^c =\cup_{k} A_k$, where $A_k = 2^{k+1} Q \backslash 2^k Q$. Then, 
    \begin{align*}
        \int_{(2Q)^c} \frac{|\diam Q|^\alpha}{|x_Q-y|^{n+\alpha}} |f(y) - f_{2Q}| dy
        &\leq \sum_{k} \frac{1}{2^{k\alpha}} \frac{C_n}{|2^{k+1} Q|} \int_{A_k} |f(y) - f_{2Q}| dy
        \\
        &\leq  \sum_{k\geq 1} \frac{C_n}{2^{k\alpha}}  \bigg( \fint_{2^{k+1} Q} |f(y)- f_{2^{k+1} Q}| dy + |f_{2^{k+1} Q} - f_{2Q}| \bigg)
        \\
        &\leq \sum_{k \geq 1} \frac{C_n}{2^{k\alpha}} \Big(\| f\| _{\BMO} +  |f_{2^{k+1} Q} - f_{2Q}|\Big).
    \end{align*} 
    It suffices to estimate $|f_{2^{k+1} Q} - f_{2Q}|$. Observe that
    \begin{align*}
        |f_{2^{k+1} Q} - f_{2^k Q}| &\leq \fint_{2^k Q} |f(x)-f_{2^{k+1} Q}| dx
        \\
        &\leq 2^n \fint_{2^{k+1}Q} |f(x)- f_{2^{k+1}Q} | dx \leq C_n\| f\| _{\BMO},
    \end{align*}
    and 
    \begin{align} \label{estitelescope}
        |f_{2^{k+1} Q} - f_{2 Q}| \leq\sum_{j=1}^k |f_{2^{j+1 } Q} - f_{2^j Q}| \leq C_nk \| f\| _{\BMO}.
    \end{align}
    Inserting the above inequality back into the far part estimate, we have
    \begin{align*}
        \fint_Q |\tilde{f}_{\text{far}} - \tilde{f}_{\text{far}}(x_Q)| \leq C_h C_n \sum_{k\geq 1} \frac{k+1}{2^{k\alpha}} \| f\| _{\BMO} \leq C \| f\| _{\BMO},
    \end{align*}
    where $C$ is a constant depending only on $n,~\alpha, ~ C_h$. Combining the estimates for both the near part and the far part, we complete the proof.
\end{proof}

\begin{lem}\label{lemBMOproduc}
    Let $f\in \C^{0,\gamma} (\RR^n)$ and $g \in \BMO(\RR^n)$. Suppose that $f$ and $g$ has a compact support in $K \subseteq \RR^n$ and
    \[
    \left|\fint_K g\right| \leq C_0.
    \]
    Then, we have
    \[
    \| f\cdot g\| _{\BMO} \leq C \| f\| _{0,\gamma}\cdot (\| g\| _{\BMO}+ C_0),
    \]
    where $C$ is a constant depending on $n$, $\gamma$ and the size of $K$. 
\end{lem}
\begin{proof}
    For each cube $Q \subseteq K$, we test the following mean oscillation 
    \[
    I := \fint_Q |fg- f_Qg_Q|.
    \]
    We split the integrand as follows: 
    \[
    |fg-f_Qg_Q| \leq |f-f_Q| |g-g_Q| + |f_Q| |g-g_Q| + |g_Q|\,|f-f_Q|.
    \]
    Then, we have
    \[
    I \leq 3\| f\| _{L^\infty} \| g\| _{\BMO} + |g_Q| \fint_{Q} |f-f_Q|.
    \]
    It suffices to deal with the second term on the right-hand side of the above inequality. Let $Q$ satisfies $Q \subseteq K \subseteq 2^k Q := Q'$. Now, if we fix the large cube $Q'$, by the ``telescope'' technique in \eqref{estitelescope}, we can control the growth of $g_Q$:
    \begin{align*}
        |g_Q| \leq C \left(1+\log \frac{|Q'|}{|Q|}\right) (\| g\| _{\BMO}+ C_0).
    \end{align*}
    By the H\"older continuity of $f$, we have
    \[
    \fint_Q |f- f_Q| \leq C  |Q|^{\gamma/n} [f]_{0,\gamma}
    \]
    Observe that $|Q|^{\gamma/n}\log |Q|^{-1} \leq C \gamma^{-1}$; hence, we complete the proof.
\end{proof}

\begin{rem}
    By John-Nirenberg Theorem, functions in BMO exhibit substantially better behavior than those in $L^p$. In particular, for any domain $B \subseteq \RR^n$ and $f \in L^p (B)$, $|\{x\in B: |f(x)|> \lambda\}| \leq c (\lambda/\| f\| _{L^p})^{-p}$. The rapid decay of the measure of super-level sets as the threshold increases guarantees that any loss of Hölder regularity remains controlled for the elliptic boundary problem, provided that $h_{b}^{\mathbf{k}, \alpha}(\tau; w, w')$ is bounded in BMO norms, where $b$ is the boundary data.
\end{rem}

The following lemma serves as the cornerstone in establishing the regularity theorem for a family of linear elliptic boundary problems, which explains the relation between the boundedness of BMO norm and the loss of regularity in H\"older spaces.

\begin{lem}[BMO--Hölder interpolation]
\label{lemBMOHolderInterpolation}
Let \(0<\beta<\alpha<1\), {let \(U\) be one of the fixed holomorphic coordinate balls,} and let \(f\) be defined on
\(U\times \partial D\). Assume
\(
\sup_{w\in U}\|f(w,\cdot)\|_{k,\alpha;\partial D}<\infty
\)
and
\[
\sup_{w\neq w'}
\left(
\|h_f^{\mathbf k,\alpha}(\cdot;w,w')\|_{\BMO; \partial D}
+
\|h_f^{\mathbf k,\alpha}(\cdot;w,w')\|_{L^2(\partial D)}
\right)<\infty.
\]
{for every \(|\mathbf k|=k\).} Then
\[
{\begin{aligned}
\|f\|_{k,\beta;U\times\partial D}
\le
\frac{C}{\alpha-\beta}
\bigg[
&\sup_{w\in U}\|f(w,\cdot)\|_{k,\alpha;\partial D} \\
&+
\max_{|\mathbf k|=k}\sup_{w\neq w'}
\left(
\|h_f^{\mathbf k,\alpha}(\cdot;w,w')\|_{\BMO; \partial D}
+
\|h_f^{\mathbf k,\alpha}(\cdot;w,w')\|_{L^2(\partial D)}
\right)
\bigg],
\end{aligned}}
\]
where $C$ is a uniform constant.
\end{lem}

\begin{proof}
We first prove for the case $k =0$.
Let \( g \) be an \(\alpha\)-H\"older continuous function on one fixed fiber $\partial D$ with \(\|g\|_{0,\alpha; \partial D} < +\infty \) and hence,  \(\|g\|_{\BMO; \partial D} + \|g\|_{L^2(\partial D)}< +\infty\). For any arc \(I\subset \partial D\),
\[
\fint_I |g-g_I|
\le
\|g\|_{\BMO},
\]
By Campanato criterion for \(\alpha\)-H\"older norm, the $\alpha$-H\"older semi-norm is equivalent to \(\alpha\)-Campanato norm: \( \displaystyle \sup_{I \subseteq \partial D} |I|^{-\alpha} \fint_I |g-g_I| \). By Lemma \ref{lemCampanatoExplicit}, we have
\[
|I|^{-\alpha}\fint_I |g-g_I|
\le
C [g]_{0,\alpha} .
\]
Therefore, we have
\( \displaystyle
\fint_I |g-g_I|
\le
C\min\{\|g\|_{\BMO},[g]_{0,\alpha; \partial D}|I|^\alpha\}.
\)
The key observation of interpolation is that for each $0<\beta< \alpha< 1$ such that  
\[
\fint_I |g-g_I|
\le
C \|g\|_{\BMO; \partial D}^{\beta/\alpha}\|g\|_{0,\alpha; \partial D}^{1-\beta/\alpha}|I|^{\alpha-\beta}.
\]
Again, by Lemma \ref{lemCampanatoExplicit}, we have
\[
[g]_{0,{\alpha-\beta};\partial D}
\le \frac{C}{\alpha-\beta}\sup_{I\subseteq\partial D}|I|^{-(\alpha-\beta)}\fint_I |g-g_I| \le
\frac{C}{\alpha-\beta}
\|g\|_{\BMO; \partial D}^{\beta/\alpha}\|g\|_{0,\alpha; \partial D}^{1-\beta/\alpha}.
\]
For $L^\infty$ norm, we can easily derive the following:
\begin{align}
\begin{split}
\|g\|_{L^\infty (\partial D)} &\le \left| \fint_{\partial D} g\right| + (2\pi)^{\alpha-\beta} [g]_{0,\alpha-\beta; \partial D} 
\\
&\le \frac{C}{\alpha-\beta}
\left(\|g\|_{\BMO; \partial D} + \|g\|_{L^2(\partial D)} \right)^{\beta/\alpha} \|g\|_{0,\alpha; \partial D}^{1-\beta/\alpha}.
\end{split}
\end{align}

Now, we consider the function defined on a family of circles, $f$ that satisfies the conditions of the lemma. 
Fix \(w\neq w'\), 
apply the preceding estimate to
\[
h_f^\alpha(\tau;w,w')
=
\frac{f(w',\tau)-f(w,\tau)}{|w-w'|^\alpha}.
\]
If we set $|w-w'| = \delta$, then
\[
\|h_f^\alpha(\tau;w,w')\|_{0,\alpha;\partial D}
\le
C \sup_{w \in U} \|f\|_{0,\alpha; \partial D} \delta^{-\alpha}.
\]
Therefore, the modulus of continuity in $U$-direction satisfies
\begin{align*}
\|&f(w, \tau) -  f(w', \tau)\|_{L^\infty(\partial D)}
\\
& \le
\frac{C}{\alpha-\beta} 
\left(
\|h_f^{\alpha}(\cdot;w,w')\|_{\BMO; \partial D}
+
\|h_f^{\alpha}(\cdot;w,w')\|_{L^2(\partial D)}
\right) ^{\beta/\alpha}
\left( \sup_{w\in U}\|f\|_{0,\alpha; \partial D}\delta^{-\alpha}\right)^{1-\beta/\alpha} \delta^\alpha
\\
& \le \frac{C}{\alpha-\beta} \delta^\beta \left( \sup_{ U}\|f\|_{0,\alpha;\partial D} +
\sup_{w\neq w'}
\left(
\|h_f^{\alpha}(\cdot;w,w')\|_{\BMO; \partial D}
+
\|h_f^{\alpha}(\cdot;w,w')\|_{L^2(\partial D)}
\right) \right),
\end{align*}
where the last step we use the Young's inequality. We obtain the $\beta$-Hölder estimate in $w$-direction. 
The Hölder estimate in $\tau$-direction follows
directly from the fiberwise bound:
\[
\sup_{w\in U}\|f(w,\cdot)\|_{0,\beta;\partial D}
\le
C \sup_{U} \|f(w,\cdot)\|_{0,\alpha;\partial D}.
\]
Therefore, 
\[
\|f\|_{0,\beta;U\times\partial D}
 \le \frac{C}{\alpha-\beta}\left( \sup_{ U}\|f\|_{0,\alpha;\partial D} +
\sup_{w\neq w'}
\left(
\|h_f^{\alpha}(\cdot;w,w')\|_{\BMO; \partial D}
+
\|h_f^{\alpha}(\cdot;w,w')\|_{L^2(\partial D)}
\right) \right).
\]

For \(k\geq1\), apply the scalar interpolation estimate above to a spatial finite difference of order \(k+1\), with step length \(\delta\) and all its points in \(U\). Repeated integration of the \(k\)-th spatial derivatives bounds its \(\BMO\) plus \(L^2\) norm by \(C\delta^{k+\alpha}\) times the second term in the bracket in the statement; its circle \(C^{0,\alpha}\) norm is bounded by \(C\) times the first term. Replace \(\beta\) in that scalar estimate by
\[
\frac{\alpha(k+\beta)}{k+\alpha}.
\]
Young's inequality then bounds its \(L^\infty\) norm by \(C\delta^{k+\beta}/(\alpha-\beta)\) times the stated bracket. By the finite-difference characterization, we have spatial \(C^{k,\beta}\) bounds. Together with the fibrewise bounds, these imply the joint \(C^{k,\beta}\) estimate, including mixed derivatives, by standard H\"older interpolation. 
\end{proof}

\begin{rem}
    In the interpolation lemma, the \(L^2\)-term can be replaced by a bound
on the average of \(h_f^{\mathbf k,\alpha}\). Indeed, on \(\partial D\),
the John--Nirenberg inequality implies
\[
\|g\|_{L^2(\partial D)}
\le
C\left(
\|g\|_{\BMO(\partial D)}
+
\left| \fint_{\partial D} g \right|
\right).
\]
Conversely, by H\"older's inequality,
\(
\left| \fint_{\partial D} g \right|
\le
C\|g\|_{L^2(\partial D)}.
\)
Hence the two quantities
\(
\|g\|_{\BMO(\partial D)}+\|g\|_{L^2(\partial D)}
\)
and
\(
\|g\|_{\BMO(\partial D)}
+
\left|\fint_{\partial D}g\right|
\)
are equivalent up to a uniform constant.

We state the lemma using the \(L^2\)-norm because the Hilbert transform is bounded on \(L^2(\partial D)\). Thus, the \(L^2\)-form of the estimate is more convenient for the
iteration in the next subsection.
\end{rem}

\begin{cor}
\label{corHilbertBMOHolder}
Let $b$ be the boundary function in the family of elliptic problems in \eqref{dbardisk} defined on $U \times \partial D$,  and 
let \(\widetilde b (w, \tau)=\Hilb_\tau b (w, \tau)\), where the Hilbert transform is taken in the
\(\tau\)-variable. Then for every \(0<\beta<\alpha\),
{\begin{align*}
 \|\tilde{b}\|_{k,\beta;U\times\partial D} \le \frac{C}{\alpha-\beta}
\bigg[
&\sup_{w\in U}\|b(w,\cdot)\|_{k,\alpha;\partial D} \\
&+
\max_{|\mathbf k|=k}\sup_{w\neq w'}
\left(
\|h_b^{\mathbf k,\alpha}(\cdot;w,w')\|_{\BMO; \partial D}
+
\|h_b^{\mathbf k,\alpha}(\cdot;w,w')\|_{L^2(\partial D)}
\right)
\bigg],
\end{align*}}
\end{cor}

\begin{proof} The proof follows directly from Lemma \ref{lembddBMOitera} and Lemma \ref{lemBMOHolderInterpolation}.
\end{proof}

\begin{rem} {For the following direct estimate, assume additionally that \(b\in C^{k,\alpha}(U\times\partial D)\).} We can estimate $h_{\tilde{b}}^{\mathbf{k}, \alpha} (\tau; w, w')$ {as} follows:
    If we denote $\delta = |w-w'|$, by \eqref{estipara}, we have
    \begin{align*}
        |h_{\tilde{b}}^{\mathbf{k}, \alpha} (\tau; w, w')|  \leq & \ C \delta^{-\alpha} \bigg\{ \int_{\tau-\delta}^{\tau + \delta} \frac{|\nabla_0^{\mathbf{k}}b(w, \theta) - \nabla_0^{\mathbf{k}} b(w, \tau) |}{|\theta -\tau| }  d\theta  + \int_{\tau-\delta}^{\tau + \delta} \frac{|\nabla_0^{\mathbf{k}}b(w', \theta) - \nabla_0^{\mathbf{k}} b(w', \tau) |}{|\theta -\tau| }  d\theta \bigg\}
         \\
         & \qquad +  \int_{\delta<|\theta-\tau|< \pi} |h_b^{\mathbf{k}, \alpha}(\theta; w, w')| \left| \cot \big(\frac{\theta-\tau}{2} \big) \right|d\theta.
    \end{align*}
{The} first two {integrals,} {including} the {factor} {\(\delta^{-\alpha}\),} are bounded by
{\[
\frac{C}{\alpha}\sup_{w\in U}
\|\nabla_0^{\mathbf k}b(w,\cdot)\|_{0,\alpha;\partial D}
\leq \frac{C}{\alpha}\|b\|_{k,\alpha;U\times\partial D}.
\]}
It {remains} to {estimate} the third term. By Lemma \ref{lemJNT}, for {\(\lambda>0\)},
    \begin{align*}
        |\{\theta\in \partial D, ~ \left|h_b^{\mathbf{k}, \alpha}(\theta; w, w') - \big( h_b^{\mathbf{k}, \alpha}(\cdot; w, w')\big)_{\partial D} \right| > \lambda \}| \leq c_1 \exp \big(-c_2 \lambda / C_0\big),
    \end{align*}
{where} {\(\big(h_b^{\mathbf k,\alpha}(\cdot;w,w')\big)_{\partial D}\)} is the {circle average, and \(C_0\) bounds the BMO and \(L^2\) norms} of {these} {quotients.} The {mean} {contributes} {at} {most} {\(CC_0(1+|\log\delta|)\).} {Applying} the layer-cake {representation} {to} the {centered} {quotient, we have}
    \begin{align*}
        (\text{The third integral})
        &\leq C \sum_{\pm}\int_\delta^\pi \frac{| h_b^{\mathbf{k}, \alpha}(\tau \pm t; w, w') |}{t} dt\\
        & \leq {CC_0(1+|\log\delta|)+} C \sum_{n = 1}^\infty  \log \Big(\frac{\delta + c_1 e^{-c_2 n /C_0}}{\delta}\Big) 
        \\
        &\leq {CC_0(1+|\log\delta|)+}
        C\sum_{1\leq n \leq \frac{C_0}{c_2}|\log 
    \delta|} 
     |\log (c_1 \delta^{-1})|
     + C \sum_{m=1}^\infty  \log (1+ c_1 e^{-c_2m / C_0})
    \\
    &\leq  C(|\log \delta|^2 +1).
    \end{align*}
    By the same observation as in Remark \ref{remregg1g2},
    \begin{align*}
        \sup_{|w - w'| \leq \delta,\ \tau \in D} \left| \nabla_0^{\mathbf{k}} \tilde{b}(w', \tau) - \nabla_0^{\mathbf{k}} \tilde{b}(w, \tau) \right| \leq C \delta^\alpha (|\log \delta|^2 + 1).
    \end{align*}
    If we estimate the $\C^{k,\beta}$ norm of $\tilde{b}$ in this way, the blow-up rate of the coefficient as $\beta \rightarrow \alpha$ is given by
    \begin{align*}
        \sup_{0< \delta\leq 1} \delta^{\alpha-\beta} |\log \delta|^2 \sim (\alpha-\beta)^{-2}.
    \end{align*}
\end{rem}

\subsection{Linearized problems near the trivial data} In this subsection, we study the linearized problem by considering a small perturbation of the trivial data $(G_0, \Lambda_0)$, where locally $G_0 = (z_0, \xi_0) = (w, \partial \rho (w))$ and $\Lambda_0$ is determined by the graph of $\partial \rho$. We consider a perturbed pair $(G, \Lambda_{\Psi})$, where the map is written in components as $G (w, \tau) = (z(w, \tau), \xi(w, \tau))$.
Defining {\(\widetilde A_{i\bar j}(w,\tau)=(\partial_i\bar\partial_j(\rho+\Psi))(z(w,\tau),\tau)\)} and {\(\widetilde B_{ij}(w,\tau)=(\partial_i\partial_j(\rho+\Psi))(z(w,\tau),\tau)\)}, the linearized problem at {\((G,\Lambda_\Psi)\)} is {written} in {a} {fixed} local {trivialization} {as} {follows:}
\begin{align}\label{holodisclinearpde'}
    \begin{split}
        &\overline{\partial}_\tau \hat{z}^i (w, \tau) = \overline{\partial}_\tau \hat{\xi}^i (w, \tau) = 0, \hspace{5.8cm} \text{ in }U \times D,\\
         &\hat{\xi}^i(w,\tau) - \widetilde{A}_{i\bar{j}}(w, \tau) \overline{\hat{z}}^j(w, \tau)- \widetilde{B}_{ij}(w, \tau) \hat{z}^j(w, \tau)=\bff^i(w,\tau),
         \hspace{.75cm}  \text{ in } U \times \partial D,\\
         &\hat{z}^i(w, -i) = 0,  \hspace{8.2cm}  \ w  \in U,
    \end{split}
\end{align}
where $\bff= (\bff^1, \ldots, \bff^n) \in \C^{\infty}(U \times \partial D, \CC^n)$.
As in the previous subsection, recall that $A_{i\bar{j}}(w) = \partial_i \overline{\partial}_j \rho$ and $B_{ij} (w) = \partial_i \partial_j \rho$ are the corresponding data at the trivial data, $(G_0, \Lambda_0)$. The small perturbation condition of $(G, \Lambda_\Psi)$ relative to $(G_0, \Lambda_0)$ implies the following estimates:
\begin{align} \label{constantsmallperturb}
\begin{split}
    \| \widetilde{A}_{i\bar{j}}-A_{i \bar{j}}\| _{{k, \alpha}, U \times \partial D} \leq \delta,\\
    \| \widetilde{B}_{ij}-B_{i {j}}\| _{{k, \alpha}, U \times \partial D} \leq \delta.
\end{split}
\end{align}
where $\delta$ is a sufficiently small constant to be determined later.

Lemma \ref{IFTlinearpart} establishes {an} inverse for the linearized operator {at} the trivial {data,} {with} a loss of {Hölder} regularity. {Write} {\(T'_0=D_{\mathcal F}T|_{(G_0,\Lambda_0)}\)} and {\(T'=D_{\mathcal F}T|_{(G,\Lambda_\Psi)}\).} {We} {perform} the {estimates} {in} {a} {fixed} local {trivialization,} where {\(T'=T'_0-P_\delta\)}. {The} formal series
$$(T')^{-1} = \sum_{n=0}^\infty \left( (T'_0)^{-1} P_\delta \right)^n (T'_0)^{-1}$$
is not well-defined. This failure arises from the loss of regularity in $(T'_0)^{-1}$: the operator maps $C^{k, \alpha}$ to $C^{k, \beta}$ (with $\beta < \alpha$), which prevents the convergence of the iterated composition $(T'_0)^{-1} P_\delta$.

The purpose of the BMO estimates from the previous subsection is to control this loss more precisely. Instead of iterating solely in Hölder spaces, we track the difference quotients in the parameter direction as BMO functions. This refined control allows the iteration to close, at the cost of a quantified loss in the final estimate, reflected by the factor
\(
(\alpha-\beta)^{-1}.
\)

\begin{lem} \label{lemlinearproblemperturb}
    Let $A$ and $B$ be the coefficient matrices associated with the trivial data $(G_0, \Lambda_0)$, with components defined locally by $A_{i\bar{j}} = \partial_i \overline{\partial}_j \rho$ and $B_{ij} = \partial_i \partial_j \rho$. We assume $A$ is uniformly nondegenerate with 
    $$\det A \geq \sigma > 0.$$ 
    
    Then, there exists \(\delta_0>0\), depending only on \(n,\ k,\ \alpha,\ \sigma^{-1},\ \|A\|_{{k,\alpha}; U },\ \|B\|_{{k,\alpha}; U},
\) such that the following holds.
    Let $\widetilde{A}$ and $\widetilde{B}$ be the coefficients of the linearized problem at the perturbed data $(G, \Lambda_\Psi)$, defined locally by \(\widetilde{A}_{i \bar{j}}(w, \tau) = (\partial_{i} \overline{\partial}_{j} (\rho +{\Psi}))(z(w, \tau), \tau) \)  and \( \widetilde{B}_{i j} (w, \tau) = (\partial_i \partial_j (\rho +{\Psi})) (z(w, \tau), \tau).\)
    We assume they satisfy the local perturbation condition \eqref{constantsmallperturb} with a sufficiently small constant $\delta$ with $0< \delta < \delta_0$. For each $\bff \in C^{\infty}(U \times \partial D, \mathbb{C}^n)$, there exists a unique global solution $\widehat{G} \in \F'_G$ to the system \eqref{holodisclinearpde'}.

Moreover, for any integer $k \geq 1$ and $0 < \beta < \alpha < 1$, there exists a uniform constant $C$ depending only on $n, \ k,\ \alpha, \ \sigma^{-1}$, and the norms $\|A\|_{k,\alpha; U \times \partial D}, \ \|B\|_{k, \alpha; U \times \partial D}$ such that the local components $(\hat{z}, \hat{\xi})$ of $\widehat{G}$ satisfy
\begin{align} \label{estiperturbholodiscs'}
    \|(\hat{z}, \hat{\xi})\|_{k,\beta; U \times D} \leq C (\alpha-\beta)^{-2} \|\bff\|_{k, \alpha; U \times \partial D}.
\end{align}
\end{lem}

\begin{proof} We first solve the problem locally. In each coordinate chart $U_i \times D$, the system \eqref{holodisclinearpde'} reduces to a set of scalar partial differential equations. Without loss of generality (and to simplify the notation), we will use $U$ to denote a local chart in this proof.

Let $T'$ be the perturbed linearized operator, defined locally by
\[
T'(\hat{z}, \hat{\xi})^i (w, \tau) = \hat{\xi}^i (w, \tau) - \widetilde{A}_{i\bar{j}} (w, \tau) \bar{\hat{z}}^j (w, \tau) - \widetilde{B}_{ij} (w, \tau) \hat{z}^{j}(w, \tau), \qquad (w, \tau) \in U \times \partial D,
\]
for $\widehat{G} = (\hat{z}, \hat{\xi}) \in \F'$. 
We define the perturbation operator $P_\delta = -( T'- T'_0)$ by 
\begin{align*}
\begin{split}
P_{\delta} (\hat{z}, \hat{\xi}) (w, \tau) = \big(&\widetilde{A}_{i\bar{j}} (w, \tau) - A_{i\bar{j}}(w)\big) \bar{\hat{z}}^{j} (w, \tau)\\
&\hspace{0.5cm}+ \big( \widetilde{B}_{ij} (w, \tau)- B_{ij} (w)\big) \hat{z}^j(w, \tau),
\end{split}
\qquad (w, \tau) \in U \times \partial D,
\end{align*}
where both $(A,\widetilde{A})$ and $(B, \widetilde{B})$ satisfy the smallness condition \eqref{constantsmallperturb}. We now introduce the notation for the iterative process. In the following iteration, we use a convention to distinguish the correction terms from the partial sums. Lower indices in parentheses denote the individual correction terms:
\(
\widehat G_{(l)}=(\hat z_{(l)},\hat\xi_{(l)}),
\)
while lower indices without parentheses denote the partial sums, or equivalently the approximations after \(l\) steps:
\(
\widehat G_l=(\hat z_l,\hat\xi_l)
:=
\sum_{\ell=0}^l \widehat G_{(\ell)}.
\)
Thus \(\widehat G_{(l)}\) is the \(l\)-th error correction, whereas \(\widehat G_l\) is the \(l\)-th approximate solution. Upper indices are reserved for local coordinate components. For example,
\(
(\hat z_{(l)}^i, \hat{\xi}_{(l)}^i)
\)
and
\(
(\hat z_{l}^i, \hat{\xi}_{l}^i).
\)
Throughout the proof, \(C_j\) \((j=0,1,2,\ldots)\) denotes a large uniform positive constant depending only on
\(
n,\ k,\ \alpha,\ \sigma^{-1},\ 
\|A\|_{{k,\alpha}; U},\ 
\|B\|_{{k,\alpha}; U},
\)

The iteration starts by setting $\bff_{(0)} = \bff$ and $\widehat{G}_{(0)} = (\hat{z}_{(0)}, \hat{\xi}_{(0)}) = (T'_0)^{-1} (\bff_{(0)})$. Inductively, given \( \widehat{G}_{(l-1)} = (\hat{z}_{(l-1)}, \hat{\xi}_{(l-1)}) \), we define the next input data as
\[
\bff_{(l)} := P_\delta \big(\widehat{G}_{(l-1)}\big),
\]
and solve
\[
\widehat{G}_{(l)} = (\hat{z}_{(l)}, \hat{\xi}_{(l)}) := \left(T_0'\right)^{-1}(\bff_{(l)}).
\]
Let  
\(
\widehat{G}_0 = \widehat{G}_{(0)}.
\)
The \( l \)-th approximation is then given by
\(
\widehat{G}_l := \widehat{G}_{l-1} + \widehat{G}_{(l)}.
\) 
It suffices to deal with the convergence of $(\widehat{G}_l)$ and the associated estimates.

If we define $(\hat{g}_{1,(l)}, \hat{g}_{2,(l)})$ by transforming the data $(\hat{z}_{(l)}, \hat{\xi}_{(l)})$ as in \eqref{defg1g2}: \(\hat{g}_{1,({l})} = \hat{\xi}_{(l)}- B\cdot\hat{z}_{(l)}- \bar{A} \cdot \hat{z}_{(l)} \) and \(\hat{g}_{2,(l)} = \hat{\xi}_{(l)}- B\cdot\hat{z}_{(l)}+ \bar{A} \cdot \hat{z}_{(l)} \), then, $(\hat{g}_{1,(l)}, \hat{g}_{2,(l)})$ solves the standard linear free boundary problem \eqref{equstandardfbpros} with the boundary data $\bff_{(l)}$. 
The solution $(\hat{g}_{1,(l)}, \hat{g}_{2,(l)})$ is obtained by applying the Hilbert transform to $\re \bff_{(l)}$ and $\im \bff_{(l)}$, respectively, and extending these boundary values harmonically into the disk $D$. 
The estimates for $(\hat{g}_{1,(l)}, \hat{g}_{2,(l)})$ can be obtained by Remark \ref{remregg1g2} and Corollary \ref{corHilbertBMOHolder}. The inverse transformation \eqref{inverserepzxibyg} then provides the estimates for $(\hat{z}_{(l)}, \hat{\xi}_{(l)})$. When $l=0$: in the $D$ direction,
\begin{align*}
    \| (\hat{z}_{(0)} (w, \cdot), \hat{\xi}_{(0)} (w, \cdot))\| _{k,\alpha; D} \leq C_0\| \bff(w, \cdot)\| _{k,\alpha; \partial D}, \quad \text{for } w\in U
\end{align*}
and in the $U$ direction, we estimate the modulus of continuity of $\nabla_0^{\mathbf{k}}(\hat{z}_{(0)}, \hat{\xi}_{(0)}) $. Recall the notation in \eqref{defckalphainN}, 
\[
h^{\mathbf{k},\alpha}_{(\hat{z}_{(0)}, \hat{\xi}_{(0)})} (\tau; w,w') = \frac{ \nabla_0^{\mathbf{k}}(\hat{z}_{(0)}, \hat{\xi}_{(0)}) (w',\tau)- \nabla_0^{\mathbf{k}}(\hat{z}_{(0)}, \hat{\xi}_{(0)}) (w,\tau) }{|w-w'|^\alpha}, \qquad \text{ for } \tau \in D, ~ w,w' \in U.
\]
Then, by Remark \ref{remregg1g2}, we have
\[
\|(\hat{g}_{1, (0)}, \hat{g}_{2,(0)})\|_{k,\beta; U \times D} \le C_0 (\alpha- \beta)^{-1} \|\bff\|_{k,\alpha; U \times \partial D} 
\]
By the inverse transformation \eqref{inverserepzxibyg}, we have 
\begin{align*}
\begin{split} 
     h^{\mathbf{k},\alpha}_{\hat{z}_{(0)}}  (\tau; w, w') = & ~ K_0(w,w', \tau) + \frac{1}{2} \bar{A}^{-1} \left( h^{\mathbf{k}, \alpha}_{\hat{g}_{2,(0)}} (\tau; w, w') - h^{\mathbf{k}, \alpha}_{ \hat{g}_{1,(0)} } (\tau; w, w') \right),\\
     h^{\mathbf{k},\alpha}_{\hat{\xi}_{(0)}}  (\tau; w, w') = & ~ K_1(w,w', \tau) + \frac{1}{2} \left( h^{\mathbf{k}, \alpha}_{\hat{g}_{1,(0)}} (\tau; w, w') + h^{\mathbf{k}, \alpha}_{ \hat{g}_{2,(0)} } (\tau; w, w') \right)\\
     & \hspace{2.2cm} + \frac{1}{2} B \bar{A}^{-1} \left( h^{\mathbf{k}, \alpha}_{\hat{g}_{2,(0)}} (\tau; w, w') - h^{\mathbf{k}, \alpha}_{ \hat{g}_{1,(0)} } (\tau; w, w') \right),
    \end{split}
\end{align*}
where $K_0(w, w', \tau)$ and $K_1 (w, w', \tau)$ are the collection of lower-order terms that consist of mixed derivatives of $\bar A^{-1}$ and $B$ up to $C^{k,\alpha}$, and derivatives of $(\hat{g}_{1, (0)}, \hat{g}_{2, (0)})$ up to order $k$. $K_{0,1}$ has the $L^\infty$ bound: 
\(
\|K_0\|_{L^\infty (U \times U \times \partial D)},~ \|K_1\|_{L^\infty (U \times U \times \partial D)} \le C_1 \alpha^{-2} \|\bff \|_{k,\alpha; U \times \partial D}.
\)
Since the Hilbert transform keeps $L^2$ norms on $\partial D$ after normalization, we have
\[
\| h^{\mathbf k, \alpha}_{(\hat{z}_{(0)},\hat{\xi}_{(0)})} (\cdot; w, w') \|_{L^2 (\partial D)} \le C_0 \|\bff\|_{k,\alpha; U \times \partial D} + \sqrt{2\pi} m_0(w, w'),
\]
where $m_0(w, w')$ is the absolute value of average of $h^{\mathbf k, \alpha}_{(\hat{z}_{(0)},\hat{\xi}_{(0)})} (\cdot; w, w')$ over $\partial D$. 
And by \eqref{bddCZOBMOLinfty} and Lemma \ref{lembddBMOitera}, BMO norm of $h^{\mathbf{k},\alpha}_{(\hat{z}_{(0)}, \hat{\xi}_{(0)})}$ satisfies
\begin{align}\label{induczxi0NBMO}
    \sup_{w,w' \in U}\|  h^{\mathbf{k},\alpha}_{(\hat{z}_{(0)}, \hat{\xi}_{(0)})} (\cdot;w, w')\| _{\BMO;~ \partial D} \leq C_0 \| \bff\| _{k,\alpha; U \times \partial D}.
\end{align}

It suffices to prove the following estimates inductively for $l \ge 1$ and $0<\delta\le \delta_0$:
\begin{align}
 \sup_{w\in U}\| \bff_{(l)}(w, \cdot)\| _{k,\alpha; \partial D} \leq &~ 2^{1-l} C_0 \delta \sup_{w\in U} \| \bff(w, \cdot)\| _{k,\alpha; \partial D} \tag{$l$1} \label{inducbnD}
\\
        \sup_{w\in U}\|(\hat{z}_{(l)} (w, \cdot), \hat{\xi}_{(l)}(w, \cdot)) \|_{k,\alpha; D} 
        &\leq 2^{1-l}  C_0^2 \delta \sup_{w\in U} \| \bff(w, \cdot) \|_{k,\alpha; \partial D}, \tag{$l$2} \label{induczxinD} 
\end{align}
and
\begin{align}
        \sup_{w, w' \in U} \| h^{\mathbf{k},\alpha}_{\bff_{(l)}}(\cdot; w,w')\| _{\BMO; \partial D} &\leq 2^{1-l} C_2 \delta \left(\| \bff\| _{k,\alpha; U \times \partial D}+m_0(w,w')\right), \tag{$l$3} \label{inducbnN}
        \\
        \sup_{w, w' \in U} \| h^{\mathbf{k},\alpha}_{\bff_{(l)}}(\cdot; w,w')\| _{L^2 (\partial D)} &\leq 2^{1-l}  C_2 \delta\left(\| \bff\| _{k,\alpha; U \times \partial D}+m_0(w,w')\right)
        \tag{$l$4} \label{inducbnintegral}
        \\
        \sup_{w, w'\in U}\|  h^{\mathbf{k},\alpha}_{(\hat{z}_{(l)}, \hat{\xi}_{(l)})} (\cdot;w, w')\| _{\BMO;~ \partial D} 
        &\leq 2^{1-l}  C_3 \delta\left(\| \bff\| _{k,\alpha; U \times \partial D}+m_0(w,w')\right), \tag{$l$5} \label{induczxinNBMO}
        \\
         \sup_{w, w'\in U}\|  h^{\mathbf{k},\alpha}_{(\hat{z}_{(l)}, \hat{\xi}_{(l)})} (\cdot;w, w')\| _{L^2(\partial D)} 
        &\leq 2^{1-l}  C_3 \delta\left(\| \bff\| _{k,\alpha; U \times \partial D}+m_0(w,w')\right). \tag{$l$6} \label{induczxinNintegral}
    \end{align}
By Lemma \ref{lemBMOHolderInterpolation}, by lowering the exponent from \(\alpha\) to \(\beta\), the following Hölder estimate holds:
\begin{align}
\begin{aligned}
\|(\hat z_{(l)},\hat\xi_{(l)})\|_{k,\beta;U\times D}
\le 2^{1-l}C_4 (\alpha -\beta)^{-1}\delta
\left(
\| \bff\| _{k,\alpha; U \times \partial D} + m_0(w,w')
\right).
\end{aligned} \tag{$l$7 \(\beta\)} \label{induczxilossHolder}
\end{align}
Here, we point out that all constants $C_2$, $C_3$ and $C_4$ are from absorbing the mixed terms appearing in the process of iterative estimates. The sufficient small constant $\delta_0$ is picking based on $C_j$, $j =0,1,\ldots, 4$ as follows:
\begin{align}\label{constantdelta0}
\delta_0 = \left(10 \sum_{j=0}^4 C_j \right)^{-1}
\end{align}

By the definition of $\bff_{(l+1)}$,
\[
\begin{aligned}
\bff_{(l+1)} &= P_{\delta} (\hat{z}_{(l)}, \hat{\xi}_{(l)}) \\
&= (\widetilde{A}(w, \tau)-A(w)) \bar{\hat{z}}_{(l)}(w,\tau)  + (\widetilde{B}(w, \tau)- B(w)) \hat{z}_{(l)} (w, \tau),
\end{aligned}
\]
together with the smallness condition \eqref{constantsmallperturb}, we have
\begin{align}
    \sup_{w\in U}\| \bff_{(l+1)}(w, \cdot)\| _{k,\alpha; \partial D} &\leq  \delta \sup_{w\in U} \| (\hat{z}_{(l)}, \hat{\xi}_{(l)}) (w, \cdot)\| _{k,\alpha; D}  \notag
    \\
    &\leq 2^{-l} C_0\delta \sup_{w\in U} \| \bff(w, \cdot)\| _{k,\alpha; \partial D}. \tag{($l$+1)1} \label{inducbn+1D}
\end{align}
The pair $(\hat{z}_{(l+1)}, \hat{\xi}_{(l+1)})$ solves the linear equation \eqref{holodisclinearpde-lemma} by setting the boundary data $\bff_{(l+1)}$. By Remark \ref{remregg1g2}, we have
\begin{align}
    \sup_{w\in U}\| (\hat{z}_{(l+1)}, \hat{\xi}_{(l+1)})(w, \cdot)\| _{k,\alpha; D} 
    &\leq C_0 \sup_{w\in U}\| \bff_{(l+1)}(w, \cdot)\| _{k,\alpha; \partial D} \notag
    \\
    &\leq 2^{-l}  C^2_0 \delta \sup_{w\in U} \| \bff(w, \cdot) \|_{k,\alpha; \partial D}. \tag{($l$+1)2} \label{induczxin+1D}
\end{align}
Hence, we complete the proof of $((l+1)1)$ and $((l+1)2)$. To prove $((l+1)3)$ and $((l+1)4)$, we analyze the difference quotient of the $k$-th derivatives for $\bff_{(l+1)}$:
\begin{align}\label{calc:hhalphabn+1}
\begin{split}
h^{\mathbf{k}, \alpha}_{{\bff}_{(l+1)}}(\tau; w, w') &= \frac{\nabla_0^{\mathbf{k}}\bff_{(l+1)} (w',\tau)- \nabla_0^{\mathbf{k}}\bff_{(l+1)} (w,\tau)}{|w-w'|^\alpha} \\
&=  K_2(w, w', \tau)+ (\widetilde{A}-A)  \cdot h^{\mathbf{k}, \alpha}_{\bar{\hat{z}}_{(l)}}(\tau; w, w') \\
&\hspace{2.55cm}
+ (\widetilde{B}-B) \cdot h^{\mathbf{k}, \alpha}_{{\hat{z}}_{(l)}}(\tau; w, w').
\end{split}
\end{align}
Here, $K_2(w, w', \tau)$ represents the collection of lower-order terms arising from the Leibniz rule; it consists of mixed derivatives of $(\widetilde{A}-A)$ and $(\widetilde{B}-B)$ up to $C^{k,\alpha}$, as well as derivatives of $(\bar{\hat{z}}_{(l)}, \hat{z}_{(l)})$ up to order $k$. By applying \eqref{constantsmallperturb}, \eqref{induczxilossHolder} for \(\beta =0\), we obtain the following $L^\infty$ bound for $K_2$:
\begin{align} \label{n+1-3Linftybd}
\| K_2(w, w', \tau)\| _{L^\infty ( U \times U \times \partial D)} \leq 2^{-l} C_1 \delta \| \bff\| _{k,\alpha; U \times \partial D}.
\end{align}
Then, by \eqref{calc:hhalphabn+1} and \eqref{n+1-3Linftybd}, 
\begin{align}
    \| h^{\mathbf{k}, \alpha}_{\bff_{(l+1)}} (\cdot; w, w')\| _{L^2 (\partial D)} &\le \sqrt{2\pi} \| K_2\| _{L^\infty(U \times U \times \partial D)} + \| \tilde{A} - A \| _{L^\infty(U \times \partial D)} \| h_{\hat{z}_{(l)}}^{\mathbf{k}, \alpha} (\cdot; w, w')\| _{L^2(\partial D)}  \notag
    \\ 
    & \hspace{4.2cm} + \| \tilde{B} - B \| _{L^\infty(U \times \partial D)} \| h_{\hat{z}_{(l)}}^{\mathbf{k}, \alpha} (\cdot; w, w')\| _{L^2(\partial D)} \notag \\
    & \le 2^{-l} (\sqrt{2\pi} C_1 \delta) \| \bff\| _{k,\alpha; U \times \partial D} + 2^{1-l}(2C_3 \delta^2) \left(\|\bff\|_{k,\alpha; U \times \partial D} + m_0 (w,w') \right) \notag
    \\
    & \le 2^{-l} C_2 \delta \left(\|\bff\|_{k,\alpha; U \times \partial D} + m_0 (w,w') \right). \tag{($l$+1)4} \label{inducbn+1integral}
\end{align}
Note that \(\|\cdot\|_{\BMO} \le 2 \|\cdot\|_{L^\infty}\), together with the product Lemma \ref{lemBMOproduc}, we have
\begin{align}
     \| h^{\mathbf{k}, \alpha}_{\bff_{(l+1)}}& (\cdot; w, w')\| _{\BMO; \partial D}  \le  2 \|K_2\|_{L^\infty (U \times U \times \partial D)} \notag
     \\ 
     & + C\left( \|\tilde{A} - A\|_{k,\alpha; U \times \partial D} +  \|\tilde{B} - B\|_{k,\alpha; U \times \partial D} \right)
     \left(\|h_{\hat{z}_{(l)}}^{\mathbf{k}, \alpha} (\cdot; w, w') \|_{\BMO; \partial D} + \|h_{\hat{z}_{(l)}}^{\mathbf{k}, \alpha} (\cdot; w, w') \|_{L^2 (\partial D)}\right) \notag
     \\
     & \hspace{2.9cm} \le 2^{-l} (2 C_1 \delta) \|\bff\|_{k,\alpha; U \times \partial D} + 2^{1-l} (4C_3\delta^2)\left(\|\bff\|_{k,\alpha; U \times \partial D} + m_0 (w,w') \right) \notag
     \\
     & \hspace{2.9cm} \le 2^{-l} C_2 \delta \left(\|\bff\|_{k,\alpha; U \times \partial D} + m_0 (w,w') \right) \tag{($l$+1)3} \label{inducbn+1BMO}
\end{align}

Now, we start to prove $((l+1)5)$ and $((l+1) 6)$. Note that the Hilbert transform in $\tau$ is bounded on $L^2$, hence the vector-valued functions, \(h^{\mathbf k, \alpha}_{\re\bff_{(l+1)}}\), \(h^{\mathbf k, \alpha}_{\im\bff_{(l+1)}}\), and their Hilbert transform \(h^{\mathbf k, \alpha}_{\Hilb \re\bff_{(l+1)}}\), \(h^{\mathbf k, \alpha}_{\Hilb\im\bff_{(l+1)}}\),
have \(L^2\)-bound by \eqref{inducbn+1integral}.
Hence, we have the estimate for the $L^2$ norm of $h^{\mathbf{k},\alpha}_{(\hat{g}_{1,(l+1)}, \hat{g}_{2,(l+1)})}$ on $\partial D$: 
\begin{align*}
\| h^{\mathbf k, \alpha}_{(\hat{g}_{1,(l+1)}, \hat{g}_{2, (l+1)})} (\cdot; w, w')\| _{L^2 (\partial D)} \le C 
&\bigl( \| h_{\bff_{(l+1)}}^{\mathbf{k}, \alpha}  (\cdot; w, w')\| _{L^2(\partial D)} 
\\ 
&+ \sqrt{2\pi}
\left|\fint_{\partial D}
h^{\mathbf k,\alpha}_{(\hat g_{1,(l+1)},\hat g_{2,(l+1)})}
(\theta;w,w')\,d\theta\right|\bigr).
\end{align*}
Then it suffices to estimate the average of \(h^{\mathbf k,\alpha}_{(\hat g_{1,(l+1)},\hat g_{2,(l+1)})}\) over \(\partial D\). One of the key difficulty is to show that the average is not accumulated in the process of iteration. 
Note the normalization provides $\hat{z}_{(l)} (w, -i) = 0$; hence, $h_{\hat z_{(l)}}^{\mathbf k,\alpha}(-i;w,w')=0,$
\[
\operatorname{Hil}
h_{\hat z_{(l)}}^{\mathbf k,\alpha}(-i;w,w')
=i\fint_{\partial D}
h_{\hat z_{(l)}}^{\mathbf k,\alpha}(\theta;w,w')\,d\theta.
\]
Note that $\left|\fint_{\partial D}
h^{\mathbf k,\alpha}_{(\hat g_{1,(l+1)},\hat g_{2,(l+1)})}
(\theta;w,w')\,d\theta\right| \le C \left(|\fint_{\partial D} h^{\mathbf k,\alpha}_{\bff_{(l+1)}}(\theta;w,w')\,d\theta| + |\operatorname{Hil}
h_{\bff_{(l+1)}}^{\mathbf k,\alpha}(-i;w,w')| \right)$.
Now, if we apply Hilbert transform to the prodects:
\[
\begin{aligned}
\operatorname{Hil}\!&\left[
(\widetilde B-B)(w',\cdot)\,
h_{\hat z_{(l)}}^{\mathbf k,\alpha}(\cdot;w,w')
\right](-i)
=
i(\widetilde B-B)(w',-i)
\fint_{\partial D}h_{\hat z_{(l)}}^{\mathbf k,\alpha}(\theta;w,w')\,d\theta\\
&\hspace{2.5cm}+
\fint_{\partial D}
\cot\!\left(-\frac{\pi}{4}-\frac{\theta}{2}\right)
\bigl[\widetilde B(w',\theta)-\widetilde B(w',-i)\bigr]\,
h_{\hat z_{(l)}}^{\mathbf k,\alpha}(\theta;w,w')\,d\theta.
\end{aligned}
\]
The second term satisfies the estimate:
\[
\begin{aligned}
\left|
\cot\!\left(-\frac{\pi}{4}-\frac{\theta}{2}\right)
\bigl[\widetilde B(w',\theta)-\widetilde B(w',-i)\bigr]
\right| 
\le C\|\partial_\theta(\widetilde B-B)(w',\cdot)\|_\infty
\le C\delta.
\end{aligned}
\]
Therefore, we have
\[
\begin{aligned}
\left|
\operatorname{Hil}\!\left[
(\widetilde B-B)(w',\cdot)\,
h_{\hat z_{(l)}}^{\mathbf k,\alpha}(\cdot;w,w')
\right](-i)
\right|
\le C\delta\fint_{\partial D}
\left|h_{\hat z_{(l)}}^{\mathbf k,\alpha}(\theta;w,w')\right|\,d\theta
\le C\delta
\|h_{\hat z_{(l)}}^{\mathbf k,\alpha}(\cdot;w,w')\|_{L^2},
\end{aligned}
\]
and
\[
\begin{aligned}
\left|\fint_{\partial D}
(\widetilde B-B)(w',\theta)\,
h_{\hat z_{(l)}}^{\mathbf k,\alpha}(\theta;w,w')\,d\theta\right|
\le\delta\fint_{\partial D}
\left|h_{\hat z_{(l)}}^{\mathbf k,\alpha}(\theta;w,w')\right|\,d\theta
\le C\delta
\|h_{\hat z_{(l)}}^{\mathbf k,\alpha}(\cdot;w,w')\|_{L^2}.
\end{aligned}
\]
Adding together, we have
\[
\begin{aligned}
&\left|\fint_{\partial D}
(\widetilde B-B)(w',\theta)\,
h_{\hat z_{(l)}}^{\mathbf k,\alpha}(\theta;w,w')\,d\theta\right|\\
&\quad+
\left|
\operatorname{Hil}\!\left[
(\widetilde B-B)(w',\cdot)\,
h_{\hat z_{(l)}}^{\mathbf k,\alpha}(\cdot;w,w')
\right](-i)
\right|
\le C\delta
\|h_{\hat z_{(l)}}^{\mathbf k,\alpha}(\cdot;w,w')\|_{L^2}.
\end{aligned}
\]
A similar estimate holds for the term with $(\widetilde A-A)h^{\mathbf k,\alpha}_{\hat{z}_{l}}$. Using \eqref{calc:hhalphabn+1}, \eqref{induczxinNintegral} and the estimate for $K_2$, we have
\begin{align*}
\left|\fint_{\partial D} h^{\mathbf k,\alpha}_{(\hat g_{1,(l+1)},\hat g_{2,(l+1)})} (\theta;w,w')\,d\theta\right| &\le C \left(|\fint_{\partial D} h^{\mathbf k,\alpha}_{\bff_{(l+1)}}(\theta;w,w')\,d\theta| + |\operatorname{Hil}
h_{\bff_{(l+1)}}^{\mathbf k,\alpha}(-i;w,w')| \right)
\\
&\le 2^{-l} C C_2 \delta \left(\|\bff\|_{k,\alpha; U \times \partial D} + m_0(w,w') \right),
\end{align*}
where the \(\|\bff_{(l+1)}\|_{0,\alpha; U \times \partial D}\) term comes from the normalization constants in \eqref{equstandardfbpros}. It is uniformly controlled as follows. 

By the inverse transformation \eqref{inverserepzxibyg}, we have
\begin{align}
\begin{split} \label{hkalphazxin+1}
     h^{\mathbf{k},\alpha}_{\hat{z}_{(l+1)}}  (\tau; w, w') = & ~ K_3(w,w', \tau) + \frac{1}{2} \bar{A}^{-1} \left( h^{\mathbf{k}, \alpha}_{\hat{g}_{2,(l+1)}} (\tau; w, w') - h^{\mathbf{k}, \alpha}_{ \hat{g}_{1,(l+1)} } (\tau; w, w') \right),\\
     h^{\mathbf{k},\alpha}_{\hat{\xi}_{(l+1)}}  (\tau; w, w') = & ~ K_4(w,w', \tau) + \frac{1}{2} \left( h^{\mathbf{k}, \alpha}_{\hat{g}_{1,(l+1)}} (\tau; w, w') + h^{\mathbf{k}, \alpha}_{ \hat{g}_{2,(l+1)} } (\tau; w, w') \right)\\
     & \hspace{2.2cm} + \frac{1}{2} B \bar{A}^{-1} \left( h^{\mathbf{k}, \alpha}_{\hat{g}_{2,(l+1)}} (\tau; w, w') - h^{\mathbf{k}, \alpha}_{ \hat{g}_{1,(l+1)} } (\tau; w, w') \right),
    \end{split}
\end{align}
where $K_3(w, w', \tau)$ and $K_4 (w, w', \tau)$ represent the collection of lower-order terms that consist of mixed derivatives of $\bar A^{-1}$ and $B$ up to $C^{k,\alpha}$,and derivatives of $(\hat{g}_{1, (l+1)}, \hat{g}_{2, (l+1)})$ up to order $k$. By Lemma \ref{lemBMOHolderInterpolation} and Corollary \ref{corHilbertBMOHolder}, together with \eqref{inducbn+1D}, \eqref{inducbn+1BMO} and \eqref{inducbn+1integral}, we have 
\begin{align} \label{estig12kbeta}
\begin{split}
    \|(\hat{g}_{1, (l+1)}, \hat{g}_{2, (l+1)})\|_{k,\beta; U \times \partial D} \le &~  C (\alpha-\beta)^{-1} \bigg( \sup_{w \in U} \| \bff_{(l+1)}(w,\cdot) \|_{k,\alpha;\partial D}
     \\ 
     &+ \sup_{w,w' \in U} \left( \| h^{\mathbf{k}, \alpha}_{\bff_{(l+1)}} (\cdot; w,w')\| _{\BMO;~ \partial D} + \| h^{\mathbf{k}, \alpha}_{\bff_{(l+1)}} (\cdot; w,w')\| _{L^2(\partial D)}  \right)  \bigg)\\
     \le & ~ 2^{-l} C (\alpha-\beta)^{-1} (C_0 + 2 C_2) \delta (\|\bff\|_{k,\alpha; U \times \partial D} + m_0(w,w')).
\end{split}
\end{align}
Hence, we have the $L^\infty$ bound for $K_3$ and $K_4$:
\begin{align*}
    \|K_3 \|_{L^\infty(U \times U \times \partial D)}, ~ \|K_4 \|_{L^\infty(U \times U \times \partial D)} \le 2^{-l}\alpha^{-1}C_1 C_0(C_0 + 2C_2) \delta\|\bff\|_{k,\alpha; U \times  \partial D}
\end{align*}
Together with \eqref{hkalphazxin+1}, we have 
\begin{align}
\| h^{\mathbf{k},\alpha}_{(\hat{z}_{(l+1)}, \hat{\xi}_{(l+1)})}  (\cdot; w, w') \| _{L^2 (\partial D)} &\le \sum_{j=3,4} \sqrt{2\pi} \| K_j  \|_{L^\infty (U \times U \times \partial D)} + 2C_0 \| h^{\mathbf k , \alpha}_{(\hat{g}_{1, (l+1)}, \hat{g}_{2, (l+1)})} (\cdot; w,w') \|_{L^2 (\partial D)} \notag
\\
& \le 2^{-l} C C_0 (C_0 + 5 C_2 + 2C_4) (1 + \alpha^{-1}) \delta (  \|\bff\|_{k,\alpha; U \times \partial D} + m_0(w,w') ) \notag
\\
& \le 2^{-l} C_3 \delta(    \|\bff\|_{k,\alpha; U \times \partial D} + m_0(w,w') )
\tag{($l$+1)6} \label{induchzxin+1L2}
\end{align}
To prove $((l+1) 5)$, note that the normalization doesn't contribute to the BMO norms: 
\[
\| h^{\mathbf k, \alpha}_{(\hat{g}_{1,(l+1)}, \hat{g}_{2, (l+1)})} (\cdot; w, w')\|_{\BMO, \partial D} \le C  \| h_{\bff_{(l+1)}}^{\mathbf{k}, \alpha}  (\cdot; w, w')\| _{\BMO ;\partial D} .
\]
Then, by \eqref{hkalphazxin+1} and Lemma \ref{lemBMOproduc}, we have
\begin{align}
    \| h^{\mathbf{k},\alpha}_{(\hat{z}_{(l+1)}, \hat{\xi}_{(l+1)})}  (\cdot; w, w') \| _{\BMO;\partial D} \le & \sum_{j=3,4} 2\| K_j  \|_{L^\infty (U \times U \times \partial D)} \notag
    \\
    & + 2C_0 \bigg( \| h^{\mathbf k , \alpha}_{(\hat{g}_{1, (l+1)}, \hat{g}_{2, (l+1)})} (\cdot; w,w') \|_{\BMO; \partial D} \notag
    \\
    & \hspace{1.5cm}+ \| h^{\mathbf k , \alpha}_{(\hat{g}_{1, (l+1)}, \hat{g}_{2, (l+1)})} (\cdot; w,w') \|_{L^2 (\partial D)} \bigg) \notag
    \\  
    \le & ~ 2^{-l} C C_0 (C_0 + 2 C_2) (2 + \alpha^{-1}) \delta (\|\bff\|_{k,\alpha; U \times \partial D} + m_0(w,w') )  \notag
    \\
    \le &~ 2^{-l} C_3 \delta (\|\bff\|_{k,\alpha; U \times \partial D} + m_0(w,w') )  
    \tag{($l$+1)5}  \label{induchzxin+1BMO}
\end{align}
The $((l+1)7~ \beta)$ follows directly from Lemma \ref{lemBMOHolderInterpolation}, \eqref{inducbn+1D}, \eqref{inducbn+1BMO} and \eqref{inducbn+1integral}. We have:
\begin{align}
    \|(\hat{z}_{(l+1)}, \hat{\xi}_{(l+1)})\|_{k,\beta; U \times \partial D} & \le 2^{-l} C C_0 (\alpha-\beta)^{-1} (C_0 + 2C_2) \delta (\|\bff\|_{k,\alpha; U \times \partial D} + m_0(w,w') )\notag
    \\
    & \le 2^{-l} C_4 (\alpha-\beta)^{-1}\delta (\|\bff\|_{k,\alpha; U \times \partial D} + m_0(w,w') ) \tag{($l$+1)7 $\beta$}
\end{align}

It remains to prove that the approximate solutions $\widehat{G}_l$, denoted by $\widehat{G}_l = (\hat{z}_l,\hat{\xi}_l)$, converge to the solution to \eqref{holodisclinearpde'}. Note that
\begin{align*}
\bff- T'(\hat{z}_l,\hat{\xi}_l)  = \bff- \sum_{s=0}^l T'(\hat{z}_{(s)}, \hat{\xi}_{(s)}) = \bff_{(l+1)}.
\end{align*}
By the above discussion, $\bff_{(l+1)} = P_\delta (\widehat{G}_{(l)})$, and by \eqref{induczxilossHolder}, we have
\begin{align*}
    \| \bff- T'(\hat{z}_{l},\hat{\xi}_{l})\| _{k,\beta; U \times \partial D} &= \| \bff_{(l+1)}\| _{k,\beta; U \times \partial D} 
    \\
    &\leq (\alpha- \beta)^{-1} 2^{-l} C_4 \delta (\| \bff\| _{k,\alpha; U \times \partial D } + m_0(w,w') )\to 0.
\end{align*}
The solution to \eqref{holodisclinearpde'} is given by $ \displaystyle \widehat{G} = \lim_{l\rightarrow \infty} \widehat{G}_l$, satisfying the estimate:
\begin{align*}
    \| (\hat{z}, \hat{\xi})\| _{k,\beta; U \times D} \leq \sum_{l=0}^\infty \| (\hat{z}_{(l)}, \hat{\xi}_{(l)})\| _{k,\beta; U \times  D} \leq 2(\alpha-\beta)^{-1}C(\| \bff\| _{k,\alpha; U \times \partial D} + \sup_{w,w'} m_0(w,w') ).
\end{align*}
Therefore, $(\widehat{G}_l) $ converges to the solution to \eqref{holodisclinearpde'}, $\widehat{G} = (\hat{z}, \hat{\xi})$, in $\C^{k,\beta}$ norm.
Now, if we pick $\alpha' = (\alpha + \beta)/2$, we have the estimate:
\begin{align*}
    \| (\hat{z}, \hat{\xi})\| _{k,\beta; U \times D} \leq 2(\alpha-\beta)^{-1}C(\| \bff\| _{k,\alpha; U \times \partial D} + \sup_{w,w'} m_0(w,w') ) \le C (\alpha-\beta)^{-2} C \|\bff\|_{k,\alpha; U \times \partial D}.
\end{align*}

It remains to prove uniqueness. Suppose that there exists a bounded nonzero solution to 
\(T'(\hat z,\hat\xi)=0\),
and \(\hat z(w,-i)=0\). 
We consider the following mixed norm of \((\hat{z}, \hat{\xi}) \):
\begin{align*}
\| (\hat{z}, \hat{\xi}) \|_{\mathrm{w},k,\alpha;U \times \partial D} &=  \sup_{w \in U} \| (\hat{z}, \hat{\xi})(w,\cdot) \|_{k,\alpha;\partial D}
     \\ 
     &+ \sup_{w,w' \in U} \left( \| h^{\mathbf{k}, \alpha}_{(\hat{z}, \hat{\xi})} (\cdot; w,w')\| _{\BMO;~ \partial D} + \| h^{\mathbf{k}, \alpha}_{(\hat{z}, \hat{\xi})} (\cdot; w,w')\| _{L^2(\partial D)}  \right) 
\end{align*}
Since we assume $(\hat{z}, \hat{\xi})$ to be nonzero, after scaling, we normalize so that \(\| (\hat{z}, \hat{\xi}) \|_{\mathrm{w},k,\alpha;U \times \partial D}  =1\).

Let \(\mathbf \bff' = T'_0(\hat z,\hat\xi)\). Note that 
\(\mathbf \bff' = T'_0(\hat z,\hat\xi) =
(T'_0-T')(\hat z,\hat\xi)
=
P_\delta(\hat z,\hat\xi).\)
By \eqref{calc:hhalphabn+1}, the smallness conditions \eqref{constantsmallperturb},  the product Lemma \ref{lemBMOproduc} and the interpolation Lemma \ref{lemBMOHolderInterpolation}, we have
\begin{align*}
\| \bff'\|_{\mathrm{w},k,\alpha;U \times \partial D} &= \sup_{w \in U} \| \bff'(w,\cdot) \|_{k,\alpha;\partial D}
     \\ 
     &+ \sup_{w,w' \in U} \left( \| h^{\mathbf{k}, \alpha}_{\bff'} (\cdot; w,w')\| _{\BMO;~ \partial D} + \| h^{\mathbf{k}, \alpha}_{\bff'} (\cdot; w,w')\| _{L^2(\partial D)}  \right)  
\\
&\le C \delta (C_1 \alpha^{-1} + 2) \bigg( \sup_{w \in U} \|\hat{z}(w,\cdot)\|_{k,\alpha; \partial D} + \sup_{w,w' \in U} \big(\|h^{\mathbf k,\alpha}_{\hat{z}} (\cdot; w, w') \|_{L^2(\partial D)} 
\\ 
&\hspace{6.1cm}+\|h^{\mathbf k, \alpha}_{\hat{z}} (\cdot; w, w')\|_{\BMO; \partial D}
\big)\bigg) 
\\
&\le C_5 \delta  \|(\hat{z}, \hat{\xi})\|_{\mathrm w, k,\alpha; U \times \partial D} \le C_5\delta
\end{align*}
To solve $T_0' (\hat{z}, \hat{\xi}) =  \bff'$, we follow the proof of \eqref{induczxin+1D}, \eqref{induchzxin+1BMO} and \eqref{induchzxin+1L2}. There 
\begin{align*}
    \|(\hat{z}, \hat{\xi})\|_{\mathrm{w}, k,\alpha; U \times \partial D} \le C C_0 (C_1 \alpha^{-1} +2 ) \|\bff'\|_{\mathrm w;k,\alpha; U \times \partial D} \le C_6 \delta.
\end{align*}
Choosing \(\delta_0\) sufficiently small so that \(C_6\delta_0<1/2\), we obtain
\[
1=\|(\hat z,\hat\xi)\|_{\mathrm w;k,\alpha;U\times\partial D}< \frac12,
\]
a contradiction.

Finally, we address the global well-definedness of the solution. Let $(\hat{z}, \hat{\xi})_j$ and $(\hat{z}, \hat{\xi})_k$ denote the solutions to \eqref{holodisclinearpde'} obtained on overlapping holomorphic charts $U_j$ and $U_k$, respectively. Recall that on the intersection $U_j \cap U_k$, the coefficient matrices satisfy the relation $(\widetilde{B})_j = (\widetilde{B})_k - W_{jk}$. The uniqueness of the solution to the linearized problem implies that the local solutions must agree on the overlap according to the transition law of the bundle $\mathcal{W}'$:
$$(\hat{z}, \hat{\xi})_k = (\hat{z}, \hat{\xi} + (W)_{jk} \hat{z})_j.$$
Consequently, the collection of local solutions patches together to define a unique global section $\widehat{G} \in \Gamma (U \times D, Z^* \W')$.
\end{proof}

\section{Existence and Local Uniqueness of HCMA Disc Family}\label{secholodiscfoliexist}

In this section, we solve the nonlinear problem \eqref{holodiscpde} for boundary data sufficiently close to the trivial one. We work near the trivial pair \((G_0,\Lambda_0)\), where, in each local
holomorphic chart, 
\(
G_0(w,\tau)=(w,\partial\rho(w)),
\)
and \(\Lambda_0\subset \W\) is the graph of \(\partial\rho\).
For a perturbed boundary potential \(\Psi\), we denote by \(\Lambda_\Psi\) the family of exact Lagrangian submanifolds whose fibre over \(\tau\in\partial D\) is locally given by the graph of 
\(
\partial (\rho + \Psi (\cdot, \tau)).
\)
The difference of the boundary data is measured by $\partial \Psi$. We assume that 
\begin{align}\label{con:smallperturb_bdry}
\|\partial\Psi\|_{k+2,\alpha;N'\times\partial D}\leq \varepsilon_0,
\qquad k\geq 1,\quad \alpha\in(0,1),
\end{align}
where \(\varepsilon_0>0\) is a sufficiently small uniform constant to be chosen below. We first prove a local uniqueness statement for disc families near \(G_0\), and then establish existence by an iteration scheme based on the linear estimates obtained in the previous section.

\subsection{Local uniqueness of foliation by holomorphic discs with a fixed potential} \label{secuniq2ndver} 

In this subsection, we prove local uniqueness for HCMA disc families near
the trivial family \(G_0\). More precisely, we show that, for a fixed
boundary potential \(\Psi\), there is at most one solution of
\eqref{holodiscpde} in a sufficiently small \(C^{k,\beta}\)-neighborhood
of \(G_0\).

Recall that \(\mathcal F\) denotes the space of smooth families of
holomorphic discs defined in \eqref{defspaceholodiscs}. We now
introduce the local neighborhood of \(G_0\) used in the uniqueness
argument. For \(\delta_1>0\), define
\begin{align}\label{familyholodiscsngbhcondi}
\mathcal F^{k,\beta}_{\delta_1}(G_0)
:=
\left\{
G\in\mathcal F \;:\;
\operatorname{dist}_{k,\beta;N\times D}(G,G_0)\leq \delta_1
\right\}.
\end{align}
Here \(\operatorname{dist}_{k,\beta;N\times D}\) is the \(C^{k,\beta}\)
distance measured in the fixed local trivializations used to define the
norms on \(\W\). Thus, if on a holomorphic chart we write
\(
G(w,\tau)=(z(w,\tau),\xi(w,\tau)),
\qquad
G_0(w,\tau)=(z_0(w,\tau),\xi_0(w,\tau)),
\)
then the condition \(G\in\mathcal F^{k,\beta}_{\delta_1}(G_0)\) means, in
particular,
\[
\|z-z_0\|_{k,\beta;N\times D}
+
\|\xi-\xi_0\|_{k,\beta;N\times D}
\leq \delta_1
\]
on each such chart, with the global distance obtained by taking the
supremum over the chosen atlas.

\begin{thm}[Local uniqueness near the trivial disc family] \label{thmlocalunifixpotential}
Let \((G_0,\Lambda_0)\) be the trivial pair. Fix \(k\geq 1\), \(0<\beta<\alpha<1\).
There exist small constants \(\varepsilon_0,\delta_1>0\), depending only on
\(n,k,\alpha,\beta,\sigma^{-1}\), and
\(
\|\nabla^2_X\rho\|_{k+1,\alpha;U^+},
\)
such that the following holds.

Let \(\Psi\in C^\infty(U^+\times\partial D,\mathbb R)\) be a smooth
boundary potential satisfying \eqref{con:smallperturb_bdry}.
Then there is at most one smooth HCMA disc family
\(
G\in \mathcal F^{k,\beta}_{\delta_1}(G_0)
\)
adapted to the boundary data \(\Lambda_\Psi\). 
\end{thm}

The key input in the proof of Theorem \ref{thmlocalunifixpotential} is
the uniqueness of the linearized problem at data close to the trivial pair
\((G_0,\Lambda_0)\). Let
\(G\in \mathcal F^{k,\beta}_{\delta_1}(G_0)\), and let
\(\Lambda_\Psi\) be determined by a boundary potential \(\Psi\) satisfying
\eqref{con:smallperturb_bdry}. Locally, write
\(
G(w,\tau)=(z(w,\tau),\xi(w,\tau)).
\)
Then the linearized equation at \((G,\Lambda_\Psi)\) is locally given by \eqref{holodisclinearpde'} with $\widetilde{A}_{i\bar j} (w, \tau) = \partial_i \overline{\partial}_j (\rho + \Psi ) (z(w, \tau), \tau)$ and $\widetilde{B}_{ij} (w, \tau) = \partial_i \partial_j (\rho + \Psi) (z(w, \tau), \tau)$.

Since \(G\in\mathcal F^{k,\beta}_{\delta_1}(G_0)\) and
\(\Psi\) satisfies \eqref{con:smallperturb_bdry}, the coefficients, \( \widetilde{A}_{i \bar j}(w, \tau)\) and \(\widetilde{B}_{ij} (w, \tau)\) are small perturbations of the coefficients, $A_{i\bar j} (w)$ and $B_{ij}(w)$, 
at the trivial data. More precisely, there exists a uniform constant \(C_0\),
depending only on the background local potential through
\(
\|\nabla_X^2\rho\|_{k+1,\alpha;U^+},
\)
such that
\begin{align*}
&\left\|
\big(\partial_i\overline{\partial}_j(\rho+\Psi)\big)(z(w,\tau),\tau)
-
\partial_i\overline{\partial}_j\rho(w)
\right\|_{k,\beta;U\times\partial D}
\leq C_0(\delta_1+\varepsilon_0),\\
&\left\|
\big(\partial_i\partial_j(\rho+\Psi)\big)(z(w,\tau),\tau)
-
\partial_i\partial_j\rho(w)
\right\|_{k,\beta;U\times\partial D}
\leq C_0(\delta_1+\varepsilon_0).
\end{align*}
We choose \(\delta_1,\varepsilon_0>0\) sufficiently small so that
\begin{align}\label{para:smallconstsavgequ}
C_0(\delta_1+\varepsilon_0)\leq \delta,
\end{align}
where \(\delta>0\) is the small constant fixed with $0<\delta< \delta_0$ as in \eqref{constantdelta0}. With this choice, Lemma
\ref{lemlinearproblemperturb} can be applied here.

\begin{proof}[Proof of Theorem \ref{thmlocalunifixpotential}]
Fix a boundary potential \(\Psi\) satisfying
\eqref{con:smallperturb_bdry}. Suppose that
\[
G,\widetilde G\in \mathcal F^{k,\beta}_{\delta_1}(G_0)
\]
are two HCMA disc families adapted to the same boundary data
\(\Lambda_\Psi\). We will show that \(G=\widetilde G\).

Work in a local holomorphic trivialization of \(\W\) near \(G_0\), and write
\[
G(w,\tau)=(z(w,\tau),\xi(w,\tau)),
\qquad
\widetilde G(w,\tau)=(\widetilde z(w,\tau),\widetilde\xi(w,\tau)).
\]
Set
\[
\widehat G:=\widetilde G-G=(\hat z,\hat\xi).
\]
Since both \(G\) and \(\widetilde G\) belong to \(\mathcal F\), we have
\(
\overline{\partial}_\tau \hat z_i
=
\overline{\partial}_\tau \hat\xi_i
=0\) in \(U\times D,\) and \(\hat z_i(w,-i)=0.\)

Moreover, since both disc families satisfy the same boundary condition,
\[
T(G,\Lambda_\Psi)=0=T(\widetilde G,\Lambda_\Psi).
\]
Because \(G\) and \(\widetilde G\) are contained in a sufficiently small
coordinate neighborhood of \(G_0\), the affine path
\[
G_t:=G+t(\widetilde G-G),\qquad t\in[0,1],
\]
is well-defined in this local trivialization. By the fundamental theorem of
calculus,
\[
0
=
T(\widetilde G,\Lambda_\Psi)-T(G,\Lambda_\Psi)
=
\int_0^1
D_{\mathcal F}T\big|_{(G_t,\Lambda_\Psi)}(\widehat G)\,dt.
\]
Writing \(G_t=(z_t,\xi_t)\), we have on \(U\times\partial D\)
\begin{align}\label{equ:avglinearequlocal}
\hat\xi_i
-
(A_{\mathrm{avg}})_{i\bar j}\,\bar {\hat z}_j
-
(B_{\mathrm{avg}})_{ij}\,\hat z_j
=0,
\end{align}
where
\begin{align*}
(A_{\mathrm{avg}})_{i\bar j}(w,\tau)
&=
\int_0^1
\big(\partial_i\overline{\partial}_j(\rho+\Psi)\big)
(z_t(w,\tau),\tau)\,dt,\\
(B_{\mathrm{avg}})_{ij}(w,\tau)
&=
\int_0^1
\big(\partial_i\partial_j(\rho+\Psi)\big)
(z_t(w,\tau),\tau)\,dt.
\end{align*}

Thus \(\widehat G=(\hat z,\hat\xi)\) solves the homogeneous linearized
problem with coefficients \(A_{\mathrm{avg}}\) and \(B_{\mathrm{avg}}\).
We compare these coefficients with the trivial coefficients
\(
A_{i\bar j}=\partial_i\overline{\partial}_j\rho\) and
\(B_{ij}=\partial_i\partial_j\rho.
\)
Since \(G,\widetilde G\in \mathcal F^{k,\beta}_{\delta_1}(G_0)\), the
entire path \(G_t\) remains in the same \(\delta_1\)-neighborhood of
\(G_0\). Together with the smallness of \(\Psi\), this implies
\[
\|(A_{\mathrm{avg}})_{i\bar j}-A_{i\bar j}\|_{k,\beta;U\times\partial D}
\leq \delta,
\qquad
\|(B_{\mathrm{avg}})_{ij}-B_{ij}\|_{k,\beta;U\times\partial D}
\leq \delta,
\]
after the choice \eqref{para:smallconstsavgequ}. Therefore
\eqref{equ:avglinearequlocal} is a sufficiently small perturbation of the
trivial homogeneous linearized problem. By Lemma
\ref{lemlinearproblemperturb}, the only solution is
\[
(\hat z,\hat\xi)=0.
\]
Hence \(G=\widetilde G\) on \(U\times D\), proving the local uniqueness.
\end{proof}

\subsection{Existence via the Nash-Moser-Zehnder Theorem}

The method we employ here is Zehnder's version of the Nash-Moser implicit function theorem, building upon the classic works \cite{Nash19956Imbedding, Moser1961Tech, Zehnder1975GeneralizedIFT}. Specifically, given the particular pattern of loss of regularity established in Lemma \ref{IFTlinearpart} and Lemma \ref{lemlinearproblemperturb}, Zehnder’s framework provides the best fit for our setting. 

We now state the main theorem of this subsection:

\begin{thm} \label{thmexistbypertb} Let $(G_0, \Lambda_0)$ be the trivial {data, and fix $k\geq1$ and $0<\beta<\alpha<1$. Assume $N\subset N'$ and $d_\omega(N,X\setminus N')\geq r_0>0$, and use the coordinate balls of Section~\ref{sec:notations} with $R-r\geq R_0>0$.} Then, there exists a uniform small constant $\varepsilon >0$ depending only on $n$, $k$, $\alpha$, {$\beta$,} $\sigma^{-1},\| \nabla_X^2 \rho\| _{k+1,\alpha; N' \times \partial D}$, {$r_0^{-1},R_0^{-1}$, and a uniform upper bound for $r^{-1}$,} such that for each ${\Psi} \in \C^{\infty} (N' \times \partial D, \RR )$ satisfying {\eqref{con:smallperturb_bdry}} {with} {$\varepsilon_0\leq\varepsilon$},
there is a HCMA Discs Family $G: N \times D \rightarrow \W$ satisfying all conditions in Definition \ref{def:disc_family}.

Moreover, {we} have the tame estimate:
\begin{align} \label{estimeasureperturb}
    \| G-G_0\| _{k,\beta; N \times D} \leq C (\alpha-\beta)^{-2} \| \partial {\Psi}\| _{k, \alpha; N'\times \partial D}, 
\end{align}
where $C$ is a constant depending on $n$, $k$, $\alpha$, $\sigma^{-1}${, $\| \nabla^2_X \rho\| _{k+1, \alpha; N' \times \partial D}$}, {and a uniform upper bound for $r^{-1}$.}
\end{thm}
In {fixed local trivializations, write} the {nonlinear} {operator} {as $T:\F\times\B\rightarrow\D$}. 
To apply Zehnder's version of Nash-Moser implicit function theorem, we consider the following families of spaces with H\"older regularity $k \geq 1$ and $0< \alpha< 1$: 
\begin{align*}
    &\F^{k,\alpha} =\{ G \in \C^{k,\alpha}(N \times D, \W); \ \overline{\partial}_\tau G (w, \tau) =0,\  \pi\circ G (w, -i) =w \} \\
    \B^{k,\alpha} &=   \{ \Lambda_{{\Psi}} \subseteq  \W \times \partial D; \ {\psi}_\tau (z)= {\Psi}(z, \tau)  \in \C^{k+1,\alpha} (N' \times \partial D), \ \omega +i\partial \overline{\partial} {\psi}_\tau>0, ~ \tau \in \partial D \},
\end{align*}
{with the same projection condition on $\F^{k,\alpha}$ as in \eqref{defspaceholodiscs},} and $\D^{k,\alpha} = \C^{k,\alpha} (N \times \partial D, \CC^n)$ {in the chosen trivialization.} The operator $T$ can be automatically extended to $T : \F^{k,\alpha} \times \B^{k,\alpha} \rightarrow \D^{k,\alpha}$. By Lemma \ref{IFTlinearpart} and Lemma \ref{lemlinearproblemperturb}, for each pair of data $(G, \Lambda_{\Psi})$ satisfying the small perturbation conditions \eqref{constantsmallperturb}, the linearized operator $D_\F T$ at $(G, \Lambda_\Psi)$ admits an inverse $\eta|_{ (G, \Lambda_{\Psi})}$ such that
\begin{align}\label{inverseesti}
    \| \eta\big|_{(G, \Lambda_\Psi)} (\bff)\| _{k,\beta; N \times D} \leq C (\alpha -\beta)^{-2} \| \bff\| _{k,\alpha; N \times \partial D}, \qquad \text{for } \beta< \alpha.
\end{align}
To apply Zehnder's version of the Nash-Moser iteration, we must account for the fact that the tangent space changes as $G$ varies. Therefore, the iterative process is constructed locally. We perform the Nash-Moser iteration within each coordinate chart and verify that the resulting local solutions can be patched together.

Accordingly, we establish the required estimates for the quadratic error term $Q$ in a local setting. Let $G_1$ and $G_2$ be maps taking values in a single holomorphic coordinate chart of $\mathcal{W}$, and let $\widehat{G}_{1,2} = G_2 - G_1$ denote their vector difference in these coordinates. We define the quadratic error term as:$$Q(G_2, G_1; \Lambda_\Psi) = T(G_2, \Lambda_\Psi) - T(G_1, \Lambda_\Psi) - D_\mathcal{F} T|_{(G_1, \Lambda_\Psi)} (\widehat{G}_{1,2}).$$

\begin{lem} \label{lem:quadratic_error} Let $G_1, G_2 \in \mathcal{F}^{k,\alpha}$ {lie} {in the prescribed small $C^{k,\alpha}$ neighborhood of $G_0$} in a local holomorphic chart of {$\mathcal W$}, and let $\widehat{G}_{1,2} = G_2 - G_1$. Suppose that {$\| \widehat{G}_{1,2} \| _{k,\alpha; N \times D} \leq 1$}. Then, there is a uniform constant $C$ such that
    \begin{align} \label{estiTalorformula}
        \| Q(G_2, G_1;  \Lambda_\Psi)\| _{k,\alpha; N \times \partial D}
        \leq C \| \widehat{G}_{1,2}\| _{k,\alpha; N \times D}^{2}
    \end{align}
\end{lem}
\begin{proof}
    In this proof, we work in local coordinates where the maps are represented as $G_1 = (z_1, \xi^1)$ and $G_2 = (z_2, \xi^2)$. We denote the difference by $\widehat{G} = (\hat{z}, \hat{\xi}) = G_2 - G_1$. Additionally, to simplify the notation, we define the total potential locally as $\widetilde{\Psi} = \rho + \Psi$.
    
    Substituting the expressions for $T(G, \Lambda_\Psi)$ and its linearization $D_\F T|_{(G, \Lambda_\Psi)} (\widehat{G})$ into the definition of $Q (G_2, G_1; \Lambda_\Psi)$, we obtain:
    \begin{align*}
    \begin{split}
        Q(G_2, G_1;  \Lambda_\Psi)_i &= - \partial_i \widetilde{\Psi} (z_2, \tau) + \partial_i \widetilde{\Psi}(z_1, \tau)  + \partial_i \overline{\partial}_{j} \widetilde{\Psi}(z_1, \tau) \bar{\hat{z}}_{j} + \partial_i \partial_j \widetilde{\Psi} (z_1, \tau) \hat{z}_j\\
        & = - \int_0^1\partial_i \overline{\partial}_j \widetilde{\Psi} (z_t, \tau) \bar{\hat{z}}_j dt +  \partial_i \overline{\partial}_j \widetilde{\Psi} (z_1, \tau) \bar{\hat{z}}_j \\
        & \quad - \int_0^1 \partial_i \partial_j \widetilde{\Psi} (z_t, \tau) \hat{z}_j dt + \partial_i \partial_j \widetilde{\Psi} (z_1, \tau) \hat{z}_j,
    \end{split}
    \end{align*}
    where $z_t = t z_1 + (1-t) z_2$ for $t \in [0,1]$. By mean value theorem again, together with \cite[Lemma~2.2]{yao2025asymptotic}, we have 
    \begin{align*} 
        \| Q(G_2, G_1; \Lambda_\Psi)\| _{k,\alpha; N \times \partial D}
        &\leq C \| (\hat{z}, \hat{\xi})\| ^2_{k,\alpha; N\times D}
    \end{align*}
    where $C$ depends only on $n$, $k$, $\alpha$, $\| \nabla^2 \rho\| _{k+1,\alpha; N' \times \partial D}$ and $\varepsilon$.
\end{proof}

\begin{proof}[Proof of Theorem \ref{thmexistbypertb} (local version)]
Zehnder's version of the Nash-Moser iteration can be carried out using estimates \eqref{inverseesti} and \eqref{estiTalorformula}. Starting from the trivial data given in $(G_0, \Lambda_0)$ with  $G_0 = (z_0, \xi^0) = (w, \partial \rho(w))$, the sequence $(G_n) = (z_n, \xi^n)$ can be defined inductively by
\begin{align*}
    G_n = G_{n-1} + \widehat{G}_{n}, \qquad \widehat{G}_{n} = -\eta|_{(G_{n-1},\Lambda_\Psi)} (T(G_{n-1},\Lambda_\Psi) ).
\end{align*}
Using the fact that $\eta|_{(G, \Lambda_\Psi)}$ serves as the inverse of $D_{\mathcal{F}} T |_{(G, \Lambda_\Psi)}$, we obtain the following key formula:
\begin{align*}
    T(G_n,\Lambda_\Psi) = T(G_n,\Lambda_\Psi)- T(G_{n-1},\Lambda_\Psi) - D_\F T \big|_{(G_{n-1},\Lambda_\Psi)} (\widehat{G}_{n}) = Q(G_{n}, G_{n-1};\Lambda_\Psi)
\end{align*}
with the estimate from \eqref{inverseesti} to \eqref{estiTalorformula}, 
$$\| T(G_n, \Lambda_\Psi)\| _{k,\beta; N \times \partial D} \leq C \| \widehat{G}_{n}\| _{k,\beta; N \times D}^2 \leq C^3 (\alpha- \beta)^{-4} \| T(G_{n-1}, \Lambda_\Psi)\| _{k, \alpha; N \times \partial D}^2$$
for arbitrary $0<\beta < \alpha$. Enlarge $C\ge2$ to absorb the composition constants and the uniform H\"older embedding constants used below. Now fix $0<\beta<\alpha<1$, set $\beta_0=\alpha$, and let $\beta_n = (\alpha-\beta) \exp(-\lambda \kappa^n) + \beta $ for $n \geq 1$. Taking $\kappa = \frac{3}{2} $,  $ \lambda = - \log \frac{(\alpha -\beta)^4}{4C^3}$ and
$\varepsilon = \min \{\delta(\alpha-\beta)^{2}/(4C^3), \exp({-8} \lambda)\},$ where $\delta$ is the small constant such that Lemma \ref{lemlinearproblemperturb} holds,
the following estimates follow by induction, using
\[
 \|\widehat G_{n+1}\|_{k,\beta_{n+1};N\times D}
 \le C^2(\beta_n-\beta_{n+1})^{-2}
       \|\widehat G_n\|_{k,\beta_n;N\times D}^{\,2}
\]
and $\beta_n-\beta_{n+1}\ge
\tfrac12(\alpha-\beta)e^{-\lambda\kappa^n}$ for $n\ge1$:
\begin{align*}
    &\| T(G_n, \Lambda_\Psi) \| _{k,\beta_n; N \times \partial D} \leq \exp (-6 \lambda (\kappa^{n}-1) )\| \partial \Psi\| _{k,\alpha; N' \times \partial D}, \quad n\geq 0,\\
    &\| \widehat{G}_{n}\| _{k,\beta_n; N \times D} \leq \exp (-6 \lambda \big(\kappa^{n-1}-1) + \lambda/2 \big) \| \partial \Psi\| _{k,\alpha; N' \times \partial D}, \quad n\geq 1.
\end{align*}
By the same recurrence, we have $\|\widehat G_{n+1}\|_{k,\beta_{n+1}}\le \tfrac12\|\widehat G_n\|_{k,\beta_n}$, while $\|\widehat G_1\|_{k,\beta_1}\le \frac{{64}C}{{49}(\alpha-\beta)^{2}}\|\partial\Psi\|_{k,\alpha;N'\times\partial D}$. In particular, the sum of the corrections remains within the prescribed small neighborhood. The above estimates imply that $(G_n)$ converges in $\F^{k,\beta}$ and $T(G_n, \Lambda_\Psi) \rightarrow 0$ in $\D^{k,\beta}$ as $ n \rightarrow \infty$. Let $G = \lim G_n$ in $ \F^{k,\beta}$. Since the operator $T$ is continuous on $\F^{k,\beta}$, we have $T (G, \Lambda_\Psi) = 0$. Moreover, $G = (z, \xi)$ satisfies the following estimates:
\begin{align*}
    \| (z - z_0, \xi-\xi^0)\| _{k,\beta; N \times D} \leq C\sum_{n=1}^{\infty} \| \widehat{G}_n\| _{k,\beta_n; N \times D} \leq C_*(\alpha-\beta)^{-2} \| \partial \Psi\| _{k,\alpha; N' \times \partial D}.
\end{align*}
Here $C_*$ absorbs the fixed embedding constants. Smallness in a fixed sufficiently high finite-order norm, together with \(\Psi\in C^\infty\), is enough to obtain a smooth \(G\). The linearized boundary problem for each disc is invertible in a fixed Hölder space. Boundary regularity and the implicit function theorem then imply smoothness in both the disc variable and the spatial parameters.
\end{proof}

\subsection{A Patching Lemma}\label{subsecpatthm} From the previous subsection, we know that the HCMA disc family exists in each local chart. In this subsection, we use a patching argument to complete the full proof of Theorem \ref{thmexistbypertb}.

Let {\(\{N_i'\}_{i\in\mathcal I}\)} be {the} {holomorphic} {charts covering \(N\), and choose the inner charts \(N_i\) with the uniform coordinate margin \(R_0\) from Section~\ref{sec:notations}. Thus \(N=\bigcup_iN_i\) and \(\overline N_i\subset N_i'\). Decrease the smallness constant in Theorem~\ref{thmexistbypertb}, if necessary, so that the} local {disc displacement is smaller than this margin after the uniform} coordinate {comparisons.}

For each chart $N_i'$, let $\rho_i$ denote the local potential of the Kähler metric. We consider the local family of holomorphic discs $G_i: N_i \times D \to E|_{N_i'}$ satisfying the boundary value problem defined by the local potential $\rho_i$ and the perturbation $\Psi$.

\begin{lem} \label{lem:patching}
Suppose that for each {\(i\in\mathcal I\)} there {is} {an} HCMA disc {family} {\(G_i:N_i\times D\to E|_{N_i'}\)}, with {boundary} on {\(\Lambda_\Psi\), satisfying \eqref{estimeasureperturb}}. Assume {the smallness constant is chosen as above and is decreased, if necessary, so} that {the} {transformed} {families} {lie} {in} the {uniqueness} {neighborhoods} {of Theorem~\ref{thmlocalunifixpotential}}.
Let $G_i(w, \tau) = (z_i(w, \tau), \xi^i(w, \tau))$ denote the representation of the solution in the local coordinates of $E|_{N_i'} \cong N_i' \times \mathbb{C}^n $.
Then, for any two indices $i, j \in \mathcal{I}$ with $N_i \cap N_j \neq \emptyset$, the local solutions satisfy the following compatibility conditions on the overlap $(N_i \cap N_j) \times D$:
\begin{itemize}
    \item The base maps coincide:
    \begin{align*}
        z_i(w, \tau) = z_j(w, \tau).
    \end{align*}
    \item The fiber components are related by the transition of fibre bundle $\W$:
    \begin{align*}
        \xi^{j}(w, \tau) = \xi^{ i}(w, \tau) + \partial \left( \rho_j - \rho_i \right)(z_i(w, \tau)).
    \end{align*}
\end{itemize}

\end{lem}

\begin{proof}
Fix {\(i,j\)} {with \(N_i\cap N_j\ne\varnothing\). The uniform margin and the displacement bound imply
\[
z_i(w,D)\subset N_i'\cap N_j',
\qquad w\in N_i\cap N_j.
\]
Indeed, the coordinate displacement bounds the metric displacement uniformly, while every point of \(N_j\) is a fixed positive distance from the complement of \(N_j'\). Hence the full bundle transition is defined along every} such {leaf.} {Expressing} {both} {sides} in the {same} {base} {coordinates,} define {\(\widetilde G_j=(\tilde z_j,\tilde\xi_j)\)} {by}
\begin{align}\label{eq:transition_candidate}
    \tilde{z}_j(w, \tau) &= z_i(w, \tau), \nonumber \\
    \tilde{\xi}_{j}(w, \tau) &= \xi^{i}(w, \tau) + {\partial} (\rho_j - \rho_i)(z_i(w, \tau)),
\end{align}
for $(w, \tau) \in (N_i \cap N_j) \times D$. We claim that $\widetilde{G}_j$ defines a smooth family of holomorphic discs parameterized by $N_i \cap N_j$ satisfying the boundary condition determined by $\Lambda_\Psi$.

Since $z_i$ is holomorphic in $\tau$, $\tilde{z}_j$ is trivially holomorphic. We check $\tilde{\xi}_j$. Let $\rho_{ij} = \rho_j - \rho_i$. Applying the operator $\overline{\partial}_\tau$:
\begin{align*}
    \overline{\partial}_{\tau} \tilde{\xi}_{j, p} &= \overline{\partial}_{\tau} \xi^{i, p} + \overline{\partial}_{\tau} \left[ \frac{\partial \rho_{ij}}{\partial z_p} (z_i(w, \tau)) \right] \\
    &= 0 + \sum_{q} \frac{\partial^2 \rho_{ij}}{\partial z_q \partial z_p} \frac{\partial z_{i, q}}{\partial \bar{\tau}} + \sum_{q} \frac{\partial^2 \rho_{ij}}{\partial \bar{z}_q \partial z_p} \frac{\partial \bar{z}_{i, q}}{\partial \bar{\tau}}.
\end{align*}
The first term in the sum vanishes because $z_i$ is holomorphic. The second term vanishes because $\rho_i$ and $\rho_j$ are potentials for the same K\"ahler metric $\omega$, implying that their difference $\rho_{ij}$ is pluriharmonic. Thus, $\overline{\partial}_\tau \tilde{\xi}_j = 0$.

On the boundary $\partial D$, $G_i$ satisfies $\xi^{i} (w, \tau) = \partial (\rho_i + \Psi)(z_i(w, \tau), \tau)$. Substituting this into (\ref{eq:transition_candidate}):
\begin{align*}
    \tilde{\xi}_{j} (w, \tau) &= {\partial (\rho_i + \Psi)} (z_i(w, \tau), \tau) + {\partial (\rho_j - \rho_i)} (z_i (w, \tau)) \\
    &= {\partial (\rho_j + \Psi)} (\tilde{z}_j(w, \tau), \tau).
\end{align*}
This is precisely the boundary condition required for the solution in the $j$-th chart.

The {chain} {rule,} {the full coordinate transition, and the uniform chart and potential bounds imply
\[
\|\widetilde G_j-G_{0,j}\|_{k,\beta;(N_i\cap N_j)\times D}
\leq C\|G_i-G_{0,i}\|_{k,\beta;N_i\times D}.
\]
Here the potential terms can be estimated using differences of \(\partial\rho_i\) and \(\partial\rho_j\), so their second-derivative bounds suffice. By the choice of the smallness constant, the right-hand side} is {less} {than} {\(\delta_1\), as is the corresponding norm} of {\(G_j-G_{0,j}\).} {Both} {families} {therefore} {lie} in the {same} {uniqueness} {neighborhood.}

By local uniqueness, Theorem \ref{thmlocalunifixpotential}, there is at most one such solution in the neighborhood. Therefore, $\widetilde{G}_j$ must coincide with the solution $G_j$ constructed directly in the chart $N_j$. Hence, the local solutions glue together consistently.
\end{proof}

We are now ready to complete the proof of Theorem \ref{thmexistbypertb}.

\begin{proof}[Proof of Theorem \ref{thmexistbypertb}]Let $\{N_i\}_{i \in \mathcal{I}}$ be the open covering of $N$ constructed above. By the local existence result established in the previous subsection (via the Nash-Moser iteration in local charts), for each $i \in \mathcal{I}$, there exists a smooth family of holomorphic discs $G_i: N_i \times D \to E|_{N_i'}$ satisfying the boundary condition determined by the local data $(N_i, N_i', \rho_i, \Psi)$. Consider any two indices $i, j \in \mathcal{I}$ such that $N_i \cap N_j \neq \emptyset$. By Lemma \ref{lem:patching}, the local solutions $G_i$ and $G_j$ satisfy the compatibility conditions on the overlap $(N_i \cap N_j) \times D$. Specifically, their base components coincide, and their fiber components transform according to the transition functions of the bundle $\mathcal{W}$.

Consequently, the collection of local families $\{G_i\}_{i \in \mathcal{I}}$ patches together to define an HCMA disc family $G: N \times D \to \mathcal{W}$ satisfying the boundary condition determined by $\Lambda_\Psi$.
\end{proof}

\section{Construction of \texorpdfstring{$\Omega$}{Ω}-psh subsolutions}\label{secmaxrkthm}

In this section, we construct an $\Omega$-psh subsolution to the HCMA equation with the boundary data, $\Psi$, satisfying \eqref{con:smallperturb_bdry}, and complete the proof of the main Theorem. The process of construction relies on the following two points:
\begin{itemize}
    \item the existence of an HCMA disc family on $N$ by Theorem \ref{thmexistbypertb};
    \item the existence of convex K\"ahler potentials for the reference K\"ahler form $\omega$ in holomorphic coordinate charts covering $N$.
\end{itemize}

\subsection{Local Construction}

Recall from Section \ref{sec:notations} that we fix a pair of nested open domains $\widetilde{N} \Subset N \subset X$, covered by a locally finite atlas with nested coordinate charts $\widetilde{N}_i \Subset N_i$ for $i \in \mathcal{I}$, satisfying uniform geometric bounds. We denote by $R_0 > 0$ the uniform lower bound for the difference between the outer and inner radii of the chart images corresponding to $N_i$ and $\widetilde{N}_i$. {The enlarged charts $N_i'$ below are chosen within the original controlled charts. The} potentials satisfy uniform convexity and regularity estimates {there,} controlled by constants $\lambda$ and $C$; specifically, $\nabla^2 \rho_i \ge \lambda I_{2n}$ and {$\sup_{i \in \mathcal{I}} \|\nabla^2 \rho_i\|_{k+2,\alpha; N_i'} \le C$}.

We first recall the existence of the HCMA disc families on $N$ established in Theorem \ref{thmexistbypertb}. For each point $x \in N$, we denote the associated holomorphic disc by the map $G(x, \cdot): {D} \to \mathcal{W}$, which lifts the curve in the base manifold. In local coordinates, we write $G(x, \tau) = (z(x, \tau), \xi(x, \tau))$. We define the leaf of the holomorphic disc at $x \in N$ as the set $$\mathcal{L}_x \coloneqq \big\{ (z(x, \tau), \tau) \in N' \times \bar{D} : \tau \in \bar{D} \big\}.$$
To rigorously control the geometry of these disc families and avoid boundary issues, we expand upon the base domains $\widetilde{N} \Subset N$ introduced in Section  \ref{sec:notations}. Write $h_\tau(x):=z(x,\tau)$. We introduce a three-tiered hierarchy of nested open sets. For the outer domain, we define slightly shrunken and slightly enlarged regions such that $N'' \Subset N \Subset N'$ satisfying $N'' \subseteq h_\tau(N) \subseteq N' \quad \text{for all } \tau \in \bar{D}.$ Similarly, for the inner domain, we define $\widetilde{N}'' \Subset \widetilde{N} \Subset \widetilde{N}'$ such that $\widetilde{N}'' \subseteq h_\tau(\widetilde{N}) \subseteq \widetilde{N}' \quad \text{for all } \tau \in \bar{D}.$
Corresponding to our locally finite atlas, we define associated nested coordinate charts for each $i \in \mathcal{I}$:\begin{itemize}\item For the intermediate layers: $N''_i \Subset N_i \Subset N'_i$.\item For the inner layers: $\widetilde{N}''_i \Subset \widetilde{N}_i \Subset \widetilde{N}'_i$.\end{itemize}These safety buffers ensure that for any base point inside a given chart, its corresponding holomorphic disc leaf does not escape the enlarged boundary. 

With these notations, we proceed to construct a family of real functions, $\{L_{x_0}\}$, indexed by the base point $x_0 \in \widetilde{N}$, which will serve as local building blocks for our $\Omega$-psh subsolution. For a fixed base point $x_0$ contained in a chart $\widetilde{N}_i$, the function $L_{x_0}$ is defined locally on the product domain $N_i' \times D$. We define $L_{x_0}$ to be affine in the spatial coordinates $z$ (relative to the holomorphic chart $N_i'$) such that its spatial gradient matches the fiber coordinate $\xi$ of the lifted disc. Explicitly, $L_{x_0}$ satisfies the following conditions:
\begin{align}\label{equ:defineLx0}
\begin{split}
\frac{\partial}{\partial \tau} \frac{\partial}{\partial\bar{\tau}} \big\{ L_{x_0} (z(x_0, \tau), \tau) \big\}  = 0, \quad \qquad  &\tau \in D; \\
L_{x_0} (z(x_0, \tau), \tau) = (\rho_i+\psi_\tau)(z(x_0, \tau)), \qquad  &\tau \in \partial D; \\
\partial_{z_p} L_{x_0} (z, \tau) =  \xi^p (x_0, \tau), \hspace{1.3cm} & (z,\tau) \in N_i' \times  D. 
\end{split}
\end{align}
In the following lemma, we will verify that $L_{x_0}$ defines a real pluriharmonic function on $N_i' \times D$.

\begin{lem} \label{lempshLcheck} Let $L_{x_0}$ be the function defined as in (\ref{equ:defineLx0}). Then, $L_{x_0}$ is a smooth pluriharmonic function in $N_i'\times D$.
\end{lem} 
\begin{proof} We can write down $L_{x_0}$ explicitly in the local holomorphic coordinates as a family of affine functions varying with respect to $\tau$,
\begin{align} \label{explicitL}
    L_{x_0} (z, \tau) = 2 \sum_{p=1}^n \re \big\{ \xi^p(x_0,\tau) (z_p - z_p(x_0, \tau)) \big\} + L_{x_0} (z(x_0, \tau), \tau), \quad (z , \tau) \in N_i' \times D.
\end{align}
 By taking derivatives for $L_{x_0} (z,\tau)$, we have
 \begin{align}
         (\partial_{z_p} L_{x_0})(z, \tau) = & \ \xi^p (x_0, \tau),\label{1stderivofLspace}\\
         \begin{split} \label{1stderivofLtime}
         (\partial_\tau L_{x_0}) (z, \tau) = & \sum_{p =1}^n \partial_{\tau} \xi^p(x_0, \tau) (z_p - z_p(x_0, \tau))  -  \sum_{p=1}^n  \xi^p (x_0,\tau) \frac{\partial z_p}{\partial \tau}  (x_0, \tau)   \\& + \frac{\partial}{\partial \tau} \big\{ L_{x_0} (z(x_0, \tau), \tau) \big\}.
         \end{split}
 \end{align}
 The smoothness of $L_{x_0}$ is obvious.
By (\ref{1stderivofLspace}), and noting that $\xi^p$ is holomorphic with respect to $\tau$, we can compute the second derivatives of $L_{x_0}$,
\begin{align*}
    (\overline{\partial}_{z_q}\partial_{z_p}L_{x_0}) (z, \tau)  &=0,\\
    (\overline{\partial}_{\tau} \partial_{z_p}L_{x_0}) (z, \tau) & = 0.
\end{align*}
It suffices to check $\overline{\partial}_\tau \partial_\tau L_{x_0} =0$. By taking the derivative of (\ref{1stderivofLtime}), we have
\begin{align}
\begin{split}\label{pshttleaf}
    \big(\overline{\partial}_\tau \partial_\tau L_{x_0}\big)(z, \tau) 
    =& \sum_{p=1}^n \overline{\partial}_\tau \partial_{\tau} \xi^p(x_0,\tau) (z_p - z_p(x_0, \tau))  - \sum_{p=1}^n \partial_\tau \xi^p(x_0,\tau) \frac{\partial z_p}{\partial \overline{\tau}}  (x_0, \tau)
    \\& -  \sum_{p=1}^n \frac{\partial}{\partial \overline{\tau}} \big\{ \xi^p (x_0, \tau)  \frac{\partial z_p}{\partial \tau} (x_0, \tau)  \big\} 
    + \frac{\partial}{\partial \overline{\tau}}\frac{\partial}{\partial \tau} \big\{ L_{x_0} (z(x_0, \tau), \tau) \big\}. 
    \\ = & \frac{\partial}{\partial \overline{\tau}}\frac{\partial}{\partial \tau} \big\{ L_{x_0} (z(x_0, \tau), \tau) \big\} =0.
\end{split}
\end{align}
Hence, we complete the proof that $L_{x_0}$ is pluriharmonic.
\end{proof}

We now establish the estimates for the local pluriharmonic functions $L_{x_0}$. Recall that for each chart index $i \in \mathcal{I}$, the difference in radii between the base outer domains $N_i$ and the inner domains $\widetilde{N}_i$ is bounded from below by $R_0$. Let $R$ and $r$ denote the radii of the base outer and inner coordinate balls, respectively, so that $R-r\ge R_0$. We choose the disc displacement and the collar width so that their sum is at most $R_0/8$.

Let $\lambda > 0$ denote the uniform convexity constant of the reference potentials $\rho_i$, and let $M > 0$ be the $L^\infty$ bound of the solution $\Phi$. To proceed with the proof, we require the domain size and convexity to be sufficiently large relative to this $L^\infty$ bound. Specifically, we assume the condition:
\begin{align} 
\label{condi:relationparas}
\lambda R_0^2 \geq 16M.
\end{align}

\begin{lem} \label{lemsubsolucon}
Fix an index $i \in \mathcal{I}$ and the associated nested domains $\widetilde{N}_i \subset N_i'$. For any base point $x_0 \in \widetilde{N}_i$, the local pluriharmonic function $L_{x_0}$ defined on $N_i' \times D$ satisfies the following estimates:
\begin{align}
L_{x_0} (z, \tau) - \rho(z) &\leq \psi_\tau (z),  \hspace{1cm} (z, \tau) \in N_i' \times \partial D, \label{subsoluconbdry}\\
L_{x_0} (z, \tau) - \rho(z) &\leq -2M, \hspace{1cm} (z, \tau) \in (\partial N_i')_{\delta} \times D. \label{subsoluconcpt}
\end{align}
where $(\partial N_i')_{\delta}$ is the $\delta$-neighborhood of the boundary relative to the chart, defined as the set of points in $N_i'$ within Euclidean distance $\delta$ of $\partial N_i'$.
\end{lem}

\begin{proof}The estimates follow from the strict convexity of the potential $\rho$. Note that $\mathcal{L}_{x_0}$ be the leaf of the HCMA disc passing through $x_0$. For a fixed $\tau \in \partial D$, let $z_{0} = z(x_0, \tau)$ denote the spatial position of the leaf on the boundary.

Recall the construction of $L_{x_0}$ given in (\ref{equ:defineLx0}). Since $L_{x_0}(\cdot, \tau)$ is affine in the spatial coordinates $z$ and satisfies the matching conditions for value and gradient at $z_0$ along the leaf $\mathcal{L}_{x_0}$, the function $z \mapsto L_{x_0}(z, \tau)$ represents the tangent plane to the function $\rho + \psi_\tau$ at $z_0$ for $\tau \in \partial D$. Specifically, the difference function
$$F_\tau(z) \coloneqq L_{x_0}(z, \tau) - \rho(z) - \psi_\tau(z), \qquad \tau \in \partial D$$
satisfies $F_\tau(z_0) = 0$ and $\partial_z F_\tau(z_0) = 0$.

Using the uniform convexity of the potential $\rho$ (Eq. \ref{est:convexity}) and the linearity of $L_{x_0}$ in $z$, the Hessian of the difference is strictly negative definite:
\begin{align*}
    D^2_z F_\tau(z, \tau) = -D^2_z \big(\rho(z)+\psi_\tau(z)\big) \leq -\frac{\lambda}{2} I_{2n}.
\end{align*}
Then, we obtain the quadratic decay estimate for $\tau \in \partial D$:
\begin{align} \label{est:taylor_bound}
    L_{x_0}(z, \tau) - \rho(z) - \psi_\tau(z) \leq -\frac{\lambda}{4} |z - z_0|^2.
\end{align}
Since the quadratic term is non-positive, this immediately implies the boundary estimate (\ref{subsoluconbdry}):
$$L_{x_0}(z, \tau) - \rho(z) \leq \psi_\tau(z) \quad \text{for all } (z, \tau) \in N_i' \times \partial D.$$

Next, we establish the interior estimate (\ref{subsoluconcpt}). Consider a point $z \in (\partial N_i')_{\delta}$. By our geometric setup, the base point $x_0$ lies in the inner domain $\widetilde{N}_i$, while $z$ lies near the boundary of the outer domain $N_i'$. Since the holomorphic discs are small perturbations, the leaf point $z_0 = z(x_0, \tau)$ remains within a small neighborhood of $B_r$. Therefore, by the Uniform Size assumption ($R - r \ge R_0$) and our choice of disc displacement and collar width, we have:
$$|z-z_0|\ge R_0-\frac{R_0}{8}=\frac{7R_0}{8}.$$
Substituting this into (\ref{est:taylor_bound}) and using the bound $\|\psi\|_{L^\infty} \le M$, we have for any $\tau \in \partial D$:
\begin{align*}
    L_{x_0}(z, \tau) - \rho(z) &\leq M - \frac{49\lambda}{256} R_0^2.
\end{align*}
Invoking the assumption $\lambda R_0^2 \geq 16M$, we obtain:
\begin{align*}
    L_{x_0}(z, \tau) - \rho(z) &\leq M-\frac{49}{16}M=-\frac{33}{16}M\leq-2M.
\end{align*}
Thus, the estimate holds on the boundary of the disc parameter space $\partial D$.

To extend this to the interior $\tau \in D$, fix the spatial point $z \in (\partial N_i')_{\delta}$. By Lemma \ref{lempshLcheck}, $L_{x_0}(z, \cdot)$ is pluriharmonic in $N_i'\times D$, while $\rho(z)$ is constant. Therefore, the function $\tau \mapsto L_{x_0}(z, \tau) - \rho(z)$ is harmonic on $D$. Applying the Maximum Principle, we conclude:
\begin{align*}
    L_{x_0}(z, \tau) - \rho(z) \leq \sup_{\zeta \in \partial D} \big( L_{x_0}(z, \zeta) - \rho(z) \big) \leq -2M, \qquad \text{for all } \tau \in D.
\end{align*}
This proves (\ref{subsoluconcpt}).
\end{proof}

Having established the properties of the individual building blocks, we now proceed to construct a local $\Omega$-plurisubharmonic function on the product domain $N_i' \times D$. We define the function $F_i$ as the upper envelope of the family $\{L_{x}\}_{x \in \widetilde{N}_i}$ relative to the local potential $\rho_i$. Specifically, we set:
\begin{align} \label{def:Fi}
F_i (z, \tau) = \sup_{x \in \widetilde{N}_i} \big( L_{x, i} (z, \tau) - \rho_i (z) \big), \qquad (z, \tau) \in N_i' \times D,
\end{align}
where $L_{x, i}$ denotes the pluriharmonic function $L_x$ constructed in Lemma \ref{lemsubsolucon} using the chart potential $\rho_i$.

\begin{lem} \label{lemconstructpshlocal}
The function $F_{i}$ defined in (\ref{def:Fi}) is a continuous $\Omega$-plurisubharmonic function on the local domain $N_i' \times D$. Furthermore, for any $x \in \widetilde{N}_i$, the function $F_i$ coincides with its generator along the leaf of the holomorphic disc $\mathcal{L}_x$. That is,
\begin{align*}
F_{i}(z(x, \tau), \tau) = L_{x, i} (z(x, \tau), \tau) - \rho_{i} (z(x, \tau)), \qquad \tau \in D.
\end{align*}
\end{lem}

\begin{proof}
For the first statement, it suffices to prove that $F_{i} (z, \tau)$ is continuous on $N'_i \times D$. By Theorem \ref{thmexistbypertb}, there exists a constant $M_i$ depending on $N_i$ such that $|\xi^p(x, \tau)| \leq M_i$ for all $x \in \widetilde{N}_i$. In view of formula \eqref{1stderivofLspace}, this implies $|\partial_{z_p} L_{x, i} (z, \tau)| \leq M_i$ for all $(z, \tau) \in N'_i \times D$. Similarly, formula \eqref{1stderivofLtime} ensures that $|\partial_\tau L_{x, i}|$ is uniformly bounded by a constant depending on $N_i$ and the diameter of $N'_i$. These uniform bounds imply that the family of smooth pluriharmonic functions $\{L_{x,i} : x \in \widetilde{N}_i\}$ is equicontinuous on $N'_i \times D$. Consequently, the envelope function $F_i(z, \tau) = \sup_{x \in \widetilde{N}_i} L_{x, i} (z, \tau)-\rho_i(z)$ is continuous.

For the second statement, fix $x_0 \in \widetilde{N}_i$. We aim to show that
\begin{align*}
    \sup_{x \in \widetilde{N}_i} {L_{x,i} (z(x_0, \tau), \tau)} = L_{x_0, i} (z(x_0, \tau), \tau).
\end{align*}
Consider the restriction to the holomorphic disc centered at $x_0$. The function $\tau \mapsto L_{x,i}(z(x_0, \tau), \tau)$ is harmonic with respect to $\tau$ for $\tau \in D$. By Lemma \ref{lemsubsolucon} and equation \eqref{equ:defineLx0}, we have
\begin{align*}
    L_{x, i} (z(x_0, \tau), \tau) \leq (\rho_i+ \psi_\tau)(z(x_0, \tau)) = L_{x_0, i} (z(x_0, \tau), \tau)
\end{align*}
for all $\tau \in \partial D$ and $x \in \widetilde{N}_i$. By the classical maximum principle for harmonic functions, this inequality holds for all $\tau \in D$, i.e., $L_{x, i} (z(x_0, \tau), \tau) \leq L_{x_0, i} (z(x_0, \tau), \tau)$, which completes the proof.
\end{proof}

\subsection{Patching and Global Construction}
\label{secconstructF}
Lemma \ref{lem:patching}, the patching lemma for HCMA disc families, ensures that the leaves of the holomorphic discs coincide on the intersection $N_i \cap N_j$. We now examine the relationship between $F_i$ and $F_j$ on the intersection $N'_i \cap N'_j$.

\begin{lem} \label{lempshpatching}
Let $x \in \widetilde{N}_i$. For any index $j \in \mathcal{I}$ such that $x \in N_j$, we have\begin{align*}F_i (z(x,\tau), \tau) \geq F_j (z(x, \tau), \tau).\end{align*}Furthermore, if $x \in \widetilde{N}_j$, then
\begin{align*}F_i (z(x,\tau), \tau) = F_j (z(x, \tau), \tau).\end{align*}
\end{lem}

\begin{proof}We first prove the case when $x \in \widetilde{N}_j$. By Lemma \ref{lemconstructpshlocal}, we have the identity $F_{i} (z(x, \tau), \tau) = L_{x, i} (z(x,\tau), \tau) - \rho_i (z(x, \tau))$. Consider the function $u(z(x, \tau), \tau) = F_i(z(x, \tau), \tau) + \rho_j(z(x, \tau))$ defined on $(N_i \cap N_j) \times D$. Restricting $u$ to the leaf of the holomorphic disc centered at $x$, we compute the Laplacian with respect to $\tau$:
\begin{align*}
\partial_\tau \overline{\partial}_\tau \big\{ F_i (z(x, \tau), \tau) + \rho_j(z(x, \tau)) \big\}
= &~\partial_\tau \overline{\partial}_\tau \big\{ L_{x, i} (z(x, \tau), \tau) \big\} 
\\
&+ \partial_{z_p} \overline{\partial}_{z_{q}} (\rho_{j} - \rho_{i}) (z(x, \tau)) \frac{\partial z_{p}}{\partial \tau} \frac{\partial \bar{z}_{q}}{\partial \bar{\tau}}.
\end{align*}
The first term on the right-hand side vanishes because $L_{x, i}$ is harmonic along the leaf (by construction). The second term also vanishes because the complex Hessians coincide, i.e., $i\partial \overline{\partial} \rho_i = i\partial \overline{\partial} \rho_j$ on $N'_i \cap N'_j$. Consequently, $F_{i} + \rho_j$ is harmonic on the leaf of the holomorphic disc centered at $x$. On the boundary of this leaf, we observe that
\begin{align*}
F_i (z(x, \tau), \tau) + \rho_j(z(x, \tau)) = L_{x, j} (z(x, \tau), \tau), \qquad \text{for } \tau \in \partial D.\end{align*}Thus, by the uniqueness of harmonic extensions, we have $F_i + \rho_j = L_{x,j}$ on the leaf. Since $x \in \widetilde{N}_i \cap \widetilde{N}_j$, it follows that $F_i = L_{x,j} - \rho_{j} = F_j$ on the leaf.

Next, consider the case when $x \in N_j$ (but not necessarily in $\widetilde{N}_j$). The function defined above, $F_i (z(x, \tau), \tau) + \rho_j(z(x, \tau))$, remains harmonic with respect to $\tau$. Since $x\in \widetilde{N}_i$, we have the boundary value:
\begin{align*}
    F_{i} (z(x, \tau), \tau) + \rho_{j}(z(x, \tau)) = (\rho_j + \psi_\tau) (z(x, \tau)) \quad \text{for } \tau \in \partial D.
\end{align*}
By Lemmas \ref{lempshLcheck} and \ref{lemsubsolucon}, for each $y \in \widetilde{N}_j$, we have the inequality
\begin{align*}
    F_{i} (z(x, \tau), \tau) + \rho_{j}(z(x, \tau)) \geq L_{y,j} (z(x, \tau), \tau).
\end{align*}
Taking the supremum over all $y \in \widetilde{N}_j$, we obtain $F_i(z(x, \tau), \tau) + \rho_j(z(x, \tau)) \geq F_j(z(x, \tau), \tau) + \rho_j(z(x, \tau))$, which implies $F_i(z(x, \tau), \tau) \geq F_{j}(z(x, \tau), \tau)$. This completes the proof.
\end{proof}

We now {construct} a global $\Omega$-plurisubharmonic subsolution. {Choose} {the} {outer} {charts $\{N_i'\}_{i\in\mathcal I}$} locally finite {in} {$X$.} {Replace each $F_i$ on $N_i'\times D$ by $\max\{F_i,-M\}$ and extend it by $-M$ outside this chart, retaining} the {notation} {$F_i$.} {By} {\eqref{subsoluconcpt},} {this} {extension is constant near the chart boundary} and {is} {continuous} {and} {$\Omega$-plurisubharmonic.} {Define}
\begin{align}\label{constructdim2}
F(z, \tau) = \max \big\{\sup_{i \in \mathcal{I}} F_i (z, \tau), -M\big\}, \qquad (z, \tau ) \in X \times D.
\end{align}
{Local} {finiteness} {of the outer charts implies} that $F$ is a {continuous} $\Omega$-plurisubharmonic {function} on {$X\times D$}.

We are now in a position to state the main theorem of this section. Recall that $(X, \omega)$ is a complete Kähler manifold and that $\Phi$ is defined as the upper envelope of the class of subsolutions $\mathcal{B}_{\Omega, \Psi}$. We also recall that $\Phi$ admits a uniform $L^\infty$ bound.
\begin{thm}\label{thmconstpshsolu}
There exists a globally defined function $F$ on $X \times D$ satisfying the following properties:
\begin{itemize}
\item[(i)] $F \in \mathcal{B}_{\Omega, \Psi}$ is continuous on $X \times D$ and smooth on {$\widetilde{N}'' \times D$}.
\item[(ii)] Let $\Phi$ be the unique bounded continuous solution to the HCMA equation given by the upper envelope of $\mathcal{B}_{\Omega, \Psi}$. Then $F = \Phi$ on {$\widetilde{N}'' \times D$}.
\end{itemize}
\end{thm}

\begin{proof} 
The continuity of $F(z, \tau)$ follows directly from Lemmas \ref{lemsubsolucon}, \ref{lemconstructpshlocal}, and \ref{lempshpatching}, as established in the construction. Thus, $F(z,\tau)$ is a continuous $\Omega$-psh function on $X \times D$. 

It remains to check the boundedness of $F$ on $X \times D$. It is clear from the definition \eqref{constructdim2} that $F(z, \tau) \geq -M$, so it is bounded from below. To show that $F$ is bounded from above, observe that for each fixed $z \in X$, the function $\tau \mapsto F(z, \tau)$ is subharmonic on $D$. By Lemma \ref{lemsubsolucon}, we have the boundary condition $F(z, \tau) \leq \psi_\tau(z)$ for $\tau \in \partial D$. Applying the maximum principle to the subharmonic function in $\tau$, we obtain
\begin{align*}
    F(z, \tau) \leq \sup_{(z, \tau) \in X\times \partial D} \Psi(z, \tau).
\end{align*}
Thus, $F \in \mathcal{B}_{\Omega, \Psi}$. Since $\Phi$ is the upper envelope of $\mathcal{B}_{\Omega, \Psi}$, it follows immediately that
\begin{align*}
    F(z, \tau) \leq \Phi(z, \tau) \quad \text{for all } (z, \tau) \in X \times D.
\end{align*}

We now prove {equality} on {$\widetilde{N}''\times D$}. Let {$x\in\widetilde N_i$}. {On} the leaf $\mathcal L_x$, the function {$\Phi+\rho_i$} is {subharmonic} {and} {has} boundary {values} {at} {most} {those} {of} the harmonic function {$L_{x,i}$,} {by} {\eqref{equ:defineLx0}.} {The} {maximum principle,} the {definition} {of} {$F_i$,} {and} {its} {truncation} {therefore} {imply
\begin{align*}
\Phi(z(x,\tau),\tau)
&\leq L_{x,i}(z(x,\tau),\tau)-\rho_i(z(x,\tau))\\
&\leq F_i(z(x,\tau),\tau)\leq F(z(x,\tau),\tau)\leq\Phi(z(x,\tau),\tau).
\end{align*}
Thus} {equality} {holds along every} leaf {labelled by $\widetilde N$}. Since {$\widetilde N''\subseteq h_\tau(\widetilde N)$}, {it} {follows} {that $F=\Phi$} on {$\widetilde N''\times D$}. {In the smooth foliation coordinates, this common value is smooth by the harmonic extension in \eqref{equ:defineLx0}, proving the remaining assertion in (i).}
\end{proof}

\section{Proof of the Main Theorem and Applications}

\subsection{Proof of Theorem \ref{mthmlocalreg}}

The proof follows by combining several results established in the preceding sections. {Fix a chart $(V_i,U_i,b_i)$ from Section~\ref{sec:notations} and work in its coordinates. Set $N''=B_r(0)=b_i(V_i)$, $N=B_{r+R_0/32}(0)$, and $N'=B_R(0)$. For sufficiently small $\varepsilon$, the disc family below satisfies $N''\subset z(N,\tau)\subset N'$ for every $\tau\in\overline D$.} By the smallness condition and Theorem~\ref{thmexistbypertb}, there exists a family of holomorphic discs $G: N\times D \rightarrow T^* N'$. If we write $G (w, \tau)= (z(w, \tau), \xi(w,\tau))$, where $w$ denote the standard holomorphic coordinates on $N'$ and $\xi$ denote the complex bundle coordinates of $T^* N'$, then we have the following displacement estimate:
\begin{align*}
    \| (z-z_0, \xi-\xi^0)\| _{k, \beta; N \times D} \leq C \| \partial \Psi\| _{k, \alpha; N' \times \partial D},
\end{align*}
where $G_0(w, \tau)= (z_0, \xi^0)(w, \tau) $ is the trivial foliation given by $(z_0, \xi^0) (w, \tau) = (w, \partial \rho (w))$. After projection via $\pi$, the family of holomorphic discs $G$ induces a holomorphic discs foliation on {its image.} Let $H(w, \tau) =\pi \circ G (w, \tau) $. Then, for each $x_0 \in N$, the corresponding leaf of foliation is given by $\mathcal{L}_{x_0} = \{(H(x_0, \tau),\tau); ~\tau \in D\}$.

A global subsolution $F \in \B_{\Omega, \Psi}$ can be constructed based on the data $(G, \Psi)$, following the construction developed in Section \ref{secconstructF}. Precisely, for each $x_0$, there exist a real pluriharmonic function $L_{x_0}$ (see Lemma \ref{lempshLcheck}) on $N' \times D$ defined as follows:
\begin{itemize}
    \item Along the leaf $\LL_{x_0}$, $L_{x_0}$ is harmonic in $\tau$ and matches the boundary data $\tilde{\psi}_\tau = \rho + \psi_\tau$.
    \item In the spatial direction of $N' \times D$, $L_{x_0}$ is affine in the spatial coordinates $z_p$, with the derivative in each $z_p$ given by $\xi^{p}$.
\end{itemize}
If we choose $\varepsilon$ in \eqref{cond:localsmall} small enough, the Euclidean Hessian satisfies $D_z^2(\rho+\psi_\tau)\geq \frac{\lambda}{2}I_{2n}$. {Set $\widetilde N=B_{R-R_0/32}(0)$ and decrease $\varepsilon$ if necessary.} By the Taylor estimate in Lemma~\ref{lemsubsolucon}, we have, for $x\in N$, $z\in N'\setminus\widetilde N$, and $\tau\in\partial D$,
\[
 L_x(z,\tau)-\rho(z)
 \le M-\frac{\lambda}{4}|z-z(x,\tau)|^2
 \le M-\frac{49\lambda R_0^2}{256}
 \le -2M,
\]
where we used \eqref{thresholdcon}. The maximum principle extends this collar bound to every $\tau\in D$.
Then, the global subsolution $F$ can be defined as follows: 
\begin{align*}
    F (z, \tau) = \max \big\{-M, \sup_{x \in N}( L_x- \rho) (z, \tau) \big\}.
\end{align*}
{Extend $F$ by $-M$ outside $N'\times D$.} According to the proof of Theorem \ref{thmconstpshsolu} (general case), we have that $F \in \B_{\Omega, \Psi}$, and moreover, $F \equiv \Phi$ on {$N'' \times D$}.

It suffices to prove the estimate \eqref{estilocal}. The method follows the same strategy as in \cite[Section~9.1]{yao2025asymptotic}. In particular, we solve Dirichlet problems on the family of holomorphic discs defined by the foliation, as follows:
\begin{align*}
    \begin{cases}
        \Delta_\tau u (w, \tau) = 0, \qquad & \text{ in } N \times D; \\
        u (w, \tau) = f(w, \tau), & \text{ in } N \times \partial D.
    \end{cases}
\end{align*}
where $f$ is given by $f(w, \tau) = \Psi (z(w, \tau), \tau) + \rho(z(w, \tau)) - \rho(w)$. Using the notation $\mathfrak{e} = \rho(z(w, \tau)) - \rho(w)$, we have $$\| \mathfrak{e}\| _{k,\beta; N \times D} \leq \| z-z_0\| _{k,\beta; N \times D} \cdot \| \nabla \rho\| _{k,\beta; N \times D} \leq  C \| \partial\Psi\| _{k,\alpha; N' \times \partial D}.$$
By \cite[Lemma~2.8]{yao2025asymptotic} and \cite[Lemma~2.2]{yao2025asymptotic}, we have
\begin{align*}
 \| u\| _{k,\beta; N \times D} 
 &\leq C\| f\| _{k,\beta; N \times \partial D} \\
 & \leq C \| \Psi\| _{k+1, \alpha; N' \times \partial D}.
 \end{align*}
Note that along the leaf $\LL_x$ for $x \in N$, $F(z(w, \tau), \tau) = u(w, \tau) - \mathfrak{e}(w, \tau)$ for $(w, \tau) \in N \times D$. Since $\Phi \equiv F$ in {$N'' \times D$ and the maps $z(\cdot,\tau)$ have uniformly bounded inverses in $C^{k,\beta}$ on $N''$}, we obtain
\begin{align*}
    \| \Phi\| _{{k,\beta; N'' \times D}} \leq C \| \Psi\| _{k+1, \alpha; N'\times \partial D}.
\end{align*}
{The sets $b_i^{-1}(N'')=V_i$ cover $V$, giving \eqref{estilocal}.}

\subsection{Geometric applications of the local higher regularity theorem}
\label{subsec:geometric_applications}

We discuss several applications of Theorem~\ref{mthmlocalreg}. {Fix $k\ge1$ and $0<\beta<\alpha<1$.} Throughout this subsection, the smooth boundary data have positive K\"ahler slices, and $\Phi$ is the corresponding bounded weak solution of \eqref{introHCMA}, given by the upper envelope of $\mathcal{B}_{\Omega,\Psi}$ as in that theorem. Set $M=\|\Psi\|_{L^\infty(X\times\partial D)}$, the global $L^\infty$ norm of the boundary data. Then
\[
 -M\leq\Phi\leq M:
\]
the constant $-M$ belongs to $\mathcal{B}_{\Omega,\Psi}$, while the upper bound follows from the maximum principle on each vertical disc.

For the applications at a fixed geometric scale, choose nested holomorphic charts $V_p\Subset U_p$, centered at $p$, with coordinate maps $b_p:U_p\to B_{2a}$ satisfying $b_p(V_p)=B_a\Subset B_{2a}\subset\CC^n$. Assume that their local potentials satisfy $\nabla_0^2\rho_p\geq\lambda I$ and have uniform higher derivative bounds. Such charts can be chosen on a compact K\"ahler manifold, or on a complete K\"ahler manifold with positive injectivity radius and bounded curvature derivatives of every order: use the holomorphic bounded geometry of \cite[Theorem~9(2)]{wuyau2020invariant}, normalize the potentials by subtracting their pluriharmonic quadratic terms at the centers, and shrink the charts uniformly. All norms below are measured in these coordinates. By recentering at $x\in V_p$ and using $B_{a/2}(b_p(x))\Subset B_a(b_p(x))\subset B_{2a}$, we have the uniform margin $a/2$. If
\begin{equation}\label{eq:application_barrier}
 16M<\lambda a^2/4,
\end{equation}
then the barrier hypothesis holds uniformly. Consequently, sufficiently small $\|\Psi\|_{{k+3,\alpha;U_p\times\partial D}}$ implies smoothness on $V_p\times D$ and
\begin{equation}\label{eq:application_local_estimate}
 \|\Phi\|_{k,\beta;V_p\times D}
 \leq C\|\Psi\|_{k+1,\alpha;U_p\times\partial D},
 \qquad 0<\beta<\alpha,
\end{equation}
with $C$ independent of $p$. By applying the theorem on these relatively compact pairs, we have regularity on the region they cover.

\smallskip
\noindent\textbf{On Compact K\"ahler manifolds.}
Let $(X,\omega)$ be compact. Fix $V\Subset U\subset X$ and holomorphic-ball data satisfying \eqref{cond:localcx}, with uniform convexity constant $\lambda>0$ and chart margin $R_0>0$. The corresponding local barrier threshold is
\[
 \vartheta(V,U;\omega):=\lambda R_0^2.
\]
Here the threshold is determined by $\omega$ and the chosen local geometric data, independently of $\Psi$. If $\vartheta(V,U;\omega)\geq 16M$, by Theorem~\ref{mthmlocalreg}, there is a constant $\varepsilon_{V,U}>0$ such that
\[
 \|\Psi\|_{{k+3,\alpha;U\times\partial D}}\leq\varepsilon_{V,U}
\]
implies that $\Phi$ is smooth on $V\times D$, with the estimate \eqref{estilocal}. Thus the local criterion compares the geometric barrier threshold with the global boundary norm $M$, together with the required higher-norm smallness on $U\times\partial D$.

By local normalization of the K\"ahler potentials and compactness, choose a finite family
\[
 \mathcal U=\{V_i\Subset U_i\}_{i=1}^N,
 \qquad X=\bigcup_{i=1}^N V_i,
\]
with holomorphic coordinates $b_i:U_i\longrightarrow B_{2a_i}$, $b_i(V_i)=B_{a_i}$, and local potentials $\rho_i$ satisfying $\nabla_0^2\rho_i\geq\lambda_i I$ on $B_{2a_i}$, where $\lambda_i>0$. For every $x\in V_i$, the concentric balls
\[
 B_{a_i/2}(b_i(x))\Subset B_{a_i}(b_i(x))\subset B_{2a_i}
\]
satisfy the holomorphic-ball conditions of Theorem~\ref{mthmlocalreg}, with the uniform margin $R_{0,i}=a_i/2$. Define
\[
 \vartheta_i:=\lambda_i R_{0,i}^2,
 \qquad
 \vartheta_{\mathcal U}:=\min_{1\leq i\leq N}\vartheta_i
 =\min_{1\leq i\leq N}\frac{\lambda_i a_i^2}{4}>0.
\]
The number $\vartheta_{\mathcal U}$ depends only on $\omega$ and the chosen finite family of charts and potentials.

Suppose that $\vartheta_{\mathcal U}>16M$. For each pair $V_i\Subset U_i$, let $\varepsilon_i>0$ be the smallness constant supplied by Theorem~\ref{mthmlocalreg}, and set $\varepsilon_{\mathcal U}:=\min_{1\leq i\leq N}\varepsilon_i>0$. We obtain the following quantitative criterion:
\begin{equation}\label{eq:compact_global_criterion}
 16\|\Psi\|_{L^\infty(X\times\partial D)}<\vartheta_{\mathcal U},
 \qquad
 \max_{1\leq i\leq N}\|\Psi\|_{{k+3,\alpha;U_i\times\partial D}}
 \leq\varepsilon_{\mathcal U}.
\end{equation}
Under \eqref{eq:compact_global_criterion}, Theorem~\ref{mthmlocalreg} makes the global weak solution $\Phi$ smooth on every $V_i\times D$. Since the $V_i$ cover $X$, $\Phi$ is a smooth global solution of \eqref{introHCMA} on $X\times D$. By the local estimates on this finite cover, we also have
\[
 \|\Phi\|_{k,\beta;X\times D}
 \leq C_X\|\Psi\|_{k+1,\alpha;X\times\partial D},
 \qquad 0<\beta<\alpha,
\]
where the norms are measured using the fixed finite atlas. Hence, under the standing weak-solution assumptions, \eqref{eq:compact_global_criterion} is a quantitative sufficient criterion for the existence of a smooth global solution to the HCMA equation.

\smallskip
\noindent\textbf{On ALE K\"ahler manifolds.}
We first specify the asymptotic geometry used in this application. A complete noncompact K\"ahler manifold $(X,J,g)$ of complex dimension $n\geq2$ is \emph{asymptotically locally Euclidean (ALE) of order $\mu>0$} if there exist a compact set $K\subset X$, a radius $R_0>0$, a finite subgroup $\Gamma\subset U(n)$ acting freely on the unit sphere, and a diffeomorphism
\[
 e:X\setminus K\longrightarrow
 (\CC^n\setminus\overline{B}_{R_0})/\Gamma
\]
such that, in the lifted asymptotic coordinates,
\begin{equation}\label{eq:application_ale_definition}
 |\nabla_0^j(g-g_0)|_{g_0}
 +|\nabla_0^j(J-J_0)|_{g_0}
 \leq C_j r^{-\mu-j},
 \qquad j=0,1,\ldots.
\end{equation}
Here $g_0$ and $J_0$ are the Euclidean metric and standard complex structure, $\nabla_0$ is the Euclidean connection, and the same symbols $g,J$ denote the transported tensors. We write $\omega=g(J\cdot,\cdot)$ and let $r=|e(\cdot)|$ on the end, extended to a positive smooth function on $X$. For sufficiently large $L$, set
\[
 B_L:=\{x\in X:r(x)\leq L\},
 \qquad X_L:=X\setminus B_L=\{r>L\}.
\]
The asymptotic map $e$ need not be holomorphic. The argument below uses local holomorphic charts, so it applies also in complex dimension two.

Fix $\kappa=1/(40|\Gamma|)$. By the local holomorphic-chart construction of \cite[Lemma~2.10]{yao2025asymptotic}, there exists $L_0$ such that each $p\in X_{L_0}$ admits an injective holomorphic chart $b_p:\widetilde U_p\longrightarrow B_{4\kappa r(p)}$, with $b_p(p)=0$. Set $V_p:=b_p^{-1}(B_{\kappa r(p)})$ and $U_p:=b_p^{-1}(B_{2\kappa r(p)})$. In these coordinates,
\[
 1/2\omega_{\mathrm{euc}}
 \leq(b_p^{-1})^*\omega
 \leq2\omega_{\mathrm{euc}},
\]
where $\omega_{\mathrm{euc}}=\frac{i}{2}\sum_{j=1}^n dz^j\wedge d\bar z^j$. By \cite[Propositions~8.2 and~8.4(i)]{yao2025asymptotic}, after increasing $L_0$ if necessary, the form $(b_p^{-1})^*\omega$ admits a local potential $\rho_p$ on $B_{2\kappa r(p)}$ satisfying $\nabla_0^2\rho_p\geq\tfrac14 I$.

After recentring the coordinate balls at points of $V_p$, we have a uniform margin $\kappa r(p)/2$, so the local barrier threshold of $V_p\Subset U_p$ is $\kappa^2r(p)^2/16$. Thus
\begin{equation}\label{eq:ale_growing_barrier}
 \vartheta_{\mathrm{ALE}}(L):=
 \inf_{p\in X_L}\frac{\kappa^2r(p)^2}{16}
 \geq\frac{\kappa^2L^2}{16}
 \longrightarrow\infty.
\end{equation}

We now specify the boundary decay and derive the weighted estimates. Suppose that, for some $\gamma>0$,
\begin{equation}\label{eq:ale_application_boundary_decay}
 |\nabla_X^j\partial_\theta^\ell\Psi|_g
 =O(r^{-\gamma-j}),
 \qquad j,\ell\geq0,
\end{equation}
uniformly in the circle variable $\theta$ on $\partial D$, where $\nabla_X$ is the connection of $g$. For the weighted estimates, introduce the rescaling map $F_p:B_2\longrightarrow U_p$, $F_p(y):=b_p^{-1}(\kappa r(p)y)$. We use the scaled-chart description of the weighted H\"older norms. 
For $u$ defined on $U_p\times D$, set
\begin{equation}\label{eq:ale_application_weighted_norm}
 \|u\|_{-\gamma;k,\alpha;U_p\times D}
 :=r(p)^\gamma\|F_p^*u\|_{k,\alpha;B_2\times D},
\end{equation}
where $F_p^*u(y,t)=u(F_p(y),t)$, so only the spatial variable is rescaled. We use the same convention on $U_p\times\partial D$, and on $V_p$ with $B_2$ replaced by $B_1$. This is equivalent to the usual weighted norm with weight $r^{\gamma+j}$ on $j$ spatial derivatives; see \cite[Lemma~2.5]{yao2025asymptotic}. In particular, \eqref{eq:ale_application_boundary_decay} implies
\[
 M_k:=\sup_{p\in X_{L_0}}\|\Psi\|_{-\gamma;k,\alpha;U_p\times\partial D}<\infty
\]
for every fixed $k$.

Given $\Psi$ with global $L^\infty$ norm $M$, choose $L$ so that $\vartheta_{\mathrm{ALE}}(L)>16M$. Fix $k$ and $0<\beta<\alpha$. For each $p\in X_L$, rescale the global data $(\omega,\Phi,\Psi)$ by $(\kappa r(p))^{-2}$ and use $F_p$ to work on $B_1\Subset B_2$. The rescaled global bound is $M/(\kappa r(p))^2$, and the rescaled boundary data have {$C^{k+3,\alpha}$} norm $O(r(p)^{-\gamma-2})$. Both tend to zero, so, after increasing $L$, Theorem~\ref{mthmlocalreg} applies uniformly on these charts, so we have local higher regularity on $X_L\times D$. Rescaling back and taking the supremum over $p\in X_L$ yields
\begin{equation}\label{eq:ale_application_weighted_estimate}
 \sup_{p\in X_L}\|\Phi\|_{-\gamma;k,\beta;V_p\times D}
 \leq C M_{k+1}.
\end{equation}
See \cite[Section~9.2]{yao2025asymptotic} for details of the weighted estimates.
In particular, since each spatial derivative introduces a factor $(\kappa r(p))^{-1}$ while the disc variable is unchanged,
\[
 |\nabla_X^jD_\tau^\ell\Phi|_g
 \leq C_{j,\ell}r^{-\gamma-j},
 \qquad j+\ell\leq k,
\]
where $D_\tau^\ell$ denotes any real derivative of order $\ell$ in the disc variable. The spatial Hessian therefore decays as $O(r^{-\gamma-2})$ when $k\geq2$. The decay of $\Psi$ ensures the required smallness on $B_2$ after enlarging $L$ depending on {$M_{k+3}$} and $\gamma$; together with the blowup barrier threshold as $L \to +\infty$, this removes any global smallness requirement on $M$.

\smallskip
\noindent\textbf{On ALF ends.}
Let $(X,\omega)$ be a complete K\"ahler surface with an asymptotically locally flat (ALF) end. Outside a compact set $K$, the end is a circle fibration $\pi: X\backslash K \longrightarrow \RR^3\setminus B_{R_0}$ or its quotient by $\{\pm1\}$. Write $r=|\pi(\cdot)|$, the radius function inherited from the base, and $X_L=\{r>L\}$. Fix a reference model metric
\[
 g_0=\pi^*g_{\RR^3}+\eta^2,
\]
where $\eta$ is a chosen connection form, normalized so that the model circle fibres have a fixed positive length. In the quotient case, this description is understood on the corresponding double cover. For some $0<\mu<1$, assume the derivative decay estimates
\[
 |\nabla_{g_0}^j(g-g_0)|_{g_0}=O(r^{-\mu-j}),
 \qquad j=0,1,\ldots,
\]
\[
 |\nabla_X^j\operatorname{Rm}(g)|_g=O(r^{-2-\mu-j}),
 \qquad j=0,1,\ldots.
\]
Here $\nabla_{g_0}$ and $\nabla_X$ denote the Levi-Civita connections of $g_0$ and $g$, respectively. The fixed positive model fibre length and these bounds imply noncollapse at a fixed scale. Consequently, the injectivity radius is uniformly bounded below on the end, and all curvature derivatives are uniformly bounded.

For example, the single-centre Taub--NUT metric has the form
\[
 g=H\pi^*g_{\RR^3}+H^{-1}\eta^2,
 \qquad H=1+\frac{1}{2r},
 \qquad d\eta=\pi^*(*_{\RR^3}dH),
\]
where $*_{\RR^3}$ is the Euclidean Hodge star. By replacing the single pole by a finite sum of such poles at distinct points, we obtain the multi-Taub--NUT metrics \cite[Section~1.2]{minerbe2010asymptotic}. Differentiating this expression, we have the weighted model estimates
\[
 |\nabla_{g_0}^j(g-g_0)|_{g_0}=O(r^{-1-j}),
 \qquad j=0,1,\ldots.
\]
These examples also satisfy $|\nabla_X^j\operatorname{Rm}(g)|_g=O(r^{-3-j})$ for every $j\geq0$; see \cite[Appendix~A]{minerbe2010asymptotic}. They therefore satisfy the preceding assumptions for every $0<\mu<1$. The limiting circle length is nonzero, while the base is unbounded. We use exponential lifts to obtain growing barriers.

Let $D\subset\CC$ be the unit disc, and let $\Omega$ be the pullback of $\omega$ to $X\times D$. Suppose that $\Psi\in C^\infty(X\times\partial D,\RR)\cap L^\infty(X\times\partial D)$ satisfies $\omega+i\partial_X\bar\partial_X\Psi(\cdot,\tau)>0$ for every $\tau\in\partial D$. Consider the HCMA problem
\[
 \begin{cases}
  (\Omega+i\partial\bar\partial\Phi)^3=0
     &\text{in }X\times D,\\
  \Omega+i\partial\bar\partial\Phi\ge0
     &\text{in }X\times D,\\
  \Phi=\Psi
     &\text{on }X\times\partial D.
 \end{cases}
\]
Assume that $\Phi$ is a bounded continuous weak solution in the Bedford--Taylor sense and is given by the Perron envelope \eqref{cond:perron}. Set $M:=\|\Psi\|_{L^\infty(X\times\partial D)}$. The constant subsolution $-M$ and the maximum principle on vertical discs imply $-M\le\Phi\le M$.

Fix $0<\alpha<1$, and write $r_p=r(p)$. By \cite[Theorem~9(2)]{wuyau2020invariant}, choose injective holomorphic maps $\chi_p:B_a\subset\CC^2\to X$, with $\chi_p(0)=p$, where $a>0$ is independent of $p$ and the metrics $\chi_p^*g$ are uniformly comparable to $g_{\mathrm E}$ with uniform bounds on all derivatives. Set $W_p:=\chi_p(B_a)$. These are fixed-size coordinate neighborhoods, and we use the convention
\[
 \|\Psi\|_{C^{m,\alpha}(W_p\times\partial D)}
 :=\|\Psi\circ(\chi_p,\mathrm{id})\|_{C^{m,\alpha}(B_a\times\partial D)}.
\]
For some $\gamma>0$ and a fixed $L_0$, assume the boundary decay
\begin{equation}\label{eq:alf_boundary_decay}
 M_m:=\sup_{p\in X_{L_0}}r_p^\gamma
 \|\Psi\|_{C^{m,\alpha}(W_p\times\partial D)}<\infty,
 \qquad m=0,1,\ldots.
\end{equation}
The bounds $|\nabla_X^j\partial_\theta^\ell\Psi|_g=O(r^{-\gamma})$ for all $j,\ell$ imply this condition. We require only the common decay rate $r^{-\gamma}$, which suffices for the argument in Step~3. Faster decay, such as $O(r^{-\gamma-j})$ for $j$ spatial derivatives, is also allowed. No improvement per spatial derivative is required, since differentiation along the circle fibres, whose lengths remain bounded, need not improve decay. All Euclidean balls below are centred at zero unless otherwise specified, and pullbacks leave the disc variable unchanged.

\smallskip\noindent\textit{Step 1. Large holomorphic balls on the lift.}
Fix $p$ sufficiently far out and $1\le R\le r_p/128$. The circle bundle over the contractible base ball $B(\pi(p),r_p/4)$ is trivial. By unwrapping its circle fibres, we obtain a covering by $B(\pi(p),r_p/4)\times\RR$, equipped with the pulled-back K\"ahler structure. In the quotient case, first choose the lift of the base ball disjoint from its negative. We construct the required holomorphic balls on an exponential lift of this covering.

Identify the tangent space $T_pX$ with $\CC^2$ so that its metric and complex structure become $g_{\mathrm E}$ and $J_{\mathrm{std}}$. Let $B_1 \subset \CC^2$, we define
\[
 E_p:B_1\longrightarrow X,\qquad E_p(y)=\exp_p(8Ry),
\]
where $\exp_p(v)$ is the endpoint at time one of the geodesic starting at $p$ with initial velocity $v$. Each such geodesic has length at most $8R$. Since $g\ge g_0/2$, projection to the base increases lengths by at most a factor of two. Thus $|r(E_p(y))-r_p|\le16R\le r_p/8$, and the image stays over the chosen base ball. Since $B_1$ is simply connected, $E_p$ lifts to the covering above. Curvature decay implies local invertibility and the normal-coordinate estimate
\[
 \|(8R)^{-2}E_p^*g-g_{\mathrm E}\|_{C^{m,\alpha}(B_1)}
 +\|E_p^*J-J_{\mathrm{std}}\|_{C^{m,\alpha}(B_1)}
 \le C_mR^2r_p^{-2-\mu};
\]
see \cite[proof of Theorem~9]{wuyau2020invariant}. The constants are independent of $p,R$, and the right-hand side tends to zero uniformly over the admissible radii.

Apply \cite[Proposition~2.2]{heintosatti2020higher} to $E_p^*J$ on $B_1$; the preceding $C^{1,\alpha}$ estimate verifies its smallness hypothesis. Let $H_p:(B_{3/4},E_p^*J)\to\CC^2$ be the resulting holomorphic coordinate map, close to the identity in every prescribed finite norm. Subtract its value at zero so that $H_p(0)=0$. For distant $p$, $\|DH_p-I\|_{C^0}<1/4$, so $H_p$ is injective and its image contains $B_{1/2}$. The precise holomorphic ball chart on the exponential lift is
\[
 b_p:H_p^{-1}(B_{1/2})\overset{\sim}{\longrightarrow}B_{4R},
 \qquad b_p(y)=8R H_p(y).
\]
Its domain carries the complex structure $E_p^*J$ and the metric $E_p^*g$. Its inverse is $b_p^{-1}(z)=H_p^{-1}(z/(8R))$, defined for every $z\in B_{4R}$. Projecting this chart to $X$, we have
\[
 f_p:=E_p\circ b_p^{-1}:B_{4R}\longrightarrow X,
 \qquad f_p(0)=p.
\]
The coordinate metric satisfies $(b_p^{-1})^*E_p^*g=f_p^*g$, so the preceding estimates imply $\tfrac12g_{\mathrm E}\le f_p^*g\le2g_{\mathrm E}$ and the same radius bounds. The projection $f_p$ is locally biholomorphic; it need not be injective.

Finally, by the quantitative $\partial\bar\partial$-lemma on fixed balls, there are real potentials $\rho_p$ on $B_{2R}$ such that $R^{-2}\rho_p(R\,\cdot)$ converges to $|\cdot|^2/2$ in every fixed $C^{m,\alpha}(B_2)$ norm as $r_p\to\infty$, uniformly in the admissible $R$. Consequently,
\[
 f_p^*\omega=i\partial\bar\partial\rho_p,
 \qquad D^2\rho_p\ge\tfrac12I,
\]
with uniform bounds for the normalized potentials in every prescribed finite norm. These are the geometric bounds needed in Steps~2 and~3.

\smallskip\noindent\textit{Step 2. Growth of the local barrier threshold.}
For $L$ sufficiently large, every $p\in X_L$ admits the charts of Step~1 at any common radius $1\le R\le L/128$, since $r_p>L$. Consider the source pair $B_R\Subset B_{2R}$ with the pulled-back metric. For each $y\in B_R$, the concentric balls
\[
 B(y,R/4)\Subset B(y,3R/4)\Subset B_{2R}
\]
have margin $R/2$, and the translated potential has convexity constant $\lambda=1/2$. Consequently, an available local barrier threshold is
\[
 \vartheta(R)=\lambda(R/2)^2=\frac{R^2}{8}.
\]
Choosing $R=R(L):=L/128$, we therefore have
\begin{equation}\label{eq:alf_application_barrier}
 \vartheta(R(L))=\frac{R(L)^2}{8}
 =\frac{L^2}{131072}\longrightarrow\infty
 \qquad\text{as }L\longrightarrow\infty.
\end{equation}
This threshold is measured on the lifted source balls and is uniform for $p\in X_L$.

\smallskip\noindent\textit{Step 3. Local smallness and higher regularity.}
Fix $k\ge1$ and $0<\beta<\alpha$. Use the charts of Step~1 with radius $L/128$ and set $\psi_p:=\Psi\circ(f_p,\mathrm{id})$. For $p\in X_L$ and $L\ge2L_0$ sufficiently large, \eqref{eq:alf_boundary_decay} and the uniform coordinate bounds imply
\[
 \|\psi_p\|_{C^{k+3,\alpha}(B_{L/64}\times\partial D)}
 \le C_kM_{k+3}r_p^{-\gamma}
 \le C_kM_{k+3}L^{-\gamma}.
\]
By the uniform bounds for $D^2\rho_p$ and a fixed positive lower bound for the chart margins, we can choose a smallness constant in Theorem~\ref{thmexistbypertb} independent of $L$. Choose $L$ so that the preceding norm satisfies this smallness condition and $\vartheta(L/128)\ge16(M+1)>16M$.

Apply {the disc and support construction in} the proof of Theorem~\ref{mthmlocalreg} on {$B_R\Subset B_{2R}$, with $R=L/128$. The local subsolution, truncated below by $-M$, equals $-M$ near $\partial B_{2R}\times D$ and is at most $\psi_p$ on $B_{2R}\times\partial D$. Since $\rho_p+\Phi\circ(f_p,\mathrm{id})$ is bounded plurisubharmonic with vanishing Monge--Amp\`ere measure, Bedford--Taylor maximality and this boundary ordering place} the {local} {subsolution} {below the pullback of $\Phi$; comparison with the harmonic supports along the discs} implies {the reverse inequality there. For sufficiently small displacement, the projected leaves cover $B_{R/2}\times D$, so} $\Phi\in C^{k,\beta}_{\mathrm{loc}}(X_L\times D)$. Local estimates on fixed-size source balls, with affine terms removed from $\rho_p$, transfer to $W_p$ and imply
\begin{equation}\label{eq:alf_application_weighted_estimate}
 \sup_{p\in X_L}r_p^\gamma
 \|\Phi\|_{C^{k,\beta}(W_p\times D)}
 \le C_{M,k}M_{k+1}.
\end{equation}
The constant depends only on $M,k,\alpha,\beta,\gamma$ and the fixed geometry; the cutoff $L$ may also depend on $M_{k+3}$.

\smallskip
\noindent\textbf{On asymptotically complex hyperbolic (ACH) K\"ahler metrics.}
Let $X\Subset\CC^n$, $n\ge2$, be smoothly bounded and strictly pseudoconvex. Fix a smooth boundary defining function $t\in C^\infty(\overline X)$ with the sign convention
 $t>0$ in $X$,
$ t=0$ on $\partial X$,
$ dt\ne0$ on $\partial X$,
and $-t$ is strictly plurisubharmonic on a neighborhood of $\partial X$; such a choice is possible by strict pseudoconvexity. 
In particular, $t$ is comparable to the Euclidean distance to $\partial X$ on a sufficiently small boundary collar. On this collar, set the model metrics as
\[
 \omega_0=i\partial\bar\partial(-\log t),\qquad
 g_0=\omega_0(\cdot,J\cdot).
\]
Then, $\omega_0$ is a K\"ahler form on the collar. Let $X_L = \{0<t<L^{-1}\}$. $X_L$ is a shrinking neighborhood of the boundary as $L\to\infty$. For example, on the unit ball one may take $t(z)=1-|z|^2$, which corresponds to the complex hyperbolic model $\omega_0=-i\partial\bar\partial\log(1-|z|^2)$.

Let $\omega$ be a complete K\"ahler form on $X$. Relative to the chosen defining function $t$ and its reference metric $g_0$, we impose the following asymptotic condition: for some $\mu>0$ and $0<\alpha<1$, on a fixed boundary collar and for every $m\ge0$,
\begin{equation}\label{eq:ach_metric_decay}
 \|t^{-\mu}(\omega-\omega_0)\|_{C^{m,\alpha}(g_0)}<\infty.
\end{equation}
This includes the Cheng--Yau metrics normalized by $\operatorname{Ric}(\omega)=-(n+1)\omega$, with a suitable choice of $t$; see \cite{cheng1980existence,leemelrose1982boundary} and \cite[Section~2]{matsumoto2018einstein}. We consider \eqref{introHCMA} with $\Omega=\pi^*\omega$ and the global hypotheses above. For some $\gamma>0$, prescribe boundary data $\Psi$ with $\omega+i\partial_X\bar\partial_X\Psi(\cdot,\tau)>0$ and
\begin{equation}\label{eq:cheng_yau_boundary_decay}
 M_m:=\|t^{-\gamma}\Psi\|_{C^{m,\alpha}(X_{L_0}\times\partial D,\,g_0+|d\tau|^2)}<\infty,
 \qquad m=0,1,\ldots.
\end{equation}
The global bound $M=\|\Psi\|_{L^\infty(X\times\partial D)}$ need not be small.

\smallskip\noindent\textit{Step 1. The model threshold.}
The local barrier threshold $\lambda R_0^2$ introduced in Section~1 is based on a uniform Euclidean convexity bound and is not sharp for the complex hyperbolic model. We therefore compute a separate threshold from the boundary heights of explicit model barriers, which captures their growth with the hyperbolic radius.

On the standard ball $B^n=\{w\in\CC^n:|w|<1\}$, set
\[
 \rho_{\rm CH}(w)=-\log(1-|w|^2),\qquad
 \omega_{\rm CH}=i\partial\bar\partial\rho_{\rm CH}.
\]
For $g_{\rm CH}=\omega_{\rm CH}(\cdot,J\cdot)$, the radial distance is
\[
 r(w)=d_{\rm CH}(0,w)
 =\sqrt2\int_0^{|w|}\frac{ds}{1-s^2}
 =\sqrt2\,\operatorname{arctanh}|w|, \qquad w\in B^n.
\]
Hence the potential height at hyperbolic distance $r$ from the origin is
\[
 \rho_{\rm CH}(w)=2\log\cosh(r/\sqrt2).
\]
We use nonnegative potentials, whose negatives are the lower barriers.

Write $B_{\rm CH}(r)=B_{g_{\rm CH}}(0,r)=B_{\tanh(r/\sqrt2)}$. To center the potential at $p\in B^n$, choose a holomorphic isometry $A_p$ taking $p$ to $0$ and set
\[
 h_p(y):=\rho_{\rm CH}(A_p(y))
 =2\log\cosh\frac{d_{\rm CH}(p,y)}{\sqrt2}.
\]
Then $h_p\ge0$, $h_p(p)=0$, and $i\partial\bar\partial h_p=\omega_{\rm CH}$. Let $V\Subset U\Subset B^n$ be open sets satisfying
$ \operatorname{dist}_{\omega_{\rm CH}}(V,\partial U)\ge L>0.$
For every $p\in\overline V$ and $y\in\partial U$, we have
\[
 h_p(y)
 =2\log\cosh\frac{d_{\rm CH}(p,y)}{\sqrt2}
 \ge2\log\cosh\frac{L}{\sqrt2}
 \ge\sqrt2\,L-2\log2,
\]
where the last inequality follows from $\log\cosh s\ge s-\log2$ for $s\ge0$. Consequently, an available local barrier threshold for any such pair $V\Subset U$ is
\begin{equation}\label{eq:cheng_yau_application_barrier}
 \vartheta_{\rm CH}(L):=2\log\cosh\frac{L}{\sqrt2}.
\end{equation}
This threshold depends only on $L$ and satisfies $\vartheta_{\rm CH}(L)\ge cL-C$, with $c=\sqrt2$ and $C=2\log2$.

\smallskip\noindent\textit{Step 2. Barrier thresholds in the general ACH domain.}
We now return to the general ACH data $(X,t,\omega)$ satisfying \eqref{eq:ach_metric_decay}. Recall that $X_L=\{0<t<L^{-1}\}$, and let $B_\omega$ denote balls for the distance $d_\omega$ of $g_\omega$. Fix $L\ge\max\{L_0,4\}$. We choose a deeper cutoff $L'=L'(L)$, with
$L'\ge L e^{2(L+2)}$,
and, for every $p\in X_{L'}$, consider the concentric intrinsic balls
\begin{equation}\label{eq:ach_nested_balls}
 V_p:=B_\omega(p,1)\Subset U_p:=B_\omega(p,L+1)\Subset X_L.
\end{equation}
The displayed bound on $L'$ ensures collar containment; we increase $L'$ further to obtain the following uniform potentials.

\begin{lem}\label{lem:ach_distance_barrier}
For each fixed $L\ge\max\{L_0,4\}$, there exists $L'=L'(L)$ such that, for every $p\in X_{L'}$, a neighborhood of $\overline U_p$ admits holomorphic coordinates and smooth convex potentials $h_{p,y}$, $y\in\overline V_p$, satisfying $i\partial\bar\partial h_{p,y}=\omega$, $h_{p,y}(y)=0$, and $h_{p,y}\ge0$. Moreover,
\begin{equation}\label{eq:ach_distance_barrier}
 \begin{aligned}
 h_{p,y}(z)&\ge d_\omega(y,z)-10,
 &&z\in\overline U_p,
 \end{aligned}
\end{equation}
\end{lem}

The lemma is proved in Appendix~\ref{app:ach_barriers}. Thus the constants in the distance bound can be taken as $c=1$ and $C=10$, independently of $p$ and $L$; the required cutoff $L'$ depends on $L$, the fixed defining function, and the geometric bounds in \eqref{eq:ach_metric_decay}. Convexity here refers to the stated holomorphic coordinates.

For $y\in\overline V_p$ and $z\in\partial U_p$, the triangle inequality implies $d_\omega(y,z)\ge L$, hence $h_{p,y}(z)\ge L-10$. Consequently, for $L>10$, an available local barrier threshold for every pair in \eqref{eq:ach_nested_balls} is
\begin{equation}\label{eq:ach_intrinsic_boundary_height}
 \vartheta(L):=L-10.
\end{equation}
This threshold is uniform in $p$ and tends to infinity as $L\to\infty$.
Both the outer sets and their centers approach the boundary as $L\to\infty$.

\smallskip\noindent\textit{Step 3. Higher regularity on the end.}
Recall the ACH data $(X,t,\omega)$ above: $X\Subset\CC^n$ is smoothly bounded and strictly pseudoconvex, $t$ is the chosen defining function, and $\omega$ is a complete K\"ahler form satisfying \eqref{eq:ach_metric_decay}. Let $D\subset\CC$ be the unit disc and $\Omega=\pi^*\omega$, where $\pi:X\times D\to X$ is the projection. Suppose that $\Psi\in C^\infty(X\times\partial D,\RR)\cap L^\infty(X\times\partial D)$ satisfies $\omega+i\partial_X\bar\partial_X\Psi(\cdot,\tau)>0$ for every $\tau\in\partial D$ and the weighted decay condition \eqref{eq:cheng_yau_boundary_decay} for some $\gamma>0$, with $M_m<\infty$ for all $m\ge0$.

Consider the HCMA problem \eqref{introHCMA},
\[
 \begin{cases}
  (\Omega+i\partial\bar\partial\Phi)^{n+1}=0
     &\text{in }X\times D,\\
  \Omega+i\partial\bar\partial\Phi\ge0
     &\text{in }X\times D,\\
  \Phi=\Psi
     &\text{on }X\times\partial D.
 \end{cases}
\]
Assume that $\Phi$ is a bounded continuous weak solution, with the equation understood in the Bedford--Taylor sense, and is given by the Perron envelope \eqref{cond:perron}. Set $M:=\|\Psi\|_{L^\infty(X\times\partial D)}$. The constant subsolution $-M$ and the maximum principle on each vertical disc imply $-M\le\Phi\le M$, as required in Theorem~\ref{mthmlocalreg}.

Fix $k\ge1$ and $0<\beta<\alpha$. Choose $L\ge\max\{L_0,4\}$ so that $\vartheta(L)=L-10>2M+2$, and keep $L$ fixed. By the construction in Appendix~\ref{app:ach_barriers}, there are holomorphic charts $f_p:B_{\rm CH}(S)\to X$, with a fixed $S>\sqrt2(L+2)$ and images containing $\overline U_p$, whose local potentials are uniformly strictly convex and have uniform bounds in every prescribed finite norm. Writing $\psi_p:=\Psi\circ(f_p,\mathrm{id})$, the boundary decay implies
\[
 \|\psi_p\|_{C^{k+3,\alpha}(B_{\rm CH}(S)\times\partial D)}
 \le C_{L,k}M_{k+3}t(p)^\gamma.
\]
Increasing $L'$ makes this norm uniformly small for $p\in X_{L'}$. In particular, both the local smallness condition of Theorem~\ref{mthmlocalreg} and the condition on $\partial\psi_p$ in the disc construction are satisfied.

The proof of Theorem~\ref{mthmlocalreg} now applies with the boundary-height bound \eqref{eq:ach_intrinsic_boundary_height} in place of its quadratic barrier estimate. Appendix~\ref{app:ach_barriers} verifies {the required outer-boundary} bound for the {perturbed} affine supports {and} carries out the global subsolution construction and Perron comparison. We obtain $C^{k,\beta}$ regularity on each $V_p\times D$. Since these sets cover $X_{L'}\times D$, the local estimates yield
\begin{equation}\label{eq:cheng_yau_application_weighted_estimate}
 \|t^{-\gamma}\Phi\|_{C^{k,\beta}(X_{L'}\times D,\,g_0+|d\tau|^2)}
 \le C_{M,k}M_{k+1}.
\end{equation}
There is no smallness restriction on $M$. Besides the fixed geometric parameters, this final cutoff may depend on $k$, $M$, $M_{k+3}$, and $\gamma$; the constants are uniform on the resulting collar.

\appendix

\section{Campanato estimates}
\label{app:campanato}

In this appendix we recall a simple form of the Campanato characterization
of H\"older spaces. This result is a special case of the standard
Morrey--Campanato characterization; see, for instance,
\cite[Theorem 5.5]{GiaquintaMartinazzi2012}. Since the dependence of the
constant on the H\"older exponent is important for the BMO--H\"older
estimate in Lemma \ref{lemBMOHolderInterpolation}, we include the short
proof and keep track of this dependence explicitly.

\begin{lem}
\label{lemCampanatoExplicit}
Let \(0<s<1\), and let \(f\in L^1(\partial D)\). 
Suppose that \(f\) has a finite \(s\)-Campanato seminorm:
\[
[f]_{\mathcal C^s}
:=
\sup_{I\subset \partial D}
|I|^{-s}\fint_I |f-f_I|< + \infty.
\]
Then \(f\) admits an $s$-H\"older continuous representative.
Furthermore, we have the equivalent seminorms
\[
[f]_{\mathcal C^s}\le [f]_{0,s;\partial D}
\le
\frac{C}{s}
[f]_{\mathcal C^s}.
\]
Here \(C\) depends only on \(\partial D\), and in particular is independent
of \(s\) and \(f\).
\end{lem}

\begin{proof}
For \(x\in \partial D\) and \(r>0\), let \(I_r(x)\) denote the arc centered
at \(x\) with length \(2r~ (\le 2\pi)\), and set
\[
f_{x,r}:=\fint_{I_r(x)}f.
\]
Since \(I_{r/2}(x)\subset I_r(x)\), we have
\[
\begin{aligned}
|f_{x,r/2}-f_{x,r}|
&=
\left|
\fint_{I_{r/2}(x)}(f-f_{x,r})
\right|  \\
&\le
2\fint_{I_r(x)}|f-f_{x,r}| \le
2[f]_{\C^s} (2r)^s.
\end{aligned}
\]
Using the telescope technique again, we have
\[
|f_{x,r}-f_{x,2^{-N}r}|
\le
2[f]_{\C^s}(2r)^s\sum_{j=0}^{N-1}2^{-js}.
\]
Since
\( \displaystyle
\sum_{j=0}^{\infty}2^{-js}
=
\frac{1}{1-2^{-s}}
\le
\frac{C}{s},
\)
we obtain
\[
|f_{x,r}-f_{x,2^{-N}r}|
\le
\frac{C}{s}r^s[f]_{\C^s}.
\]
By the Lebesgue differentiation theorem, for a.e. \(x\in \partial D\),
\(
\lim_{r\to 0} f_{x,r}=f(x).
\)
For $R/2<r\le R$, the same averaging estimate implies $|f_{x,r}-f_{x,R}|\le 2[f]_{\mathcal C^s}(2R)^s$. Comparing arbitrary radii with adjacent dyadic radii, the dyadic estimate above implies that \(f_{x,r}\) has a limit as
\(r\to 0\) for every \(x\). We denote this limit by
\(
\bar f(x):=\lim_{r\to 0}f_{x,r}.
\)
Then, \(\bar f=f\) almost everywhere. Replacing \(f\) by \(\bar f\), we may
assume that
\(
f(x)=\lim_{r\to 0}f_{x,r}
\)
for every \(x\in\partial D\).
Hence, we have
\[
|f(x)-f_{x,r}|
\le
\frac{C}{s}r^s[f]_{\C^s}.
\]

Now take \(x,y\in \partial D\), and set
\(
r:=d_{\partial D}(x,y)/2.
\)
Let \(I_x=I_r(x)\), \(I_y=I_r(y)\), and let \(J\) be an arc containing both
\(I_x\) and \(I_y\), with \(|J|\le 5r\). Then
\[
|f(x)-f(y)|
\le
|f(x)-f_{I_x}| + |f_{I_x}-f_J| + |f_J-f_{I_y}| + |f_{I_y}-f(y)|.
\]
The first and last terms are bounded by
\(
\frac{C}{s}r^s [f]_{\C^s}.
\)
Moreover,
\[
|f_{I_x}-f_J|
\le
\frac{|J|}{|I_x|}\fint_J |f-f_J|
\le
\frac52[f]_{\C^s}|J|^s
\le
C[f]_{\C^s}r^s,
\]
and similarly
\(
|f_{I_y}-f_J|
\le C[f]_{\C^s}r^s.
\)
Therefore
\(
|f(x)-f(y)|
\le
\frac{C}{s}[f]_{\C^s}r^s.
\)
Hence
\[
[f]_{0,s;\partial D}
\le
\frac{C}{s} [f]_{\C^s}.
\]

Conversely, if $f$ is $s$-H\"older continuous, then for any arc $I \subseteq \partial D$, we have
\[
|I|^{-s} \fint_I |f-f_I| \le |I|^{-s}\fint_I\fint_I |f(x)-f(y)|\,dx\,dy \le [f]_{0,s; \partial D}.
\]
\end{proof}

\section{Convex barriers and distance on ACH ends}
\label{app:ach_barriers}

We prove Lemma~\ref{lem:ach_distance_barrier}. Recall that
 $V_p=B_\omega(p,1)\Subset U_p=B_\omega(p,L+1)\Subset X_L$
for $p\in X_{L'}$ with $L'$ sufficiently large. We construct convex local potentials $h_{p,y}$ of $\omega$ satisfying 
\[
h_{p,y}(z)\ge c\,d_\omega(y,z)-C
\]
for $y\in\overline V_p$ and $z\in\overline U_p$, with $c=1$ and $C=10$. They are obtained by subtracting affine tangent planes from a potential close to $\rho_{\rm CH}$ in suitable holomorphic coordinates.
\begin{proof}
Increase the initial cutoff $L_0$, if necessary, so that $\tfrac12\omega_0\le\omega\le2\omega_0$ on $X_{L_0}$, and fix $L\ge\max\{L_0,4\}$. Since
\[
 \omega_0=-\frac{i\partial\bar\partial t}{t}
          +i\partial\log t\wedge\bar\partial\log t
\]
and $-t$ is strictly plurisubharmonic, $|d\log t|_{g_0}\le\sqrt2$ and hence $|d\log t|_{g_\omega}\le2$. To ensure $\overline U_p\subset X_L$ for every $p\in X_{L'}$, we want $\operatorname{dist}_{g_\omega}(X_{L'},\partial X_L)>L+1$. The gradient bound implies $\operatorname{dist}_{g_\omega}(X_{L'},\partial X_L)\ge\tfrac12\log(L'/L)$, so it suffices to choose $L'\ge Le^{2(L+2)}$.

\smallskip\noindent\textit{Step 1. Boundary coordinates.}
Extend $t$ smoothly to a neighborhood of $\partial X$. For $q\in\partial X$, let $\nu_q = \nabla t(q)/|\nabla t(q)|$ be the inward Euclidean unit normal. By strict pseudoconvexity, we can choose holomorphic coordinates $z=z_q(x)$ centered at $q$ such that
\begin{equation}\label{eq:ach_boundary_normal_form}
 t_q(z):=t(z_q^{-1}(z))
 =2\operatorname{Re}z_n-|z'|^2+O(|z_n|^2+|z_n||z'|+|z'|^3).
\end{equation}
Denote the remainder in \eqref{eq:ach_boundary_normal_form} by $E_q= O(|z_n|^2+|z_n||z'|+|z'|^3)$.
These coordinates also satisfy, for small $s$,
$z_q(q+s\nu_q)=(0,s|\nabla t(q)|/2)$.
The charts have a common coordinate radius, and the coordinate maps, their inverses, and $E_q$ have uniform $C^k$ bounds for each fixed $k$, independently of $q$.

\smallskip\noindent\textit{Step 2. Rescaling and estimating derivatives.}
Write each $p$ sufficiently close to $\partial X$ uniquely as $p=q+s\nu_q$, put $\delta=s|\nabla t(q)|$, and define
\[
 F_\delta(w)=(F_\delta'(w),F_{\delta,n}(w))=\left(
 \frac{\sqrt\delta\,w'}{1+w_n},
 \frac{\delta(1-w_n)}{2(1+w_n)}\right),\qquad
 f_\delta=z_q^{-1}\circ F_\delta.
\]
By Step 1,
\[
 F_\delta(0)=(0,\delta/2)=z_q(p),\qquad
 f_\delta(0)=p,\qquad t(p)=\delta+O(\delta^2).
\]
The map $F_\delta$ is injective and, for small $\delta$, maps $B_{\rm CH}(2L+1)$ into the coordinate neighborhood. By direct computation, we have
\[
 2\operatorname{Re}F_{\delta,n}(w)-|F_\delta'(w)|^2
 =\delta\frac{1-|w_n|^2-|w'|^2}{|1+w_n|^2} 
 =\delta \frac{1-|w|^2}{|1+ w_n|^2}
 =:\delta Q(w),
\]
and hence
\[
 t\circ f_\delta=t_q\circ F_\delta=\delta Q+e_\delta,
\]
where $e_\delta:=E_q\circ F_\delta$ is the additive error.

We estimate $e_\delta$ in $C^{m,\alpha}$ for each fixed $m$. Write $\|\cdot\|_{m,\alpha}$ for the $C^{m,\alpha}(B_{\rm CH}(2L+1))$ norm. Since $|1+w_n|$ is bounded away from zero on this fixed ball,
\[
 \|F_\delta'\|_{m,\alpha}\le C_{L,m}\sqrt\delta,
 \qquad \|F_{\delta,n}\|_{m,\alpha}\le C_{L,m}\delta,
 \qquad \|Q^{-1}\|_{m,\alpha}\le C_{L,m}.
\]
Taylor's formula in real coordinates writes $E_q$ as quadratic monomials containing a normal factor, plus cubic monomials with smooth coefficients uniformly bounded in every fixed $C^k$ norm. The product and composition estimates therefore imply, for $0<\delta\le1$,
\[
 \begin{aligned}
 \|e_\delta\|_{m,\alpha}
 &\le C_{L,m}\bigl(\|F_{\delta,n}\|_{m,\alpha}^2
       +\|F_{\delta,n}\|_{m,\alpha}\|F_\delta'\|_{m,\alpha}
       +\|F_\delta\|_{m,\alpha}^3\bigr)\\
 &\le C_{L,m}(\delta^2+\delta^{3/2})
 \le C_{L,m}\delta^{3/2}.
 \end{aligned}
\]
The relative error $\epsilon_\delta:=e_\delta/(\delta Q)$ therefore satisfies
\begin{equation}\label{eq:ach_scaled_defining_function}
 \begin{aligned}
 t\circ f_\delta&=\delta Q(1+\epsilon_\delta),\\
 \|\epsilon_\delta\|_{m,\alpha}
 &\le C_{L,m}\delta^{-1}\|e_\delta\|_{m,\alpha}
 \le C_{L,m}\sqrt\delta.
 \end{aligned}
\end{equation}
On the fixed ball, $Q$ is bounded above and away from zero. For small $\delta$, $0<t\circ f_\delta\le C_L\delta<L^{-1}$, so $f_\delta(B_{\rm CH}(2L+1))\subset X_L$. Since $(t\circ f_\delta)/t(p)=(\delta/t(p))Q(1+\epsilon_\delta)$ and $\delta/t(p)=1+O(\delta)$, the normalized defining function and its reciprocal satisfy
\begin{equation}\label{eq:ach_normalized_defining_function}
 \left\|\frac{t\circ f_\delta}{t(p)}\right\|_{m,\alpha}
 +\left\|\frac{t(p)}{t\circ f_\delta}\right\|_{m,\alpha}
 \le C_{L,m}.
\end{equation}
All constants are independent of $q$, with $L$ fixed before $\delta\to0$. These bounds imply $\|(t\circ f_\delta)^\mu\|_{m,\alpha}\le C_{L,m}\delta^\mu$, which Step 3 uses to estimate $f_\delta^*(\omega-\omega_0)$.

\smallskip\noindent\textit{Step 3. Constructing a convex potential.}
Taking logarithms in \eqref{eq:ach_scaled_defining_function}, we have
\[
 -\log(t\circ f_\delta)
 =-\log\delta+\rho_{\rm CH}+\log|1+w_n|^2+\varepsilon_0,
 \qquad \varepsilon_0=-\log(1+\epsilon_\delta).
\]
The constant and $\log|1+w_n|^2$ are pluriharmonic, and $\|\varepsilon_0\|_{C^{m+2,\alpha}}\le C_{L,m}\sqrt\delta$. Thus $f_\delta^*\omega_0=i\partial\bar\partial(\rho_{\rm CH}+\varepsilon_0)$. It follows that $f_\delta^*g_0$, its inverse, and their derivatives have uniform bounds on the fixed ball. Denote the pulled-back difference of the K\"ahler forms by
\[
 \eta_\delta:=f_\delta^*(\omega-\omega_0)
 =(t\circ f_\delta)^\mu f_\delta^*\bigl(t^{-\mu}(\omega-\omega_0)\bigr).
\]
These metric bounds turn the covariant derivative estimates in \eqref{eq:ach_metric_decay} into ordinary derivative estimates for $f_\delta^*\bigl(t^{-\mu}(\omega-\omega_0)\bigr)$. Using \eqref{eq:ach_normalized_defining_function} and the product estimate, we obtain
\[
 \|\eta_\delta\|_{C^{m,\alpha}(B_{\rm CH}(2L+1))}
 \le C_{L,m}\delta^\mu.
\]

We need a potential for $\eta_\delta$ with a small full real Hessian. Since $\eta_\delta$ is a closed real $(1,1)$-form, the real radial homotopy and the weighted $L^2$ estimate for $\bar\partial$ with weight $|w|^2$ (see \cite{hormander1973introduction}) imply the existence of a real-valued potential $\varepsilon_1$ on $B_{\rm CH}(2L+1)$ satisfying
\[
 i\partial\bar\partial\varepsilon_1=\eta_\delta,\qquad
 \|\varepsilon_1\|_{L^2(B_{\rm CH}(2L+1))}\le C_L\delta^\mu.
\]
Write $\eta_\delta=i\sum\eta_{j\bar k}\,dw_j\wedge d\bar w_k$. Taking the trace, we have $\Delta\varepsilon_1=4\sum_j\eta_{j\bar j}$. Interior elliptic regularity shows that $\varepsilon_1$ is smooth. We then apply the interior Schauder estimate with an $L^2$ term to obtain
\[
 \|\varepsilon_1\|_{C^{m+2,\alpha}(B_{\rm CH}(2L))}
 \le C_{L,m}\bigl(\|\varepsilon_1\|_{L^2(B_{\rm CH}(2L+1))}
       +\|\eta_\delta\|_{C^{m,\alpha}(B_{\rm CH}(2L+1))}\bigr)
 \le C_{L,m}\delta^\mu.
\]
Combining the two potential corrections, we have
\begin{equation}\label{eq:ach_controlled_potential}
 f_\delta^*\omega=i\partial\bar\partial(\rho_{\rm CH}+e),\qquad
 \|e\|_{C^{m+2,\alpha}(B_{\rm CH}(2L))}
 \le C_{L,m}(\sqrt\delta+\delta^\mu).
\end{equation}
Here $e=\varepsilon_0+\varepsilon_1$, and the constants are uniform in $q$. The norm bound follows from the triangle inequality and the preceding estimates for $\varepsilon_0$ and $\varepsilon_1$. Choose $\delta$ small enough that
\begin{equation}\label{eq:ach_chart_comparison}
 \tfrac12g_{\rm CH}\le f_\delta^*g_\omega\le2g_{\rm CH},\qquad
 \|e\|_\infty,\|de\|_\infty\le\tfrac14,\qquad \|D^2e\|_{\rm op}\le1.
\end{equation}
For a real vector $\eta$, by direct differentiation we have
\[
 D^2\rho_{\rm CH}(w)[\eta,\eta]
 =\frac{2|\eta|^2}{1-|w|^2}
  +\frac{4(\operatorname{Re}\langle w,\eta\rangle)^2}{(1-|w|^2)^2}
 \ge2|\eta|^2.
\]
By \eqref{eq:ach_chart_comparison}, we obtain $D^2(\rho_{\rm CH}+e)\ge I$. This uses the full $C^2$ bound for the error $e$. A positive complex Hessian alone would not imply real convexity.

\smallskip\noindent\textit{Step 4. Constructing the barriers and proving the distance estimate.}
Write $f=f_\delta$. The metric comparison implies
\[
 \overline U_p\subset f(B_{\rm CH}(2L)),\qquad
 \overline V_p\subset f(B_{\rm CH}(3/2)),\qquad
 f(B_{\rm CH}(2))\subset U_p.
\]
For fixed $L$, the error bounds in \eqref{eq:ach_scaled_defining_function} and \eqref{eq:ach_controlled_potential} have constants independent of $q$. Thus the chart construction and the metric and $C^2$ bounds in \eqref{eq:ach_chart_comparison} hold for every $q$ once $\delta$ is smaller than a single constant depending on $L$ and the fixed geometric data. Since $t(p)=\delta+O(\delta^2)$ uniformly and $t(p)<1/L'$ for $p\in X_{L'}$, increasing $L'$ makes this choice uniform for all $p\in X_{L'}$. The Taylor argument below takes place in $B_{\rm CH}(2L)$, which contains both the centers $v\in\overline{B_{\rm CH}(2)}$ and $f^{-1}(U_p)$.

For $v\in\overline{B_{\rm CH}(2)}$, define on the chart
\[
 h_{p,f(v)}(f(w))=(\rho_{\rm CH}+e)(w)-(\rho_{\rm CH}+e)(v)-d(\rho_{\rm CH}+e)|_v(w-v).
\]
The line segment from $v$ to $w\in B_{\rm CH}(2L)$ stays inside this Euclidean ball. Taylor's formula and $D^2(\rho_{\rm CH}+e)\ge I$ therefore imply
\[
 h_{p,f(v)}(f(w))
 =\int_0^1(1-s)D^2(\rho_{\rm CH}+e)(v+s(w-v))[w-v,w-v]\,ds
 \ge\tfrac12|w-v|^2.
\]
Thus $h_{p,f(v)}\ge0$ and $h_{p,f(v)}(f(v))=0$. A real affine function is pluriharmonic, so $i\partial\bar\partial h_{p,f(v)}=\omega$, and $D_w^2(h_{p,f(v)}\circ f)=D^2(\rho_{\rm CH}+e)\ge I$. This proves convexity in the $w$-coordinates. It does not assert real convexity in the original ambient coordinates.

For the distance estimate, first note that for $|w|<1$,
\[
 d\rho_{\rm CH}|_v(w-v)
 =\frac{2\operatorname{Re}\langle w-v,v\rangle}{1-|v|^2}
 \le\frac{2|v|}{1+|v|}\le1.
\]
Also $\rho_{\rm CH}(v)\le2\sqrt2$. Subtracting the affine tangent plane from $e$ leaves a remainder whose absolute value is at most $2\|e\|_\infty+2\|de\|_\infty\le1$. Thus
\[
 h_{p,f(v)}(f(w))\ge\rho_{\rm CH}(w)-2\sqrt2-2.
\]
Put $r=d_{\rm CH}(0,w)$. Comparing lengths along radial paths and using the model formula, we have
\[
 d_\omega(f(v),f(w))\le2\sqrt2+\sqrt2r,\qquad
 \rho_{\rm CH}(w)=2\log\cosh(r/\sqrt2)
 \ge\sqrt2r-2\log2.
\]
It follows that
\[
 h_{p,f(v)}(f(w))
 \ge d_\omega(f(v),f(w))-(4\sqrt2+2+2\log2)
 \ge d_\omega(f(v),f(w))-10.
\]
If $f(w)\in\partial U_p$, the same comparison along a radial path implies $L+1=d_\omega(p,f(w))\le\sqrt2r$. By the preceding calculation, we then have
\begin{equation}\label{eq:ach_extended_boundary_height}
 h_{p,f(v)}(z)\ge L-10,
 \qquad v\in\overline{B_{\rm CH}(2)},\quad z\in\partial U_p.
\end{equation}
\end{proof}

{For the application in Step~3 of the ACH discussion, fix $L$ with $L-10>2M+2$. Construct the disc family with labels in an intermediate parameter ball between $B_{\rm CH}(3/2)$ and $B_{\rm CH}(2)$. For sufficiently small displacement, the projected leaves cover $f^{-1}(V_p)\times D$ and remain in $B_{\rm CH}(2)$. Increasing $L'$ makes
\[
 D_w^2(\rho_{\rm CH}+e+\psi_p)\ge\tfrac12 I,
 \qquad
 \|\psi_p\|_\infty+2\|d_w\psi_p\|_\infty\le1
 \quad\text{on }B_{\rm CH}(2L)\times\partial D.
\]
For $\tau\in\partial D$, put $v=z(x,\tau)$. The affine support from Section~\ref{secconstructF} satisfies
\[
 L_x(w,\tau)-(\rho_{\rm CH}+e)(w)
 =-h_{p,f(v)}(f(w))+\psi_p(v,\tau)
   +d_w\psi_p(v,\tau)(w-v).
\]
Since $|w-v|\le2$, equation~\eqref{eq:ach_extended_boundary_height} implies
\[
 L_x(w,\tau)-(\rho_{\rm CH}+e)(w)
 \le -(L-10)+1<-2M-1,
 \qquad f(w)\in\partial U_p.
\]
The uniform coordinate bounds imply the weaker bound $-2M$ in a spatial collar of $\partial U_p$, and harmonicity in $\tau$ extends it to every $\tau\in D$. Convexity also implies $L_x-(\rho_{\rm CH}+e)\le\psi_p$ on the disc boundary. Thus the support-envelope construction of Section~\ref{secconstructF}, truncated below by $-M$, extends by $-M$ outside $U_p$ to a bounded global $\Omega$-psh subsolution. It lies below $\Phi$ by the Perron property. Comparing $\Phi$ with the harmonic supports along the discs, we have the reverse inequality on $V_p\times D$. This proves the comparison used in the ACH application.}

\bibliographystyle{amsplain}
\bibliography{reference}

\end{document}